\newif\ifsiam

\siamfalse

\newif\ifarxiv
\ifsiam \arxivfalse \else \arxivtrue \fi

\ifsiam
\documentclass[review,onefignum,onetabnum]{siamonline220329}
\fi

\ifarxiv
\documentclass[english]{myarticle}
\fi

\usepackage{lipsum}
\usepackage{amsfonts}
\usepackage{amsopn}
\usepackage{graphicx}
\usepackage{epstopdf}
\usepackage{algorithmic}
\usepackage{mathtools,thmtools}
\usepackage{amsmath}
\ifpdf
  \DeclareGraphicsExtensions{.eps,.pdf,.png,.jpg}
\else
  \DeclareGraphicsExtensions{.eps}
\fi

\usepackage{enumitem}
\setlist[enumerate]{leftmargin=.5in}
\setlist[itemize]{leftmargin=.5in}

\ifsiam
\newsiamremark{remark}{Remark}
\newsiamremark{hypothesis}{Hypothesis}
\newsiamthm{claim}{Claim}
\fi

\ifarxiv

\fi

\crefname{hypothesis}{Hypothesis}{Hypotheses}

\ifsiam
\headers{Periodic Normal Forms for Bifurcations of Limit Cycles DDEs}{B. Lentjes, L. Spek, M. M. Bosschaert, Yu. A. Kuznetsov}
\fi

\title{Periodic Normal Forms for Bifurcations of Limit Cycles in DDEs}

\author{Bram Lentjes\thanks{Department of Mathematics, Hasselt University, Diepenbeek Campus, 3590 Diepenbeek, Belgium \email{(bram.lentjes@uhasselt.be)}.}
\and Len Spek\thanks{Department of Applied
Mathematics, University of Twente,
7500 AE Enschede, The Netherlands \email{(l.spek@utwente.nl)}}
\and Maikel M. Bosschaert\thanks{Data Science Institute, Hasselt University, Diepenbeek Campus, 3590 Diepenbeek, Belgium \email{(maikel.bosschaert@uhasselt.be)}.}
\and Yuri A. Kuznetsov\thanks{Department of Mathematics, Utrecht University, 3508 TA Utrecht, The Netherlands and Department of Applied
Mathematics, University of Twente,
7500 AE Enschede, The Netherlands \email{(i.a.kouznetsov@uu.nl)}.}}

\DeclareMathOperator{\NBV}{NBV}
\DeclareMathOperator{\weak}{w}
\DeclareMathOperator{\AC}{AC}
\DeclareMathOperator{\Lip}{Lip}
\DeclareMathOperator{\loc}{loc}
\DeclareMathOperator{\spn}{span}
\DeclareMathOperator{\sign}{sign}

\ifpdf
\hypersetup{
  pdftitle={Periodic Normal Forms for Bifurcations of Limit Cycles in DDEs},
  pdfauthor={B. Lentjes, L. Spek, M. M. Bosschaert, and Yu. A. Kuznetsov}
}
\fi

\begin{document}

\maketitle

\begin{abstract}
A recent work \cite{Lentjes2023} by the authors on the existence of a periodic smooth finite-dimensional center manifold near a nonhyperbolic cycle in delay differential equations motivates the derivation of periodic normal forms. In this paper, we prove the existence of a special coordinate system on the center manifold that will allow us to describe the local dynamics on the center manifold near the cycle in terms of these periodic normal forms. To construct the linear part of this coordinate system, we prove the existence of time periodic smooth Jordan chains for the original and adjoint system. Moreover, we establish duality and spectral relations between both systems by using tools from the theory of delay equations and Volterra integral equations, dual perturbation theory, duality theory and evolution semigroups.
\end{abstract}

\begin{keywords}
delay differential equations, dual perturbation theory, sun-star calculus, center manifold, normal forms, Jordan chains, nonhyperbolic cycles, duality theory, evolution semigroups
\end{keywords}

\begin{MSCcodes}
34C20, 34K13, 34K17, 34K18, 34K19 
\end{MSCcodes}

\begin{sloppypar}

\markboth{B. Lentjes, L. Spek, M.M. Bosschaert, and Yu.A. Kuznetsov}{Periodic Normal Forms for Bifurcations of Limit Cycles in DDEs}

\section{Introduction} \label{sec:introduction}
Bifurcation theory allows one to analyze the behavior of complicated high dimensional nonlinear dynami\-cal systems near bifurcations by reducing the system to a low dimensional invariant manifold, called the center manifold. Using normal form theory, the dynamics on the center manifold can be described by a simple canonical equation called the normal form. These bifurcations and normal forms can be categorized, and their properties can be understood in terms of certain coefficients of the normal form, see \cite{Kuznetsov2023a} for more details. Methods to compute these normal form coefficients have been implemented in software like \verb|MatCont| \cite{Matcont} and \verb|DDE-BifTool| \cite{Engelborghs2002,Sieber2014} to study various classes of dynamical systems.

For bifurcations of limit cycles in continuous-time dynamical systems, there are three generic codimension one bifurcations: fold (or limit point), period-doubling (or flip) and Neimark-Sacker (or torus) bifurcation. These bifurcations are well understood for ordinary differentials equations (ODEs) \cite{Iooss1988, Iooss1999, Kuznetsov2023a,Kuznetsov2005,Witte2014,Witte2013}, but for delay differential equations (DDEs) the theory is still lacking. To understand these bifurcations, one should first prove the existence of a center manifold near a nonhyperbolic cycle. The authors proved in \cite{Lentjes2023} that indeed such a center manifold near a nonhyperbolic cycle exists and is sufficiently smooth. The next step is to study the local dynamics on the center manifold near a nonhyperbolic cycle in terms of normal forms.

The aim of this paper is to show for classical DDEs that the local dynamics on the center manifold near a nonhyperbolic cycle can be studied in terms of periodic normal forms. We generalize the results from Iooss (and Adelmeyer) \cite{Iooss1988,Iooss1999} on periodic normal forms from finite-dimensional ODEs towards infinite-dimensional DDEs. This task will be accomplished by using the rigorous perturbation framework of dual semigroups (sun-star calculus). In an upcoming paper, we present explicit computational formulas for the critical normal form coefficients, along the lines of the periodic normalization method \cite{Kuznetsov2005,Witte2014,Witte2013}, for all codimension one bifurcations of limit cycles in classical DDEs. Finally, we plan to implement the obtained computational formulas into a public software package.

\subsection{Background}
Consider a {\it classical delay differential equation} (DDE) 
\begin{equation} \label{eq:IntroDDE}
    \dot{x}(t) = F(x_t), \quad t \geq 0,
\end{equation}
where $x(t) \in \mathbb{R}^n$ and $x_t \in X := C([-h,0],\mathbb{R}^n)$ represents the \emph{history} of the unknown $x$ at time $t \geq 0$, denoted by $x_t \in X$, i.e.
\begin{equation} \label{eq:history}
    x_t(\theta) := x(t+\theta), \quad \forall\theta \in [-h,0].
\end{equation}
The $\mathbb{R}^n$-valued sufficiently smooth operator $F$ is defined on the Banach space $X$ consisting of $\mathbb{R}^n$-valued continuous functions defined on the compact interval $[-h,0]$, endowed with the supremum norm. Here, $0 < h < \infty$ denotes upper bound of (finite) delays. Furthermore, we assume that \eqref{eq:IntroDDE} has a $T$-periodic solution $\gamma : \mathbb{R} \to \mathbb{R}^n$ such that the periodic orbit (cycle) $\Gamma := \{ \gamma_t \in X : t \in \mathbb{R} \}$ is nonhyperbolic.

Using the perturbation framework of dual semigroups, called sun-star calculus, developed in \cite{Clement1988,Clement1987,Clement1989,Clement1989a,Diekmann1991}, the existence of a periodic smooth finite-dimensional center manifold $\mathcal{W}_{\loc}^c(\Gamma)$ near $\Gamma$ for \eqref{eq:IntroDDE} has been recently rigorously established by the authors in \cite{Lentjes2023}. We also refer to \cite{Bosschaert2020,Church2018,Church2019,Diekmann1991CM,Diekmann1995,Janssens2020,Kuznetsov2023a,Lentjes2023b} for results on the existence of (parameter-dependent) center manifolds near nonhyperbolic equilibria and cycles in finite-dimensional ODEs and classical, impulsive, and abstract infinite-dimensional DDEs. 

To study the local dynamics generated by \eqref{eq:IntroDDE} on $\mathcal{W}_{\loc}^c(\Gamma)$, we will prove the existence of a special coordinate system on the center manifold. Then we show that any solution of \eqref{eq:IntroDDE} on this invariant manifold can be locally parametrized in terms of these coordinates. It turns out that the linear part of our coordinate system will be closely related to the coordinate system invented by Hale and Weedermann \cite{Hale2004}, which was used to study perturbations of periodic orbits in classical DDEs. The existence of such a coordinate system for finite-dimensional ODEs was already established mainly by Iooss in \cite{Iooss1988,Iooss1999}, and so we will generalize his results towards the setting of infinite-dimensional classical DDEs, using the sun-star calculus framework. Iooss already indicated in \cite{Iooss1988} that some of his results should hold in the infinite-dimensional setting. However, as a reader of this paper will see, several results could not be easily extended to the setting of classical DDEs. Moreover, the following issues were even unexpected by the authors:

1) To parametrize $\mathcal{W}_{\loc}^c(\Gamma)$ in terms of this special coordinate system, we need a sufficiently smooth periodic basis of the center eigenspace. Thus, we should be able to construct sufficiently smooth $T$-periodic $X$-valued (generalized) eigenfunctions of a certain linear differential operator $\mathcal{A}$. Hence, the most suitable space to define $\mathcal{A}$ on is a subspace of $C_T(\mathbb{R},X)$ that consists of sufficiently smooth $T$-periodic $X$-valued functions. To give an additional argument on the choice of this space, many functions $\varphi \in C_T(\mathbb{R},X)$ with $X = C([-h,0],\mathbb{R}^n)$ in the sequel will be of the form $\varphi(t)(\theta) = \varphi(t+\theta)(0)$ for all $t \in \mathbb{R}$ and $\theta \in [-h,0]$. Therefore, it is natural to choose twice the space of continuous functions as the time and space variable gets intertwined in the first component of $\varphi$. Consequently, $\mathcal{A}$ emerges as an unbounded linear operator \eqref{eq:D(A)DDE} and to study its (spectral) structure, we would like to have that $\mathcal{A}$ is closed. However, it turns out that $\mathcal{A}$ is a not closed but only closable linear operator in the setting of classical DDEs (\cref{prop:closable}), whereas in the context of finite-dimensional ODEs, this operator is closed (\cref{remark:ODEclosed}). To characterize the operator $\mathcal{A}$ even further, we invoke and lift the concept of evolution semigroups \cite{Chicone1999,Engel2000,Nickel1997,Latushkin1996} towards the setting of $C_T(\mathbb{R},X)$. This provides a concrete description of the closure of $\mathcal{A}$ and allows us to prove that $\mathcal{A}$ is densely defined (\cref{prop:DAdense}). Using this density, we can uniquely define the dual operator $\mathcal{A}^\star$ of $\mathcal{A}$ on a subspace of $C_T(\mathbb{R},X^\star)$ (\cref{prop:pairingcurlyA}). The study on this particular duality will be crucial for our upcoming paper regarding the study of bifurcations of limit cycles in classical DDEs. However, it is important to note already that $C_T(\mathbb{R},X^\star)$ is not a representation of the topological dual space of $C_T(\mathbb{R},X)$ and the pairing $\langle \cdot , \cdot \rangle_T : C_T(\mathbb{R},X^\star) \times C_T(\mathbb{R},X) \to \mathbb{C}$, introduced in \eqref{eq:pairingT} that will uniquely define $\mathcal{A}^\star$, is not the natural dual pairing (\cref{remark:mazur}). Consequently, to establish a direct sum decomposition of $C_T(\mathbb{R},X)$ in terms of (adjoint) (generalized) eigenspaces and annihilators with respect to $\langle \cdot , \cdot \rangle_T$, we need to extend a portion of the Hahn-Banach and closed range theorem (\cref{lemma:pairingnondeg} and \cref{prop:decomposition}) from $X$ towards $C_T(\mathbb{R},X)$, see \eqref{eq:decompCT2} for the final result. Further information on this topic and the (dis)advantages of the choice of the space $C_T(\mathbb{R},X)$ is provided in \cref{sec:conclusions}.

2) Once the existence of the periodic smooth Jordan chains for $\mathcal{A}$ is established, it turns out that we have two bases available for the (generalized) eigenspace. One basis is in terms of history \eqref{eq:history} while the other basis is $T$-periodic. We illustrate in \cref{subsec: time periodic smooth Jordan} that there is an interesting interplay between the history and periodicity of these Jordan chains. It is worth noting that this phenomenon can not occur in finite-dimensional ODEs due to the absence of history.

3) The proof on the existence of this coordinate system on $\mathcal{W}_{\loc}^c(\Gamma)$ happened to be far more involved compared to the finite-dimensional ODE version, see for example \Cref{thm:normalformI} and especially the crucial role of the sun-star calculus machinery in the proof. However,
our findings lead to the conclusion that the (critical) periodic normal forms for bifurcations of limit cycles in classical DDEs are precisely identical to those for ODEs.

Regarding normal form theory for periodic DDEs, we mention the work by Faria \cite{Faria1997} where normal forms were derived for periodic DDEs possessing an autonomous linear part. However, when studying \eqref{eq:IntroDDE} around $\Gamma$, the linear part is given by the non-autonomous term $DF(\gamma_t)$ and so the normal forms derived by Faria do not apply. Therefore, one can see our work as a natural extension of \cite{Faria1997} towards the setting of general periodic DDEs.

\subsection{Overview}
The paper is organized as follows. In \Cref{sec:periodic center manifolds} we review and extend the results from \cite{Lentjes2023} on periodic smooth finite-dimensional center manifolds near nonhyperbolic cycles in the setting of classical DDEs. 

In \Cref{sec:periodic spectral computations} we mainly characterize the center eigenspace and its associated adjoint. To do this, we prove the existence of a time periodic smooth basis of the (adjoint) center eigenspace, see \cref{thm:eigenfunctions} and \cref{thm:adjoint eigenfunctions}. To finalize the linear part, we describe in detail the duality and spectral relations between the linear operators that generate the Jordan chains for the (adjoint) center eigenspace.

In \Cref{sec:characterization} we prove the existence of a special coordinate system on $\mathcal{W}_{\loc}^c(\Gamma)$ that will allow us to study the local dynamics on the periodic center manifold by the means of periodic normal forms, see \Cref{thm:normalformI}, \Cref{thm:normalformII} and \Cref{thm:normalformIII} for the main results.

\section{Periodic center manifolds for classical DDEs} \label{sec:periodic center manifolds}
In this section, we primarily summarize the results from \cite{Lentjes2023} and secondly recall and extend some results from (time-dependent) dual perturbation theory, for which the book \cite{Diekmann1995} together with the article \cite{Clement1988} are standard references. All unreferenced claims related to basic properties of (time-dependent) perturbations of delay equations can be found in these references. For an introduction to the general theory on (adjoint) semigroups, we refer to the well-known books \cite{Engel2000,Neerven1992}.

In the setting of classical DDEs, we work with the real Banach space $X := C([-h,0],\mathbb{R}^n)$ as the state space for some (maximal) finite \emph{delay} $h > 0$ equipped with the supremum norm $\|\cdot\|_{\infty}$. Consider for an integer $k \geq 1$ a $C^{k+1}$-smooth (nonlinear) operator $F : X \to \mathbb{R}^n$ together with the initial value problem
\begin{equation} \label{eq:DDEphi} \tag{DDE} 
    \begin{cases}
    \dot{x}(t) = F(x_t), \quad &t\geq 0,\\
    x_0 = \varphi, \quad &\varphi \in X,
    \end{cases}
\end{equation}
where $x_t \in X$ denotes the {history} of $x$ at time $t \geq 0$, as introduced in (\ref{eq:history}).
By a \emph{solution} of $\eqref{eq:DDEphi}$ we mean a continuous function $x : [-h,t_{\varphi}) \to \mathbb{R}^n$ for some $ t_\varphi \in (0,\infty]$ that is continuously differentiable on $[0,t_{\varphi})$ and satisfies \eqref{eq:DDEphi}. We say that a function $\gamma : \mathbb{R} \to \mathbb{R}^n$ is a \emph{periodic solution} of \eqref{eq:DDEphi} if $\gamma |_{[-h,\infty)}$ is a solution of \eqref{eq:DDEphi} and there exists a minimal $T > 0$, called the \emph{period} of $\gamma$, such that $\gamma_T = \gamma_0$. We call $\Gamma := \{ \gamma_t \in X : t \in \mathbb{R} \}$ a \emph{periodic orbit} or (limit) \emph{cycle}. 

We want to study \eqref{eq:DDEphi} near the periodic solution $\gamma$, and it is therefore convenient to translate $\gamma$ towards the origin. More specifically, if $x$ is a solution of \eqref{eq:DDEphi}, then $y$ defined by $x = \gamma + y$ satisfies the periodic (nonlinear) DDE 
\begin{equation} \label{eq:T-DDEphi1}
\dot{y}(t) = L(t)y_t + G(t,y_t),\\
\end{equation}
where the $C^k$-smooth bounded linear operator $L(t) := DF(\gamma_t) \in \mathcal{L}(X,\mathbb{R}^n)$ denotes the Fr\'echet derivative of $F$ evaluated at $\gamma_t$, and the $C^k$-smooth operator $G(t,\cdot) := F(\gamma_t + \cdot) - F(\gamma_t) - L(t)$ consists solely of nonlinear terms. Regarding the linear part in \eqref{eq:T-DDEphi1}, it is traditional to apply a vector-valued version of the Riesz representation theorem \cite[Theorem 1.1]{Diekmann1995} and write
\begin{equation} \label{eq:L(t)varphi}
    L(t)\varphi = \int_0^h d_2 \zeta(t,\theta) \varphi(-\theta) =: \langle \zeta(t,\cdot),\varphi \rangle, \quad \forall t \in \mathbb{R}, \ \varphi \in X.
\end{equation}
The \emph{kernel} $\zeta : \mathbb{R} \times [0,h] \to \mathbb{R}^{n \times n}$ is a matrix-valued function, $\zeta(t,\cdot)$ is of bounded variation, right-continuous on the open interval $(0,h)$, $T$-periodic in the first component and normalized by the requirement that $\zeta(\cdot,0) = 0$. The integral appearing in \eqref{eq:L(t)varphi} is of Riemann-Stieltjes type and the subscript in $d$ reflects on the fact that we integrate over the second variable of $\zeta$.

The starting point of applying the sun-star calculus construction towards the setting of classical DDEs is by studying the trivial DDE
\begin{equation} \label{eq:TDDEphi}
    \begin{cases}
    \dot{x}(t) = 0, \quad &t\geq 0,\\
    x_0 = \varphi, \quad &\varphi \in X,
    \end{cases}
\end{equation}
which has the unique solution
\begin{equation*}
    x(t) =
    \begin{cases}
    \varphi(t), \quad &-h \leq t \leq 0,\\
    \varphi(0), \quad &t \geq 0.
    \end{cases}
\end{equation*}
We define the $\mathcal{C}_0$-semigroup $T_0 := \{T_0(t)\}_{t \geq 0}$ on $X$, also called the \emph{shift semigroup}, by
\begin{equation*}
    (T_0(t)\varphi)(\theta) := 
    \begin{cases}
    \varphi(t+\theta), \quad &-h \leq t + \theta \leq 0,\\
    \varphi(0), \quad &t+\theta \geq 0,
    \end{cases}
    \quad \forall t \geq 0, \ \varphi \in X, \ \theta \in [-h,0],
\end{equation*}
and notice that $T_0$ generate solutions of \eqref{eq:TDDEphi} as $T_0(t)\varphi = x_t$ for all $t \geq 0$. The shift semigroup has (infinitesimal) generator $A_0$ satisfying
\begin{equation*}
    \mathcal{D}(A_0) = \{ \varphi \in C^{1}([-h,0],\mathbb{R}^n) : \varphi'(0) = 0 \}, \quad A_0 \varphi = \varphi',
\end{equation*}
where the accent denotes the derivative with respect to the state. For this specific combination of $X$, $T_0$ and $A_0$, the abstract duality structure can be constructed explicitly, see \cite[Section II.5]{Diekmann1995} or \cite[Section 2.4]{Bosschaert2020} for a rather compact introduction. We only summarize here the basic results that are needed for the upcoming (sub)sections.

Let us first introduce a convention. For $\mathbb{K} \in \{ \mathbb{R}, \mathbb{C} \}$ let $\mathbb{K}^n$ be the linear space of column vectors, while $\mathbb{K}^{n \star}$ denotes the linear space of row vectors, all with components in $\mathbb{K}$. A representation theorem by Riesz \cite{Riesz1914} enables us to identify the topological dual space $X^\star = C([-h,0],\mathbb{R}^n)^\star$ of $X$ with the Banach space $\NBV([0,h],\mathbb{R}^{n\star})$ consisting of functions $\zeta : [0,h] \to \mathbb{R}^{n \star}$ that are normalized by $\zeta(0) = 0$, are continuous from the right on $(0,h)$ and have bounded variation. The (canonical) natural dual pairing between $X^\star$ and $X$ will be denoted by $\langle \cdot, \cdot \rangle : X^\star \times X \to \mathbb{R}$ and note that this is a bilinear map. Recall that we already used this notation implicitly in \eqref{eq:L(t)varphi}, but in a vector-valued version fashion. The same $\langle \cdot, \cdot \rangle$-notation will be used in the sequel to denote the natural dual pairing between a Banach space $E$ and its dual $E^\star$. 

Because $X$ is not reflexive, the \emph{dual semigroup} $T_0^\star$ of $T_0$ is only weak$^\star$ continuous on $X^\star$. This is also visible on the generator level, as the adjoint $A_0^\star$ of $A_0$ is only the weak$^\star$ generator of $T_0^\star$ and takes the form
\begin{align} \label{eq:D(Astar0)}
    \mathcal{D}(A_0^\star) = \bigg\{  f \in \NBV([0,h],\mathbb{R}^{n \star}) : \ & f(\theta) =  f(0^+) + \int_0^\theta g(\sigma) d\sigma \mbox{ for all }  \theta \in (0,h], \nonumber \\
     & \ g \in \NBV([0,h],\mathbb{R}^{n\star}) \mbox{ and }  g(h) = 0 \bigg \}, \quad A_0^\star f = g,
\end{align}
where $f(0^+) := \lim_{t \downarrow 0} f(t)$ and the function $g =: f'$ is called the \emph{derivative} of $f$. The maximal subspace of strong continuity
\begin{equation*}
    X^\odot := \{ x^\star \in X^\star : t \mapsto T_0^\star(t) x^\star \mbox{ is norm continuous on } [0,\infty) \}
\end{equation*}
is a norm closed $T_0^\star(t)$-invariant weak$^\star$ dense linear subspace of $X^\star$ and we have the characterization
\begin{equation} \label{eq:Xsuncharac}
    X^\odot = \overline{\mathcal{D}(A_0^\star)},
\end{equation}
where the bar denotes the norm closure in $X^\star$. Expression \eqref{eq:Xsuncharac} enables us to compute the \emph{sun dual} $X^\odot$ by taking the closure of $\mathcal{D}(A_0^\star)$ with respect to the norm defined on $\NBV([0,h],\mathbb{R}^{n\star})$. As the space of bounded variation functions is norm dense in $L^1$, $X^\odot$ has the same description as $\mathcal{D}(A_0^\star)$, but the derivative $g$ is allowed to be in $L^1([0,h],\mathbb{R}^{n\star})$. 

Let $\AC_0([0,h],\mathbb{R}^{n \star})$ denote the space of $\mathbb{R}^{n \star}$-valued functions that are absolute continuous on $(0,h]$, have zero value at zero and zero derivative at $h$. From the description of the sun dual, it is clear that $X^\odot = \AC_0([0,h],\mathbb{R}^{n \star})$. It turns out that another characterization of the sun dual will be helpful in the sequel. As an element $f \in X^\odot$ is completely specified by $f(0^+) \in \mathbb{R}^{n \star}$ and $g \in L^1([0,h],\mathbb{R}^{n\star})$, we obtain the isometric isomorphism
\begin{equation*}
    X^{\odot} \cong \mathbb{R}^{n \star} \times L^1([0,h],\mathbb{R}^{n\star}),
\end{equation*}
and note that $\mathcal{D}(A_0^\star)$ can be identified via this isomorphism with the subspace $\mathbb{R}^{n \star} \times \NBV([0,h],\mathbb{R}^{n\star})$ of $X^\odot$. The restriction $T_0^\odot(t) := T_0^\star(t) |_{X^\odot}$ forms a $\mathcal{C}_0$-semigroup on $X^\odot$ and its generator $A_0^\odot$ is the part of $A_0^\star$ in $X^\odot$. The (restricted) natural dual pairing $\langle \cdot , \cdot \rangle : X^\odot \times X \to \mathbb{R}$ between $\varphi^\odot = (c,g) \in X^\odot$ and $\varphi \in X$ is given by
\begin{equation} \label{eq:Xsunpairing} 
    \langle \varphi^\odot, \varphi \rangle = c \varphi(0) + \int_0^h g(\theta) \varphi(-\theta) d\theta.
\end{equation}
We have at this moment a $\mathcal{C}_0$-semigroup $T_0^\odot$ with generator $A_0^\odot$ on the Banach space $X^\odot$, which are precisely the ingredients we started with. Hence, we can repeat the construction once more. A representation of the dual $X^{\odot \star}$ of $X^\odot$, and its restriction to the maximal subspace of strong continuity $X^{\odot \odot}$, is given by
\begin{equation*}
    X^{\odot \star} \cong \mathbb{R}^n \times L^\infty([-h,0],\mathbb{R}^n), \quad X^{\odot \odot} \cong \mathbb{R}^n \times C([-h,0],\mathbb{R}^n),
\end{equation*}
and the weak$^\star$ continuous generator of $T_0^{\odot \star}$ takes the form
\begin{equation*}
    \mathcal{D}(A_0^{\odot \star}) = \{ (\alpha,\varphi) \in X^{\odot \star} : \varphi \in \Lip([-h,0],\mathbb{R}^n) \mbox{ and } \varphi(0) = \alpha \}, \quad A_0^{\odot \star}(\alpha,\varphi) = (0,\varphi').
\end{equation*}
The natural dual pairing $\langle \cdot, \cdot \rangle : X^{\odot \star} \times X^\odot \to \mathbb{R}$ between $\varphi^{\odot \star} = (a,\psi) \in X^{\odot \star}$ and $\varphi^\odot = (c,g) \in X^{\odot}$ is given by
\begin{equation} \label{eq:Xsunstarpairing}
    \langle \varphi^{\odot \star}, \varphi^\odot \rangle = ca + \int_0^h g(\theta) \psi(-\theta) d\theta.
\end{equation}
The linear (canonical) embedding $j: X \to X^{\odot \star}$ has action $j\varphi = (\varphi(0),\varphi)$ for all $\varphi \in X$, mapping $X$ onto $X^{\odot \odot}$, meaning that $X$ is $\odot$-\emph{reflexive} with respect to the shift semigroup $T_0$. 

Next, we turn our attention to the periodic linear DDE
\begin{equation} \label{eq:T-LDDEphi}
    \begin{cases}
        \dot{y}(t) =  L(t)y_t, \quad &t \geq s, \\
        y_s = \varphi, \quad &\varphi \in X,
    \end{cases}
\end{equation}
in the setting of (time-dependent) dual perturbation theory. Recall that $L(t) = DF(\gamma_t)$ for all $t \in \mathbb{R}$ and $s \in \mathbb{R}$ denotes the \emph{starting time}. For $i = 1,\dots,n$ let $r_i^{\odot \star} := (e_i,0)$, where $e_i$ is the $i$th standard basic vector of $\mathbb{R}^n$. It is conventional and convenient to introduce the shorthand notation
\begin{equation*}
    wr^{\odot \star} := \sum_{i = 1}^{n} w_ir_i^{\odot \star}, \quad \forall w = (w_1,\dots,w_m) \in \mathbb{R}^{n},
\end{equation*}
and note that $wr^{\odot \star} = (w,0) \in X^{\odot \star}$. To lift the linear perturbation $L(t)$ from $X$ to $X^{\odot \star}$, let us introduce the finite rank $T$-periodic \emph{time-dependent bounded linear perturbation} $B$ by
\begin{equation} \label{eq:perturbationB}
    B(t)\varphi := [L(t)\varphi]r^{\odot \star}, \quad \forall t \in \mathbb{R}, \ \varphi \in X,
\end{equation}
and since $F \in C^{k+1}(X,\mathbb{R}^n)$, $L(t) = DF(\gamma_t)$ and $t \mapsto \gamma_t$ is $T$-periodic and of the class $C^k$, we have that $B \in C^{k}( \mathbb{R}, \mathcal{L}(X,X^{\odot \star}))$ is $T$-periodic and Lipschitz continuous. It is proven in \cite[Theorem XII.
3.1]{Diekmann1995} that there is a one-to-one correspondence between solutions of \eqref{eq:T-LDDEphi} and the time-dependent linear abstract integral equation
\begin{equation} \label{eq:T-LAIEphi}
    u(t) = T_0(t-s)\varphi + j^{-1}\int_s^t T_0^{\odot \star}(s-\tau) B(\tau)u(\tau) d\tau, \quad \varphi \in X,
\end{equation}
with $t \geq s$. Here, the integral has to be interpreted as a weak$^\star$ Riemann integral \cite[Chapter III]{Diekmann1995} and takes values in $j(X)$ under the running assumption of $\odot$-reflexivity, see \cite[Lemma 2.2]{Clement1988}. The unique solution of \eqref{eq:T-LAIEphi} on $[s,\infty)$ is generated by a strongly continuous forward evolutionary system  $U := \{U(t,s)\}_{(t,s) \in \Omega_\mathbb{R}}$ on $X$ in the sense that $u(t) = U(t,s)\varphi$ for all $t \in [s,\infty)$, where $\Omega_J := \{(t,s) \in J \times J : t \geq s \}$ for some interval $J \subseteq \mathbb{R}$, see \cite[Definition XII.2.1 and Theorem XII.2.7]{Diekmann1995}. Because $B$ is $T$-periodic, there holds
\begin{equation} \label{eq:Uperiodic}
    U(t+T,s+T) = U(t,s), \quad U(s+lT,s) = U(s+T,s)^l,
\end{equation}
for all $(t,s) \in \Omega_{\mathbb{R}}$ and $l \in \mathbb{N}$, see \cite[Corollary XIII.2.2]{Diekmann1995}. As we have defined $U(t,s)$ for all $(t,s) \in \Omega_\mathbb{R}$, we are interested in the associated (sun) dual(s). It is clear that one can define $U^\star(s,t) := U(t,s)^\star \in \mathcal{L}(X^\star):=\mathcal{L}(X^\star,X^\star)$ and that $U^\star := \{U^\star(s,t)\}_{(s,t) \in \Omega_\mathbb{R}^\star}$ forms a backward evolutionary system on $X^\star$, where $\Omega_J^\star := \{(s,t) \in J^2 : t \geq s \}$ for some interval $J \subseteq \mathbb{R}$. Furthermore, the Lipschitz continuity of $B$ ensures that the restriction $U^{\odot}(s,t) := U^\star(s,t) |_{X^\odot}$ leaves $X^\odot$ invariant and thus, by construction, $U^\odot := \{U^{\odot}(s,t)\}_{(s,t) \in \Omega_\mathbb{R}^\star}$ is a strongly continuous backward evolutionary system, see \cite[Theorem 5.3]{Clement1988}. This allows us to define $U^{\odot \star}(t,s) := (U^{\odot}(s,t))^\star$ and it is clear that $U^{\odot \star} := \{U^{\odot \star}(t,s) \}_{(t,s) \in \Omega_\mathbb{R}}$ is a forward evolutionary system on $X^{\odot \star}$ that extends $U$, which was previously defined on $X$.

Let us now characterize for all $\tau\in \mathbb{R}$ the (generalized) generators $A(\tau), A^{\star}(\tau), A^\odot(\tau)$ and $A^{\odot \star}(\tau)$ of the evolutionary systems $U,U^\star,U^\odot$ and $U^{\odot \star}$, respectively. More information on these generators can be found in \cite{Clement1988,Lentjes2023}. It turns out that the closed and densely defined unbounded linear operator $A(\tau)$ has the form
\begin{equation} \label{eq:D(As)}
    \mathcal{D}(A(\tau)) = \{ \varphi \in C^1([-h,0],\mathbb{R}^n) : \varphi'(0) = L(\tau) \varphi \}, \quad A(\tau)\varphi = \varphi'.
\end{equation}
The weak$^\star$ continuous \emph{dual generator} $A^\star(\tau) := [A(\tau)]^\star$ of $A(\tau)$ is then a well-defined closed linear operator with weak$^\star$ dense domain $\mathcal{D}(A^\star(\tau))$. A concrete description of this operator is treated in the following result.
\begin{lemma} \label{lemma:Astar}
The dual generator $A^{\star}(\tau)$ has the representation
\begin{align*}
    \mathcal{D}(A^{\star}(\tau)) &= \mathcal{D}(A_0^{\star}), \quad A^{\star}(\tau)f = f' + f(0^{+})\zeta(\tau,\cdot).
\end{align*}
\end{lemma}
\begin{proof}
The equality between the domains follows from a sun-variant of \cite[Lemma 4.2]{Clement1988}. To show the action, we first have to determine $B^{\star}(\tau) := [B(\tau)]^\star$ and notice that we traditionally restrict this map to $X^{\odot}$, see \cite[page 58]{Diekmann1995}. It follows from the uniqueness of the adjoint that $B^{\star}(\tau) : X^\odot \to X^\star$ is given by $B^\star(\tau) (c,g) = c\zeta(\tau,\cdot)$ since
\begin{equation*}
\langle B^{\star}(\tau) (c,g), \varphi \rangle = \langle c \zeta(\tau,\cdot),\varphi \rangle  = \langle (\langle \zeta(\tau,\cdot),\varphi \rangle ,0),(c,g) \rangle = \langle B(\tau)\varphi, (c,g) \rangle
\end{equation*}
for all $\tau \in \mathbb{R}, \varphi \in X$ and $(c,g) \in X^\odot$. Here, we used \eqref{eq:Xsunstarpairing} in the second equality, and the action of $B$ in the third. Hence, $A^{\star}(\tau)f = A_0^{\star}f + B^\star(\tau)f = f' + f(0^{+})\zeta(\tau,\cdot)$, which completes the proof.
\end{proof}
This allows us to introduce $A^{\odot}(\tau) := [A(\tau)]^\odot$ which is the part of $A^{\star}(\tau)$ in $X^\odot$. This results into a closed and densely defined unbounded linear operator that takes the form
\begin{equation} \label{eq:D(Asuns)}
    \mathcal{D}(A^\odot(\tau)) = \{(c,g) \in \mathcal{D}(A_0^\star) : g+c\zeta(\tau,\cdot) \in X^{\odot} \}, \quad A^\odot(\tau)(c,g) = g + c\zeta(\tau,\cdot),
\end{equation}
The weak$^\star$ continuous generator of $U^{\odot \star}$ is given by $A^{\odot \star}(\tau) := [A^\odot(\tau)]^\star$ and is a well-defined closed linear operator with weak$^\star$ dense domain $\mathcal{D}(A^{\odot \star}(\tau))$. The representation of $A^{\odot \star}(\tau)$ for periodic DDEs can be computed as in \cref{lemma:Astar} and is given by
\begin{align} \label{eq:Asunstartau}
    \mathcal{D}(A^{\odot \star}(\tau)) = \mathcal{D}(A_0^{\odot \star}), \quad A^{\odot \star}(\tau)(\alpha,\varphi) = (L(\tau)\varphi,\varphi').
\end{align}
To go full circle, note that $A(\tau)$ can also be obtained from $A^{\odot \star}(\tau)$ as the preimage under $j$ of the part of $A^{\odot \star}(\tau)$ in $j(X)$. Since $\tau \mapsto \zeta(\tau,\cdot)$ is of the class $C^{k}$, it is clear from the definitions of the generators that the maps $\tau \mapsto A(\tau) , \tau \mapsto A^{\star}(\tau), \tau \mapsto A^{\odot}(\tau)$ and $\tau \mapsto A^{\odot \star}(\tau)$ are all of the class $C^{k}$. 

Our next aim is to study how the periodic (nonlinear) delay differential equation
\begin{equation} \label{eq:T-DDEphi} 
    \begin{cases}
        \dot{y}(t) =  L(t)y_t + G(t,y_t), \quad &t \geq s, \\
        y_s = \varphi, \quad &\varphi \in X,
    \end{cases}
\end{equation}
fits naturally in the setting of (time-dependent) dual perturbation theory. Recall that $G(t,\cdot) = F(\gamma_t + \cdot) - F(\gamma_t) - L(t)$ and so $G$ is $C^k$-smooth and $T$-periodic in the first component. The $T$-periodic $C^k$-smooth \emph{time-dependent nonlinear perturbation} $R$ is given by
\begin{equation*}
    R(t,\varphi) := G(t,\varphi)r^{\odot \star}, \quad \forall (t,\varphi) \in \mathbb{R} \times X,
\end{equation*}
and it is shown in \cite[Theorem 7]{Lentjes2023} that there is a one-to-one correspondence between solutions of \eqref{eq:T-DDEphi} and the time-dependent nonlinear abstract integral equation
\begin{equation} \label{eq:T-AIE}
    u(t) = U(t,s)\varphi + j^{-1}\int_s^t U^{\odot \star}(t,\tau) R(\tau,u(\tau)) d\tau, \quad \varphi \in X,
\end{equation}
for $t \geq s$. Here, the integral has to be interpreted as a weak$^\star$ Riemann integral and takes values in $j(X)$ under the running assumption of $\odot$-reflexivity, see \cite[Lemma 1]{Lentjes2023}. The $C^k$-smoothness of the nonlinearity $R$ ensures that for any $\varphi \in X$ there exists a unique (maximal) solution $u_\varphi$ of \eqref{eq:T-AIE} on some (maximal) interval $[s,t_\varphi)$ of $\mathbb{R}$ with $s < t_\varphi \leq \infty$, see \cite[Proposition 1]{Lentjes2023}.

Let us now turn our attention towards the construction of the periodic center manifold $\mathcal{W}_{\loc}^c(\Gamma)$ around the cycle $\Gamma$ from \cite[Section 3.6]{Lentjes2023}. We start with characterizing its linear part. The spectrum $\sigma(U(s+T,s))$ of the \emph{monodromy operator} (at time $s$) $U(s+T,s) \in \mathcal{L}(X)$ is a countable set in $\mathbb{C}$, independent of the starting time $s$, consisting of $0$ and isolated eigenvalues of finite type, called \emph{Floquet multipliers}, that can possibly accumulate to $0$. It is known that $1$ is always a Floquet multiplier, called the \emph{trivial Floquet multiplier}, with associated eigenfunction $\dot{\gamma}_s$. The number $\sigma \in \mathbb{C}$ satisfying $\lambda = e^{\sigma T}$ is called the \emph{Floquet exponent} (of $\lambda$) and note that this number is uniquely determined up to an additive factor of $\frac{2 \pi i }{T}$. The \emph{(generalized) eigenspace} of $U(s+T,s)$ (at time $s$) associated to the Floquet multiplier $\lambda$ is defined by $E_\lambda(s) := \mathcal{N}((\lambda I - U(s+T,s))^{k_\lambda})$, where $1 \leq k_\lambda < \infty$ is the order of the pole $\lambda = \mu$ of the map $\mu \mapsto (\mu I - U(s+T,s))^{-1}$. Hence, $E_\lambda(s)$ is finite-dimensional, and its dimension, called the \emph{algebraic multiplicity} (of $\lambda$), will be denoted by $m_\lambda$. The \emph{geometric multiplicity} (of $\lambda$) reflects on the dimension of the associated eigenspace $\mathcal{N}(\lambda I - U(s+T,s))$. As a consequence, the set of Floquet multipliers on the unit circle $\Lambda_0 := \{ \lambda \in \sigma(U(s+T,s)) : |\lambda| = 1\}$ is finite, say $1 \leq n_0 + 1 < \infty$ when counted with algebraic multiplicity, and we define the $(n_0 +1)$-dimensional \emph{center eigenspace} (at time $s$) by $X_0(s) := \oplus_{\lambda \in \Lambda_0} E_{\lambda}(s)$. In this setting, the periodic local center manifold theorem \cite[Corollary 2]{Lentjes2023} for \eqref{eq:DDEphi} applies. To be more precise, define the \emph{center fiber bundle} as $X_0 := \{ (t,\varphi) \in \mathbb{R} \times X : \varphi \in X_0(t) \}$ and denote for any $\delta > 0$ the $\delta$-ball in $X$ centered at the origin by $B_\delta(X)$. Then there exists a $C^k$-smooth map $\mathcal{C} : X_0 \to X$ and a sufficiently small $\delta > 0$ such that the manifold 
\begin{equation} \label{eq:Wloc}
    \mathcal{W}_{\loc}^{c}(\Gamma) := \{ \gamma_t + \mathcal{C}(t, \varphi) \in X : (t,\varphi) \in X_0 \mbox{ and } \mathcal{C}(t,\varphi) \in B_\delta(X) \}
\end{equation}
is a $T$-periodic $C^k$-smooth $(n_0+1)$-dimensional locally positively invariant manifold in $X$, called the \emph{(local) center manifold around $\Gamma$}, defined in the vicinity of $\Gamma$ for a sufficiently small $\delta > 0$.

\section{Periodic spectral computations for classical DDEs} \label{sec:periodic spectral computations}
To construct in \cref{sec:characterization} a periodic smooth parametrization of $\mathcal{W}_{\loc}^{c}(\Gamma)$, it is clear that we need a periodic smooth basis of the center eigenspace $X_0(s)$. Therefore, we establish in \cref{subsec: time periodic smooth Jordan} that any (generalized) eigenspace $E_\lambda(s)$ possesses a periodic smooth basis, see \cref{thm:eigenfunctions} for the main result. To obtain such a basis, we introduce in \eqref{eq:D(A)DDE} the linear differential operators $\mathcal{A}$ and $\mathcal{A}^{\odot \star}$ defined on a subspace of $C_T(\mathbb{R},X)$ and $C_T(\mathbb{R},X^{\odot \star})$ respectively, and we prove that these operators admit a sufficiently smooth periodic Jordan chain. As the spectral structure of $\mathcal{A}$ and $\mathcal{A}^{\odot \star}$ will be important in our upcoming paper on bifurcations of limit cycles in classical DDEs, we would like to prove that both operators are closed. However, we find, to our surprise, that both operators are not closed but only closable, as detailed in \cref{prop:closable}. To gain deeper insights into the underlying structure of $\mathcal{A}$ and $\mathcal{A}^{\odot \star}$, we introduce the concept of evolution semigroups on the space $C_T(\mathbb{R},X)$ and $C_T(\mathbb{R},X^{\odot \star})$, respectively. This enables us to further characterize, for example, the domain $\mathcal{D}(\mathcal{A})$ of the operator $\mathcal{A}$. In particular, we establish its density in $C_T(\mathbb{R},X)$, see \cref{prop:DAdense}. In \cref{remark:ODEclosed}, we elaborate on the similarities and differences of the operator $\mathcal{A}$ in finite-dimensional ODEs and classical DDEs.

In our forthcoming paper on bifurcations of limit cycles, we also need a periodic smooth basis for the adjoint center eigenspace $X_0^\star(s)$. This will be discussed in \cref{subsec: dual time periodic smooth Jordan}, where we consider the unbounded linear operators $\mathcal{A}^\odot$ and $\mathcal{A}^\star$, introduced in \eqref{eq:D(Asun)DDE}, defined on subspaces of $C_T(\mathbb{R},X^\odot)$ and $C_T(\mathbb{R},X^{\star})$, respectively. While the formulation of our results in this subsection will resemble those in \cref{subsec: time periodic smooth Jordan}, some proofs will be considerably more involved.

The natural question arises how the non-closed but closable unbounded linear operators $\mathcal{A}, \mathcal{A}^\star,\mathcal{A}^\odot$ and $\mathcal{A}^{\odot \star}$ are interconnected through duality theory. This topic will be addressed in \cref{subsec: duality relations}. We demonstrate in \cref{prop:pairingcurlyA} that $\mathcal{A}^\star$ is the adjoint of $\mathcal{A}$ and $\mathcal{A}^{\odot \star}$ is the adjoint of $\mathcal{A}^{\odot}$ with respect to the pairings $\langle \cdot , \cdot \rangle_T$ defined in \eqref{eq:pairingT}. However, it is important to note that this pairing does not represent the natural dual pairing since $C_T(\mathbb{R},X^\star)$ and $C_T(\mathbb{R},X^{\odot \star})$ are not representations of the topological dual spaces of $C_T(\mathbb{R},X)$ and $C_T(\mathbb{R},X^\odot)$, respectively. Further details regarding these duality relations are provided in \cref{remark:mazur}. To finalize this subsection, we present in \cref{prop:decomposition} particular direct sum decompositions of the spaces $C_T(\mathbb{R},X)$ and $C_T(\mathbb{R},X^{\odot \star})$ which will play a significant role in our upcoming paper on bifurcations of limit cycles in classical DDEs.

Before we continue, let us recall from \cref{sec:periodic center manifolds} that the semigroups and generators on $X^{\star}$ and $X^{\odot \star}$ are only defined in a weak$^\star$ sense. Therefore, it is no surprise that we have to study abstract ODEs in a weak$^\star$ setting. Hence, let us recall the definition of the (partial) weak$^\star$ differential operators.

\begin{definition}[{\cite[Definition 4.4]{Clement1988}}] \label{def:wkstarderivative}
Let $E$ be a Banach space, $J \subseteq \mathbb{R}$ an interval and $\Omega \subseteq J \times J$. We say that a function $f : J \to E^\star$ is \emph{weak$^\star$ differentiable} with \emph{weak$^\star$ derivative} $d^\star f : J \to E^\star$ if
\begin{equation*}
    \frac{d}{dt} \langle f(t),x \rangle = \langle d^\star f (t), x \rangle, \quad \forall x \in E, \ t \in J.
\end{equation*}
If in addition $d^\star f$ is weak$^\star$ continuous, then $f$ is called \emph{weak$^\star$ continuously differentiable}. Furthermore, we say that a function $g : \Omega \to E^\star$ has \emph{partial weak$^\star$ derivatives} $\partial_t^\star g : \Omega \to E^\star$ and $\partial_s^\star g : \Omega \to E^\star$ if
\begin{align*}
    \frac{\partial}{\partial t} \langle g(t,s),x \rangle &= \langle \partial_t^\star g (t,s), x \rangle, \quad \forall x \in E, \ (t,s) \in \Omega, \\
    \frac{\partial}{\partial s} \langle g(t,s),x \rangle &= \langle \partial_s^\star g (t,s), x \rangle, \quad \forall x \in E, \ (t,s) \in \Omega.
\end{align*}
If in addition $\partial_t^\star g$ and $\partial_s^\star g$ are weak$^\star$ continuous, then $g$ is called \emph{weak$^\star$ continuously differentiable}.
\end{definition}
\begin{remark} \label{remark:complexification}
For the spectral theory on the real Banach space $X$, we have to \emph{complexify} $X$ and all discussed operators on $X$. This is not entirely trivial and is discussed in \cite[Section III.7 and Section IV.2]{Diekmann1995}. To clarify, by the spectrum of a real (unbounded) operator $L$ defined on (a subspace of) $X$, we mean the spectrum of its complexification $L_\mathbb{C}$ on (a subspace of) the complexified Banach space $X_\mathbb{C}$. For the ease of notation, we omit the additional symbols. \hfill $\lozenge$
\end{remark}

\subsection{Periodic smooth Jordan chains} \label{subsec: time periodic smooth Jordan}
Before we construct and explore the periodic smooth Jordan chains, let us first prove the following result as it will be helpful in the sequel. Its autonomous analogue, written in the language of $\mathcal{C}_0$-semigroups, can be found in \cite[Corollary III.4.8]{Diekmann1995}.
\begin{lemma} \label{lemma:Dinvariant}
For any $s \in \mathbb{R}$, the operator $U(\tau,s)$ maps $X$ into $\mathcal{D}(A(\tau))$ for all $\tau \geq s+h$. Moreover, there holds
\begin{equation*}
    \frac{\partial}{\partial \tau} U(\tau,s)\varphi = A(\tau)U(\tau,s)\varphi, \quad \forall \tau \geq s+h, \ \varphi \in X,
\end{equation*}
where the partial derivative should be interpreted in the norm topology of $X$.
\end{lemma}
\begin{proof}
Let $\varphi \in X$ and $s \in \mathbb{R}$ be given. Recall from \cref{sec:periodic center manifolds} that $U(\tau,s)\varphi = y_\tau$ for $\tau \geq s$, where $y$ is the unique solution of \eqref{eq:T-LDDEphi}. Applying the methods of steps \cite[Theorem 1.2.2]{Hale1993} onto \eqref{eq:T-LDDEphi} shows that $\tau \mapsto y(\tau)$ is continuously differentiable on $[s,\infty)$, i.e. $\tau \mapsto y_\tau = U(\tau,s)\varphi \in C^1([-h,0],\mathbb{C}^n)$ for all $\tau \geq s + h$. It remains to show that $y_\tau$ satisfies the second condition of the domain defined in \eqref{eq:D(As)}. We obtain from \eqref{eq:history} that 
\begin{equation} \label{eq:DDEtranslation}
    \frac{d}{d \theta} y_\tau(\theta) \bigg |_{\theta = 0} = \frac{d}{d \theta} y(\tau + \theta) \bigg |_{\theta = 0} = \dot{y}(\tau) = L(\tau)y_t, \quad \forall \tau \geq s,
\end{equation}
which proves that $U(\tau,s)\varphi \in \mathcal{D}(A(\tau))$ for all $\tau \geq s+h$.

To prove the second claim, consider in addition $\psi = U(\tau,s)\varphi \in \mathcal{D}(A(\tau))$ for some $\tau \geq s+h$. Then,
\begin{align*}
    \bigg\| \frac{1}{h}(U(\tau+h,s)\varphi - U(\tau,s)\varphi) - A(\tau)U(\tau,s)\varphi \bigg\| &= \bigg\| \frac{1}{h}(U(\tau+h,\tau)U(\tau,s)\varphi - U(\tau,s)\varphi) - A(\tau)U(\tau,s)\varphi \bigg\| \\
    &= \bigg\| \frac{1}{h}(U(\tau+h,t)\psi - \psi) - A(\tau)\psi \bigg\| 
    \to 0, \quad h \downarrow 0,
\end{align*}
due to \cite[Lemma 2.4]{Clement1988} as $\psi \in \mathcal{D}(A(\tau))$.
\end{proof}
Let us now focus on a specific Floquet multiplier $\lambda \in \sigma(U(s+T,s))$ for a fixed $s \in \mathbb{R}$. By the construction given in \cite[Section IV.4]{Diekmann1995}, it is possible to find a basis of $E_\lambda(s)$ that is in \emph{Jordan normal form}. That is, there exists an ordered basis $\{\phi_s^0,\dots,\phi_s^{m_\lambda - 1} \}$ of $E_\lambda(s)$, called a \emph{Jordan chain} (of $E_\lambda(s)$), satisfying
\begin{equation} \label{eq:generalizedeigenvecU(s+T,s)}
    (U(s+T,s) - \lambda I)\phi_s^i =
    \begin{cases}
    0, \quad &i=0,\\
    \phi_s^{i-1}, \quad &i=1,\dots,m_\lambda-1,
    \end{cases}
\end{equation}
where $\phi_s^i$ should be interpreted as history \eqref{eq:history}. As the map $U_\lambda(\tau,s) := U(\tau,s)|_{E_\lambda(s)} : E_\lambda(s) \to E_\lambda(\tau)$ is a topological isomorphism \cite[Theorem XIII.3.3]{Diekmann1995}, we know that $\{\phi_\tau^0,\dots,\phi_\tau^{m_\lambda - 1} \}$ is an ordered basis of $E_\lambda(\tau)$, where $\phi_\tau^i := U_\lambda(\tau,s)\phi_s^i$ for all $\tau \in \mathbb{R}$ and $i=0,\dots,m_\lambda - 1$. The following lemma shows that this specific basis of $E_\lambda(\tau)$ has additional structure.

\begin{lemma} \label{lemma:basisD(A)}
The ordered basis $\{\phi_\tau^0,\dots,\phi_\tau^{m_\lambda -1} \} \subseteq \mathcal{D}(A(\tau))$ consists of $C^{k+1}$-smooth functions and forms a Jordan chain for $E_\lambda(\tau)$ for all $\tau \in \mathbb{R}$.
\end{lemma}
\begin{proof}
Since $\phi_\tau^i = U_\lambda(\tau,s)\phi_s^i$ for all $\tau \in \mathbb{R}$, it is clear from the computation
\begin{equation*}
    (U(\tau+T,\tau) -\lambda I)\phi_\tau^i = U(\tau,s)(U(s+T,s) - \lambda I) \phi_s^i = 
    \begin{dcases}
            0, \quad &i = 0, \\
            \phi_\tau^{i-1}, \quad &i=1,\dots,m_\lambda-1,
    \end{dcases}
\end{equation*}
that $\{\phi_\tau^0,\dots,\phi_\tau^{m_\lambda -1} \}$ is a Jordan chain for $E_\lambda(\tau)$. Let us now prove by induction that $\phi_\tau^i \in \mathcal{D}(A(\tau))$ for a fixed $\tau \in \mathbb{R}$. If $i = 0$, choose $m \in \mathbb{N}$ large enough to guarantee that $\tau + mT \geq s + h$ because then $\phi_{\tau + mT}^0 \in \mathcal{D}(A(\tau))$ as a consequence of \cref{lemma:Dinvariant} and recalling that $\tau \mapsto A(\tau)$ is $T$-periodic, see \eqref{eq:D(As)}. It follows from 
\begin{equation*}
   \lambda^m \phi_\tau^0 = U(\tau+T,\tau)^m\phi_\tau^0 = U(\tau+mT,\tau)\phi_\tau^0 = \phi_{\tau + mT}^0,
\end{equation*}
where we used \eqref{eq:Uperiodic} in the second equality, that $\phi_\tau^0 \in \mathcal{D}(A(\tau))$ due to linearity and the fact that the Floquet multiplier $\lambda \neq 0$. Now, assume that $\phi_\tau^0,\dots,\phi_\tau^{i-1}$ are in $\mathcal{D}(A(\tau))$ for some $i \in \{1,\dots,m_\lambda - 1\}$. Choose again a $m \in \mathbb{N}$ such that $\tau + mT \geq s + h$. By the Jordan chain structure \eqref{eq:generalizedeigenvecU(s+T,s)}, there exist $c_{m,0},\dots,c_{m,i} \in \mathbb{C}$ such that
\begin{equation*}
\lambda^m \phi_\tau^{i} + \sum_{l=0}^{i-1}c_{m,l} \phi_\tau^l = U(\tau+mT,\tau)\phi_\tau^i = \phi_{\tau + mT}^i \in \mathcal{D}(A(\tau)).
\end{equation*}
By the induction hypothesis, $\phi_\tau^0,\dots,\phi_\tau^{i-1}$ are all in $\mathcal{D}(A(\tau))$ and so we conclude that $\lambda^m \phi_\tau^{i} \in \mathcal{D}(A(\tau))$ which proves that $\phi_\tau^{i} \in \mathcal{D}(A(\tau))$ since the Floquet multiplier $\lambda \neq 0$.

The same arguments can be used now, by choosing $m \in \mathbb{N}$ large enough to guarantee that $\tau + mT \geq s + (k+1)h$ and employing the method of steps, as performed in the proof \cref{lemma:Dinvariant}, to conclude that all maps $\tau \mapsto \phi_\tau^i$ are all $C^{k+1}$-smooth.
\end{proof}
Let us now take a look at the $T$-(anti)periodicity of the Jordan chain defined in \eqref{eq:generalizedeigenvecU(s+T,s)}. Therefore, let us first recall the notion of $T$-(anti)periodicity for vector-valued functions. For a (complex) Banach space $E$, a function $\varphi : \mathbb{R} \to E$ is called $T$\emph{-periodic} if $\varphi(\tau) = \varphi(\tau + T)$ for all $\tau \in \mathbb{R}$, and $T$\emph{-antiperiodic} if $\varphi(\tau) = -\varphi(\tau + T)$ for all $\tau \in \mathbb{R}$. Notice that $T$-antiperiodic maps are $2T$-periodic. In this $T$-antiperiodic setting, it is convenient to define for a Floquet multiplier $\lambda$ its associated Floquet exponent $\sigma \in \mathbb{C}$ by $e^{\sigma T} = |\lambda|$. Let us now turn our attention to the Jordan chains from \eqref{eq:generalizedeigenvecU(s+T,s)}. It is clear from the computation
\begin{equation*}
    \phi_{s+T}^i - \phi_{s}^i  = 
    \begin{dcases}
    (\lambda - 1) \phi_s^i, \quad &i=0, \\
    (\lambda - 1) \phi_s^i + \phi_{s}^{i-1}, \quad &i=1,\dots,m_\lambda-1,
    \end{dcases}
\end{equation*}
that $\tau \mapsto \phi_\tau^i$ is $T$-periodic if and only if $\lambda = 1$ and $i=0$, and $T$-antiperiodic if and only if $\lambda = -1$ and $i=0$. However, in the upcoming characterization of the center manifold, we explicitly need a $T$-(anti)periodic $C^{k+1}$-smooth (generalized) eigenbasis. To construct this basis for the $T$-periodic setting, let us first introduce some notation. For a (complex) Banach space $E$ and $l \in \mathbb{N}$, we define $C_T^l(\mathbb{R},E)$ as the (complex) Banach space consisting of $C^l$-smooth $T$-periodic $E$-valued functions defined on $\mathbb{R}$ equipped with the standard $C^l$-norm. Furthermore, let us introduce the unbounded linear operators $\mathcal{A} : \mathcal{D}(\mathcal{A})  \to C_T(\mathbb{R},X)$ and $\mathcal{A}^{\odot \star} : \mathcal{D}(\mathcal{A}^{\odot \star}) \to C_T(\mathbb{R},X^{\odot \star})$ by
\begin{align} \label{eq:D(A)DDE}
\begin{split}
    \mathcal{D}(\mathcal{A}) &:= \{ \varphi \in C_T^1(\mathbb{R},X) : \varphi(\tau) \in \mathcal{D}(A(\tau)) \mbox{ for all $\tau \in \mathbb{R}$} \} \subseteq C_T(\mathbb{R},X), \\
    \mathcal{D}(\mathcal{A}^{\odot \star}) &:= \{ \varphi^{\odot \star} \in C_T^1(\mathbb{R},X^{\odot \star}) : \varphi^{\odot \star}(\tau) \in \mathcal{D}(A^{\odot \star}(\tau)) \mbox{ for all $\tau \in \mathbb{R}$} \} \subseteq C_T(\mathbb{R},X^{\odot \star}),
\end{split}
\end{align}
with action
\begin{equation} \label{eq:curlyAsunstar}
    \mathcal{A}\varphi := ( \tau \mapsto A(\tau)\varphi(\tau) - \dot{\varphi}(\tau)), \quad \mathcal{A}^{\odot \star} \varphi^{\odot \star} := ( \tau \mapsto A^{\odot \star}(\tau)\varphi^{\odot \star}(\tau) - \dot{\varphi}^{\odot \star}(\tau)),
\end{equation}
and define the linear (canonical) embedding $\iota : C_T(\mathbb{R},X) \to C_T(\mathbb{R},X^{\odot \star})$ by $\iota \varphi := (\tau \mapsto j \varphi(\tau))$. Note that $\iota \varphi$ takes values in $X^{\odot \odot}$ due to the $\odot$-reflexivity of $X$ with respect to the shift semigroup $T_0$. These operators can also be introduced on the space of $T$-antiperiodic continuous $X$-valued functions as this will be imported when we will study Floquet multipliers $\lambda \in \mathbb{R}_{-} := (-\infty,0)$, see \cref{prop:eigenfunctions2}. For the sake of simplicity, we will work from now on with $T$-periodic functions and remark at the end of this subsection on the $T$-antiperiodic setting. The following result is a generalization from finite-dimensional ODEs \cite[Proposition III.1]{Iooss1999} towards classical DDEs.

\begin{theorem} \label{thm:eigenfunctions}
Let $\lambda \in \mathbb{C} \setminus \mathbb{R}_{-}$ be a Floquet multiplier of algebraic multiplicity $m_\lambda$ with $\sigma$ its associated Floquet exponent. Then there exist $\varphi_i \in C_T^{k+1}(\mathbb{R},X)$ satisfying
\begin{equation} \label{eq:ODEeig2}
    ( \mathcal{A} - \sigma I)\varphi_i=
    \begin{cases}
    0, \quad &i = 0,\\
    \varphi_{i-1}, \quad &i=1,\dots,m_\lambda-1,
    \end{cases}
\end{equation}
or equivalently
\begin{equation} \label{eq:ODEeig}
    (\mathcal{A}^{\odot \star} - \sigma I)\iota\varphi_i=
    \begin{cases}
    0, \quad &i = 0,\\
    \iota\varphi_{i-1}, \quad &i=1,\dots,m_\lambda-1,
    \end{cases}
\end{equation}
such that the set of functions $\{ \varphi_0(\tau),\dots,\varphi_{m_\lambda - 1}(\tau) \}$ is an ordered basis of $E_\lambda(\tau)$ for all $\tau \in \mathbb{R}$.
\end{theorem}

\begin{proof}
Let $s \in \mathbb{R}$ be a starting time and consider the basis $\{ \phi_s^0, \dots, \phi_s^{m_\lambda - 1} \}$ from \eqref{eq:generalizedeigenvecU(s+T,s)}. We show the claim by induction on $i \in \{ 0,\dots, m_\lambda - 1 \}$. For the base case ($i=0$), consider the initial value problem
\begin{equation} \label{eq:IVPzeta1}
    \begin{dcases}
    (d^\star - A^{\odot \star}(\tau) + \sigma)j\varphi_0(\tau) = 0, \quad \tau \geq s,\\
    \varphi_0(s) = \phi_s^0,
    \end{dcases}
\end{equation}
where $d^\star$ denotes the weak$^\star$ differential operator from \Cref{def:wkstarderivative}. It follows from \eqref{eq:IVPzeta1} that
\begin{align*}
    d^\star (j \circ e^{\sigma(\cdot -s)} \varphi_0)(\tau) &= \sigma e^{\sigma(\tau-s)} j\varphi_0(\tau) + e^{\sigma(\tau-s)} d^\star j\varphi_0(\tau) \\
    &= e^{\sigma(\tau-s)} ( d^\star + \sigma ) j\varphi_0(\tau) = A^{\odot \star}(\tau)j(e^{\sigma(\tau-s)}\varphi_0(\tau)).
\end{align*}
This differential equation is of the form \cite[Equation (4.10)]{Clement1988} and thus its unique solution \cite[Theorem 4.14]{Clement1988} on $[s,\infty)$ is given by
\begin{equation} \label{eq:zeta1}
    e^{\sigma(\tau-s)}\varphi_0(\tau) = U(\tau,s)\varphi_0(s),
\end{equation}
whenever $\varphi_0(s) \in j^{-1} \mathcal{D}(A_0^{\odot \star})$. Since $\varphi_0(s) = \phi_s^0$ the claim follows from \Cref{lemma:basisD(A)} because $\phi_s^0 \in \mathcal{D}(A(s)) \subseteq j^{-1} \mathcal{D}(A_0^{\odot \star})$, recall \cref{sec:periodic center manifolds}. Let us now prove the $T$-periodicity of $\varphi_0$. Choosing $\tau = s+T$ in \eqref{eq:zeta1} and using \eqref{eq:generalizedeigenvecU(s+T,s)} yields
\begin{equation} \label{eq:Tperiodiceig}
    e^{\sigma T}\varphi_0(s+T) = U(s+T,s)\varphi_0(s) = \lambda \varphi_0(s).
\end{equation}
Because $\lambda = e^{\sigma T}$ is nonzero we get $\varphi_0(s + T) = \varphi_0(s)$, and so $\varphi_0$ is $T$-periodic and extends to $\mathbb{R}$. To prove the smoothness assertion, recall from \Cref{lemma:basisD(A)} that $ \tau \mapsto \phi_\tau^0 = U(\tau,s)\phi_s^0$ is $C^{k+1}$-smooth, and because $\tau \mapsto e^{-\sigma(\tau-s)}$ is analytic, it is clear from \eqref{eq:zeta1} that the map $\varphi_0$, given by $\varphi_0(\tau) = e^{-\sigma(\tau-s)}\phi_\tau^0$ for all $\tau \in \mathbb{R}$, is $C^{k+1}$-smooth. Hence, the weak$^\star$ differential operator $d^\star$ in \eqref{eq:IVPzeta1} can be replaced by a standard differential operator for $i=0$, proving \eqref{eq:ODEeig}. By linearity, we have that $\varphi_0(\tau) \in \mathcal{D}(A(\tau))$ for all $\tau \in \mathbb{R}$ which proves the base case for \eqref{eq:ODEeig} and \eqref{eq:ODEeig2}.

To complete the induction, assume that the maps $\varphi_0,\dots,\varphi_{i-1} \in C_T^{k+1}(\mathbb{R},X)$ are constructed for some $i \in \{1,\dots,m_\lambda - 1 \}$ and consider the initial value problem
\begin{equation} \label{eq:IVPzetai}
    \begin{dcases}
    (d^\star - A^{\odot \star}(\tau) + \sigma)j\varphi_i(\tau) = -j\varphi_{i-1}(\tau), \quad \tau \geq s, \\
    \varphi_i(s) = \sum_{k = 0}^{i} \alpha_{ik} \phi_s^k,
    \end{dcases}
\end{equation}
where $\phi_s^0,\dots,\phi_s^i \in \mathcal{D}(A(s))$ are from the Jordan chain \eqref{eq:generalizedeigenvecU(s+T,s)}. The goal is to find scalars $\alpha_{ik}$ such that $\varphi_i$ becomes $T$-periodic. A similar computation as done for the base case tells us by using \eqref{eq:IVPzetai} that
\begin{equation} \label{eq:AODEinhom}
    d^\star (j \circ e^{\sigma(\cdot-s)} \varphi_i)(\tau) =  e^{\sigma(\tau-s)}( d^\star  + \sigma ) j\varphi_i(\tau) = A^{\odot \star}(\tau)j(e^{\sigma(\tau-s)}\varphi_i(\tau)) - j(e^{\sigma (\tau-s)}\varphi_{i-1}(\tau)).
\end{equation}
We will show that this differential equation admits a unique solution on $[s,\infty)$. Consider the function $w_i : [s,\infty) \to X$ defined by
\begin{equation} \label{eq:w_i}
    w_i(\tau) := U(\tau,s) \sum_{k=0}^{i} \frac{ (s-\tau)^k}{k!} \varphi_{i-k}(s), \quad \forall \tau \in [s,\infty).
\end{equation}
Because $\phi_s^0,\dots,\phi_s^i \in \mathcal{D}(A(s))$, we obtain from \eqref{eq:IVPzetai} that $\varphi_0(s),\dots,\varphi_{i}(s) \in \mathcal{D}(A(s))$ by linearity. It is clear from \cref{lemma:basisD(A)} that
\begin{equation*}
    U(\tau,s)\varphi_{j}(s) = \sum_{k=0}^j \alpha_{jk}\phi_\tau^k \in \mathcal{D}(A(\tau)), \quad \forall \tau \in [s,\infty), \ j \in \{0,\dots,i\},
\end{equation*}
and so it follows that $\tau \mapsto w_i(\tau)$ takes values in $\mathcal{D}(A(\tau)) \subseteq j^{-1} \mathcal{D}(A_0^{\odot \star})$ and is $C^{k+1}$-smooth, which implies the weak$^\star$ differentiability of $w_i$. Clearly, $w_i(s) = \varphi_i(s)$ and notice that
\begin{align*}
    d^\star (j \circ w_i)(\tau) &= A^{\odot \star}(\tau)jU(\tau,s) \sum_{k=0}^i \frac{(s-\tau)^k}{k!} \varphi_{i-k}(s)  - jU(t,s) \sum_{k=0}^{i-1} \frac{(s-\tau)^k}{k!}\varphi_{i-k}(s) \\
    &= A^{\odot \star}(\tau)jw_i(\tau) - jw_{i-1}(\tau),
\end{align*}
which proves that $w_i$ is a solution of \eqref{eq:AODEinhom} on $[s,\infty)$. Since $w_{i-1}$ is at least continuous, it follows from \cite[Proposition 13]{Lentjes2023} and by construction of $w_i$ that \eqref{eq:AODEinhom} admits a unique solution $j \circ w_i$ on $[s,\infty)$, where $w_i = e^{\sigma(\cdot - s)}\varphi_i$. As a consequence, $\varphi_i = e^{-\sigma(\cdot - s)} w_i$ is $C^{k+1}$-smooth.

Let us now turn our attention towards proving $T$-periodicity. We see from \eqref{eq:w_i} that $\varphi_i(s) = \varphi_i(s+T)$ if and only if
\begin{equation} \label{eq:thmU(s+T,s)-lambda}
    (U(s+T,s) - \lambda I)\varphi_i(s) = U(s+T,s) \sum_{k=1}^{i} \frac{(-1)^{k+1} T^k}{k!} \varphi_{i-k}(s). 
\end{equation}
Recall from \eqref{eq:IVPzetai} that $\varphi_i(s)=\sum_{k = 0}^{i} \alpha_{ik}\phi_s^k$ and retrieving \eqref{eq:generalizedeigenvecU(s+T,s)} yields
\begin{align*}
    \sum_{k=1}^{i} \alpha_{ik} \phi_s^{k-1} &= U(s+T,s) \sum_{l=1}^{i} \sum_{k=1}^{i-l} \alpha_{i-l,k}\frac{(-1)^{l+1} T^l}{l!} \phi_{s}^k \\
    &= \sum_{l=1}^{i} \sum_{k=0}^{i-l} \alpha_{i-l,k}\frac{(-1)^{l+1} T^l}{l!} \begin{cases}
    \lambda \phi_s^k, \quad &k =0, \\
    \lambda \phi_s^k + \phi_s^{k-1}, \quad &k=1,\dots,i-l.
    \end{cases}
\end{align*}
Because the right-hand side is a known element in the subspace spanned by $\phi_s^0,\dots,\phi_s^{i-1}$, the $\alpha_{ik}$'s are uniquely determined for $k=0,\dots,i$ and so we have proven that $\varphi_i$ is $T$-periodic. Hence, $\tau \mapsto \varphi_i(\tau)$ extends to a $C^{k+1}$-smooth solution on $\mathbb{R}$ and takes values in $\mathcal{D}(A(\tau))$. Similarly, the weak$^\star$ differential operator from \eqref{eq:IVPzetai} can be replaced by a standard differential operator and so the formulas \eqref{eq:ODEeig} and \eqref{eq:ODEeig2} hold. One can easily verify from the relations above that $\alpha_{ii} = (\lambda T)^{i}$ for all $i \in \{0,\dots,m_\lambda-1\}$. Hence, $\varphi_0(s),\dots,\varphi_{m_\lambda - 1}(s)$ are all linearly independent and thus form a basis of $E_\lambda(s)$. Therefore, the functions $\varphi_0,\dots,\varphi_{m_\lambda-1}$ are all linearly independent, as they are all solutions to the equation $(\mathcal{A} - \sigma I)^{m_\lambda} \varphi = 0$.
\end{proof} 
Let us now take the time to discuss the connection between the $T$-periodicity and history \eqref{eq:history} of the (generalized) eigenfunctions for Floquet multipliers $\lambda \in \mathbb{C} \setminus \mathbb{R}_{-}$. It is clear that $\{\phi_\tau^{0},\dots,\phi_\tau^{m_\lambda - 1} \}$ is (in general) a non-$T$-periodic basis of $E_\lambda(\tau)$ in terms of history, as they satisfy the periodic linear DDE \eqref{eq:T-LDDEphi}. On the other hand, \Cref{thm:eigenfunctions} shows us that $\{\varphi_{0}(\tau),\dots,\varphi_{m_\lambda-1}(\tau) \}$ is a $T$-periodic basis of $E_\lambda(\tau)$, but how is this basis related to history \eqref{eq:history}? Note that a function $\varphi \in C_T(\mathbb{R},X)$ can be written in terms of history, that is $\varphi(\tau) = \varphi_\tau$, if and only if it satisfies the \emph{translation property} $\varphi(\tau)(\theta) = \varphi(\tau + \theta)(0)$ for all $\tau \in \mathbb{R}$ and $\theta \in [-h,0]$. Moreover, if $\varphi \in C_T(\mathbb{R},X)$ is differentiable in both arguments, $\varphi$ has the translation property if and only if it satisfies the transport equation
\begin{equation*}
    \frac{\partial}{\partial \tau} \varphi(\tau)(\theta) = \frac{\partial}{\partial \theta} \varphi(\tau)(\theta), \quad \forall \tau \in \mathbb{R}, \ \theta \in [-h,0].
\end{equation*}
However, \eqref{eq:ODEeig2} in combination with \eqref{eq:D(As)} shows us directly that 
\begin{equation*}
    \bigg( \frac{\partial}{\partial \theta} - \frac{\partial}{\partial \tau} \bigg)\varphi_i(\tau)(\theta)=
    \begin{cases}
    \sigma \varphi_i(\tau)(\theta), \quad &i = 0,\\
    \varphi_{i-1}(\tau)(\theta) +  \sigma \varphi_i(\tau)(\theta), \quad &i=1,\dots,m_\lambda-1,
    \end{cases}
\end{equation*}
which proves that $\varphi_i$ has the translation property if and only if $\sigma = 0$ and $i=0$. Hence, the only $T$-periodic (generalized) eigenfunction that satisfies the translation property is the map $\tau \mapsto \dot{\gamma}_\tau$. It is however the $T$-periodic basis $\{\varphi_{0}(\tau),\dots,\varphi_{m_\lambda-1}(\tau) \}$ of $E_{\lambda}(\tau)$ that is needed for the characterization of $\mathcal{W}_{\loc}^c(\Gamma)$ in \Cref{sec:characterization}.

\begin{remark} \label{remark:identification}
Recall from \cref{sec:periodic center manifolds} that the (complexified) Banach spaces $X^{\odot}, X^{\odot \star}$ and $X^{\odot \odot}$ can be represented as $E = E_1 \times E_2,$ where $E_{1,2}$ are (complex) Banach spaces with respect to some norm $\|\cdot\|_{E_{1,2}}$. It is convenient to make for any $l \in \mathbb{N}$ the identification
\begin{equation*}
    C_T^l(\mathbb{R},E) \cong C_T^l(\mathbb{R},E_1)  \times C_T^l(\mathbb{R},E_2), \\
\end{equation*}
so that we can write any $\varphi \in C_T^l(\mathbb{R},E)$ as $\varphi = (\varphi_1,\varphi_2)$ where $\varphi_i \in  C_T^l(\mathbb{R},E_i)$ for $i=1,2$. With this identification, we define the norm $\| \cdot \|_{C^l}$ on $C_T^l(\mathbb{R},E_1)  \times C_T^l(\mathbb{R},E_2)$ by
\begin{equation*}
    \| \varphi \|_{C^l} := \| \varphi_1 \|_{C_1^l} + \| \varphi_2 \|_{C_2^l}, \quad \| \varphi_i \|_{C_i^l} := \sum_{k=0}^l \sup_{\tau \in \mathbb{R}} \| \varphi_i^{(k)}(\tau) \|_{E_i},
\end{equation*}
where $\varphi_i^{(k)}$ denotes the $k$th derivative of $\varphi_i$ for $i=1,2$. When $l = 0$, we simply write $\| \cdot \|_\infty := \| \cdot \|_{C^0}$. The same notation and norms will be used when we work on subspaces of a Banach space $E$, recall for example $\mathcal{D}(\mathcal{A}^{\odot \star})$ that is a subspace of $X^{\odot \star}$. \hfill $\lozenge$
\end{remark}

\cref{thm:eigenfunctions} tells us that $\{ \sigma \in \mathbb{C} : \sigma \mbox{ is a Floquet exponent} \}$ is a subset of the point spectrum of the unbounded linear operators $\mathcal{A}$ and $\mathcal{A}^{\odot \star}$. It turns out that the spectral structure of these, and two other upcoming operators from \cref{subsec: dual time periodic smooth Jordan}, will be of great importance in our upcoming paper on the study of bifurcations of limit cycles in classical DDEs. Therefore, we will study the spectral structure of all four operators in \cref{subsec: duality relations}. To have a rather reasonable spectral theory for unbounded linear operators, one has to prove that the operators of interest are at least closable, recall for example \cite[Appendix B]{Engel2000} for definitions and more information on this topic.

\begin{proposition} \label{prop:closable}
    The unbounded linear operators $\mathcal{A}$ and $\mathcal{A}^{\odot \star}$ are not closed but closable. 
\end{proposition}
\begin{proof}
Let us first prove that $\mathcal{A}$ is closable. Let $(\varphi_m)_m$ be a sequence in $\mathcal{D}(\mathcal{A})$ converging in norm to zero, and assume that the sequence $(\mathcal{A}\varphi_m)_m$ converges in norm to some $\psi \in C_T(\mathbb{R},X)$. We have to show that $\psi = 0$. Set $\psi_m = \mathcal{A}\varphi_m$ and notice that this is equivalent to
\begin{equation*}
    \bigg(\frac{\partial}{\partial \theta} - \frac{\partial}{\partial \tau} \bigg) \varphi_m(\tau)(\theta) = \psi_m(\tau)(\theta),  \quad \forall \tau \in \mathbb{R}, \ \theta \in [-h,0],
\end{equation*}
where $\varphi_m \in \mathcal{D}(\mathcal{A})$ for all $m \in \mathbb{N}$. This is an (inhomogeneous) transport equation that has the solution
\begin{equation*}
    \varphi_m(\tau)(\theta) = \varphi_m(\tau+\theta)(0) + \int_0^\theta \psi_m(\tau+\theta-s)(s) ds, \quad \varphi_m(\tau)'(0) = L(\tau) \varphi_m(\tau).
\end{equation*}
Since the convergence of $\varphi_m \to 0$ and $\psi_m \to \psi$ is uniform in the first and second variable as $m\to \infty$, we obtain
\begin{align}
    \begin{split} \label{eq:integralpsi}
    \int_0^\theta  \psi(\tau+\theta-s)(s) ds &= \lim_{m \to \infty} \int_0^\theta  \psi_m(\tau+\theta-s)(s) ds \\
    &= \lim_{m \to \infty} (\varphi_m(\tau)(\theta) - \varphi_m(\tau+\theta)(0)) = 0, \quad \forall \tau \in \mathbb{R}, \ \theta \in [-h,0].
    \end{split}
\end{align}
Let $\psi_i$ denote the $i$th component of $\psi$ for $i=1,\dots,n$. To prove the claim, let $\tau_1 \in \mathbb{R}$ and $\theta_1,\theta_2 \in [-h,0]$ with $\theta_1 \neq \theta_2$ be given. If $\tau_2 = \tau_1 + \theta_1 - \theta_2$, then according to \eqref{eq:integralpsi} we obtain
\begin{align*}
    \int_{\theta_1}^{\theta_2} \psi_i(\tau_1+\theta_1-s)(s) ds &= \int_{0}^{\theta_1} \psi_i(\tau_1+\theta_1-s)(s) ds + \int_{\theta_1}^{\theta_2} \psi_i(\tau_1+\theta_1-s)(s) ds \\
    &= \int_{0}^{\theta_2} \psi_i(\tau_1+\theta_1-s)(s) ds = \int_{0}^{\theta_2} \psi_i(\tau_2+\theta_2-s)(s) ds = 0
\end{align*}
since the second and last integral vanish due to \eqref{eq:integralpsi}. Assume now without loss of generality that $\theta_2 > \theta_1$. Since $\psi$ is continuous, the real-valued function $ s \mapsto \operatorname{Re} \psi_i(\tau_1+\theta_1-s)(s)$ is continuous on the compact set $[\theta_1,\theta_2]$. Hence, the mean value theorem for integrals applies and there exists a $c_i \in (\theta_1,\theta_2)$ such that $\operatorname{Re} \psi_i(\tau_1+\theta_1-c_i)(c_i) = 0$ since $\theta_2 \neq \theta_1$. By the squeeze theorem, $c_i \downarrow \theta_1$ when $\theta_2 \downarrow \theta_1$ and we obtain $ \operatorname{Re} \psi_i(\tau_1)(\theta_1) = 0$. We obtain similarly that $ \operatorname{Im} \psi_i(\tau_1)(\theta_1) = 0$ and thus $\psi_i(\tau_1)(\theta_1) = 0$. Because the index $i$ and values $\tau_1$ and $\theta_1$ were chosen arbitrary, we obtain $\psi = 0$, i.e. $\mathcal{A}$ is a closable operator. 

Let us now prove that $\mathcal{A}^{\odot \star}$ is closable. Let $(\varphi_m^{\odot \star})_m$ be a sequence in $\mathcal{D}(\mathcal{A}^{\odot \star})$ converging in norm to zero, and assume that the sequence $(\mathcal{A}^{\odot \star}\varphi_m^{\odot \star})_m$ converges in norm to some $\psi^{\odot \star} \in C_T(\mathbb{R},X^{\odot \star})$. Set $\psi_m^{\odot \star} = \mathcal{A}^{\odot \star}\varphi_m^{\odot \star}$ with $\varphi_m^{\odot \star} = (\alpha_m,\varphi_m)$ and $\psi_m^{\odot \star} = (\beta_m,\psi_m)$, recall \cref{remark:identification}. Using the definition of $\mathcal{A}^{\odot \star}$ from \eqref{eq:Asunstartau}, we obtain along the same lines of the proof on the closability of $\mathcal{A}$, that
\begin{equation} \label{eq:varphiclosable}
    \varphi_m(\tau)(\theta) = \alpha_m(\tau+\theta) + \int_0^\theta \psi_m(\tau+\theta-s)(s) ds,
\end{equation}
where the sequence $(\alpha_m)_m$ in $C_T^1(\mathbb{R},\mathbb{C}^n)$ satisfies the periodic linear (inhomogeneous) DDE
\begin{equation} \label{eq:alphan}
    \dot{\alpha}_m(\tau) = L(\tau)[\alpha_m]_\tau + H_m(\tau), \quad H_m(\tau) = L(\tau) \bigg[\theta \mapsto \int_0^\theta \psi_m(\tau+\theta-s)(s) ds \bigg] - \beta_m(\tau).
\end{equation}
The convergence of $\varphi_m^{\odot \star} \to 0$ implies that $\alpha_m \to 0$ and $\varphi_m \to 0$ in norm as $m \to \infty$. Hence, we obtain from \eqref{eq:varphiclosable}, and a similar argument as in the proof on the closability of $\mathcal{A}$, that
\begin{equation} \label{eq:integralpsin}
    \int_{\theta_1}^{\theta_2} \psi_m(\tau-\theta_1+s)(s) ds \to 0, \quad m \to \infty, \quad \forall \tau \in \mathbb{R}, \ \theta_1,\theta_2 \in [-h,0].
\end{equation}
Assume without loss of generality that $\theta_2 > \theta_1$ and consider for a fixed $\tau \in \mathbb{R}$ the map $r_m \in L^\infty([\theta_1,\theta_2],\mathbb{C}^{n})$ defined by $s \mapsto \psi_m(\tau-\theta_1+s)(s)$ for all $m \in \mathbb{N}$. The sequence of functions $(r_m)_m$ converges in $L^\infty$ towards $r \in L^\infty([\theta_1,\theta_2],\mathbb{C}^n)$ defined by $s \mapsto \psi(\tau+\theta_1-s)(s)$. Hence, $(r_m)_m$ converges almost everywhere pointwise towards $r$ as $m \to \infty$. By the convergence obtained in \eqref{eq:integralpsin}, the sequence of functions $(r_m)_m$ is dominated by some integrable function and so the dominated convergence theorem applies to \eqref{eq:integralpsin}. Thus,
\begin{equation*}
    \int_{\theta_1}^{\theta_2} \psi_i(\tau-\theta_1+s)(s) ds = 0, \quad \forall \tau \in \mathbb{R}, \ \theta_1,\theta_2 \in [-h,0],
\end{equation*}
where $\psi_i$ denotes the $i$th component of $\psi$ for $i=1,\dots,n$. Since $r \in L^{\infty}([\theta_1,\theta_2], \mathbb{C}^{n}) \subseteq L^{1}([\theta_1,\theta_2], \mathbb{C}^{n})$ because $[\theta_1,\theta_2] \subseteq [-h,0]$ is bounded, Lebesgue's differentiation theorem tells us that $\psi_i = 0$ almost everywhere for all $i=1,\dots,n$ and so $\psi = 0$ almost everywhere. It remains to show that $\beta = 0$. Since $\alpha_m,\dot{\alpha}_m,\psi_m \to 0$ in norm as $m \to \infty$, it follows directly from \eqref{eq:alphan} that $\beta = 0$ and so $\psi^{\odot \star} = 0$.

Before we prove that $\mathcal{A}^{\odot \star}$ is not closed, let us first show that the unbounded linear operator $\mathcal{A}_0^{\odot \star} : \mathcal{D}(\mathcal{A}_0^{\odot \star}) \to C_T(\mathbb{R},X^{\odot \star})$ defined by
\begin{equation*}
    \mathcal{D}(\mathcal{A}_0^{\odot \star}) = \mathcal{D}(\mathcal{A}^{\odot \star}), \quad \mathcal{A}_0^{\odot \star}(\alpha,\varphi) = (-\dot{\alpha}, \tau \mapsto \varphi(\tau)' - \dot{\varphi}(\tau)),
\end{equation*}
is not closed. To do this, let $(f_m)_m$ be a sequence in $C_T^1(\mathbb{R},\mathbb{C}^n)$ converging uniformly to some $f \in C_T(\mathbb{R},\mathbb{C}^n) \setminus C_T^1(\mathbb{R},\mathbb{C}^n)$. Define the sequence $(\varphi_m)_m$ in $C_T^1(\mathbb{R},C^1([-h,0],\mathbb{C}^n))$ by $\varphi_m(\tau)(\theta) = \theta^2 f_m(\tau+\theta)$ and let $\varphi(\tau)(\theta) = \theta^2 f(\tau+\theta)$. Then $((0,\varphi_m))_m$ is a sequence in $\mathcal{D}(\mathcal{A}_0^{\odot \star})$, and $\| (0,\varphi_m) - (0,\varphi) \| \leq h^2 \| f_m - f \|_\infty \to 0$ as $m \to \infty$ and $(\mathcal{A}_0^{\odot \star}(0,\varphi_m))(\tau) = (0, \theta \mapsto 2 \theta f_m(\tau+\theta))$. Hence, $(\mathcal{A}_0^{\odot \star}(0,\varphi_m))_m$ converges in norm to $(0, \tau \mapsto (\theta \mapsto 2 \theta f(\tau+\theta)))$. However, $\varphi$ is not continuously differentiable in the first variable and so $\varphi \notin \mathcal{D}(\mathcal{A}_0^{\odot \star})$, i.e. $\mathcal{A}_0^{\odot \star}$ is not closed.

To prove that $\mathcal{A}^{\odot \star}$ is not closed, introduce the perturbation $\mathcal{B} : C_T(\mathbb{R},X^{\odot \star}) \to C_T(\mathbb{R},X^{\odot \star})$ by $\mathcal{B} \varphi^{\odot \star} = (\tau \mapsto L(\tau)\varphi(\tau),0)$, where $\varphi^{\odot \star}(\tau) = (\varphi(\tau)(0), \varphi(\tau)) \in X^{\odot \star}$ and notice that $\mathcal{B} \in \mathcal{L}(C_T(\mathbb{R},X^{\odot \star}))$ since $L \in C_T(\mathbb{R},\mathcal{L}(X,\mathbb{R}^n))$. Recall from \eqref{eq:curlyAsunstar} that $\mathcal{A}^{\odot \star} = \mathcal{A}_{0}^{\odot \star} + \mathcal{B}$, and since $\mathcal{A}_{0}^{\odot \star}$ is the sum of a non-closed linear operator and a bounded linear operator, we have that $\mathcal{A}^{\odot \star}$ is not closed. Since this counterexample is constructed with a sequence of $C^1$-smooth functions, the same strategy can be applied to show that $\mathcal{A}^{\odot \star}$ with restricted domain $\mathcal{D}^{\odot \star} := \iota(C_T^1(\mathbb{R},C^{1}([-h,0],\mathbb{C}^n)))$ is not closed, where $\iota$ is introduced below \eqref{eq:curlyAsunstar}.

Let us now prove that $\mathcal{A}$ is not closed. Since $\mathcal{A}^{\odot \star}$ with restricted domain $\mathcal{D}^{\odot \star}$ is not closed, there exists a sequence $(\iota\varphi_m)_m$ in $\mathcal{D}^{\odot \star}$ with $\varphi_m \in C_T^1(\mathbb{R},C^{1}([-h,0],\mathbb{C}^n))$ such that 
$\varphi_m \to \varphi \in C_T(\mathbb{R},X)$ and $\mathcal{A}^{\odot \star} \iota \varphi_m = \psi_m^{\odot \star} \to \psi^{\odot \star}$ in norm as $m \to \infty$, but $\iota \varphi \notin \mathcal{D}^{\odot \star}$ or $\mathcal{A}^{\odot \star}\iota \varphi \neq \psi^{\odot \star}$, where $\psi^{\odot \star} = (\beta,\psi)$ and $\psi_m^{\odot \star} = (\beta_m,\psi_m)$. Our aim is to show that $\varphi \notin  C_T^1(\mathbb{R},C^{1}([-h,0],\mathbb{C}^n))$. Recall first that $\varphi_m$ and $\alpha_m = \varphi_m(\cdot)(0)$ satisfy \eqref{eq:varphiclosable} and \eqref{eq:alphan} respectively. Since $\psi_m$ and $\psi$ are $C^1$-smooth in the first component and essentially bounded in the second component, the maps $\theta \mapsto \int_0^\theta \psi_m(\tau+\theta-s)(s) ds$ and $\theta \mapsto \int_0^\theta \psi(\tau+\theta-s)(s) ds$ are in $X$ for all $m \in \mathbb{N}$ and fixed $\tau \in \mathbb{R}$. Introduce the sequence $(H_m)_m$ as in \eqref{eq:alphan} and notice that $H_m \to H \in C_T(\mathbb{R},\mathbb{C}^n)$, where $H(\tau) = L(\tau)[\theta \mapsto \int_0^\theta \psi(\tau+\theta-s)(s) ds] - \beta(\tau)$ is well-defined. Since in addition $L \in C_T(\mathbb{R},\mathcal{L}(X,\mathbb{R}^n))$ and $H_m$ are all continuous, we obtain $\dot{\alpha}(\tau) = L(\tau)\alpha_\tau + H(\tau)$. Thus, $\| \dot{\alpha}_m - \dot{\alpha} \|_{\infty} \to 0$ as $m \to \infty$ and so $\alpha \in C_T^1(\mathbb{R},\mathbb{C}^n)$. Applying the dominated convergence theorem on $(\psi_m)_m$ as performed similarly to the proof on the closability of $\mathcal{A}^{\odot \star}$ below \eqref{eq:integralpsin}, and recalling the uniform convergence of $(\alpha_m)_m$ to $\alpha$, yields
\begin{equation} \label{eq:phipsiDDE}
     \varphi(\tau)(\theta) = \alpha(\tau+\theta) + \int_0^\theta \psi(\tau+\theta-s)(s) ds,
\end{equation}
and so $\varphi(\cdot)(0) = \alpha$, i.e. $\iota \varphi = (\alpha,\varphi)$. If $\iota \varphi \in \mathcal{D}^{\odot \star}$, then it follows from \eqref{eq:phipsiDDE} that $\mathcal{A}^{\odot \star}\iota \varphi = \psi^{\odot \star}$, which is not possible by the construction of $\varphi$. Thus, we conclude that $\varphi \notin C_T^1(\mathbb{R},C^{1}([-h,0],\mathbb{C}^n))$. To prove that $\mathcal{A}$ is not closed, take the sequence $(\varphi_m)_m$ from before with uniform limit $\varphi \notin C_T^1(\mathbb{R},C^{1}([-h,0],\mathbb{C}^n))$. Then $(\varphi_m)_m$ is by construction a sequence in $\mathcal{D}(\mathcal{A})$ since $\varphi_m(\tau)'(0) = L(\tau)\varphi_m(\tau)$ for all $\tau \in \mathbb{R}$ is equivalent to having $T$-periodic solutions $\alpha_m = \varphi_m(\cdot)(0)$ of \eqref{eq:alphan} for all $m \in \mathbb{N}$. However, $\varphi \notin C_T^1(\mathbb{R},C^{1}([-h,0],\mathbb{C}^n))$ and so it follows from \eqref{eq:D(As)} and \eqref{eq:D(A)DDE} that $\varphi \notin \mathcal{D}(\mathcal{A})$, i.e. $\mathcal{A}$ is not closed.
\end{proof}
To get more insight into the structure of the linear operators $\mathcal{A}$ and $\mathcal{A}^{\odot \star}$, observe that we can decompose $\mathcal{A}$ as the difference of the unbounded linear operators
\begin{equation} \label{eq:Adecomp}
    \mathcal{D}(\mathcal{A}) = \mathcal{D}(M_A) \cap \mathcal{D}(D), \quad \mathcal{A} = M_A - D,
\end{equation}
where $M_A : \mathcal{D}(M_A) \to C_T(\mathbb{R},X)$ is the (complexified) multiplication operator (associated to $A$) and $D : \mathcal{D}(D) \to C_T(\mathbb{R},X)$ is the (complexified) differential operator defined by
\begin{alignat}{2}
    \mathcal{D}(M_A) &=\{ \varphi \in C_T(\mathbb{R},X) : \varphi(\tau) \in \mathcal{D}(A(\tau)) \mbox{ for all } \tau \in \mathbb{R} \} \subseteq C_T(\mathbb{R},X), \quad &&M_A\varphi = A(\cdot)\varphi(\cdot),  \label{eq:Arecht}\\
    \mathcal{D}(D) &= C_T^1(\mathbb{R},X) \subseteq C_T(\mathbb{R},X), \quad &&D\varphi = \dot{\varphi} \label{eq:Drecht}.   
\end{alignat}
\begin{lemma} \label{lemma:MADclosed}
The operators $M_A$ and $D$ are densely defined and closed.
\end{lemma}
\begin{proof}
It is well-known that $D$ is a closed and densely defined linear operator on $C_T(\mathbb{R},X)$. To prove that $M_A$ is closed, let $(\varphi_m)_m$ be a sequence in $\mathcal{D}(M_A)$ converging in norm to some $\varphi \in C_T(\mathbb{R},X)$ and let $(M_A\varphi_m)_m$, where $M_A\varphi_m = \psi_m$, converge in norm to some $\psi \in C_T(\mathbb{R},X)$. It follows immediately from \eqref{eq:D(As)} that $\varphi_m$ has the form
\begin{equation*}
    \varphi_m(\tau)(\theta) = \varphi_m(\tau)(0) + \int_0^\theta \psi_m(\tau)(s) ds, \quad \varphi_m(\tau)'(0) = \psi_m(0) = L(\tau)\varphi_m(\tau),
\end{equation*}
for all $\tau \in \mathbb{R}$ and $\theta \in [-h,0]$. Similar to the argument provided in \eqref{eq:integralpsi}, we obtain
\begin{equation*}
    \varphi(\tau)(\theta) = \varphi(\tau)(0) + \int_0^\theta \psi(\tau)(s) ds, \quad \varphi(\tau)'(0) = L(\tau)\varphi(\tau),
\end{equation*}
where the last equality follows from the fact that $L \in C_T(\mathbb{R},\mathcal{L}(X,\mathbb{C}^n))$. Hence, $\varphi \in \mathcal{D}(M_A)$ and $M_A\varphi = \psi$, which proves the claim. It will be proven in \cref{prop:DAdense} that the subspace $\mathcal{D}(\mathcal{A})$ of $\mathcal{D}(M_A)$ is dense in $C_T(\mathbb{R},X)$ and thus the density of $\mathcal{D}(M_A)$ follows.
\end{proof}
Note that \cref{lemma:MADclosed} does not contradict \cref{prop:closable} as the difference of two closed operators is not necessarily closed, see \cite[Chapter III and Appendix B]{Engel2000} for more information.

The reader with knowledge on evolution semigroups \cite{Chicone1999,Engel2000,Nickel1997,Latushkin1996}, first introduced by Howland \cite{Howland1974}, probably recognizes the decomposition \eqref{eq:Adecomp}, as the generator of the associated evolution semigroup. It turns out in the sequel that the concept of evolution semigroups will provide us a better understanding of the operators $\mathcal{A}$ and $\mathcal{A}^{\odot \star}$. Therefore, let us first define for any $s \geq 0$ the operator $\mathcal{U}(s) : C_T(\mathbb{R},X) \to C_T(\mathbb{R},X)$ by
\begin{equation} \label{eq:evolutionsemi}
    (\mathcal{U}(s)f)(\tau) := U(\tau,\tau-s)f(\tau-s), \quad \forall f \in C_T(\mathbb{R},X),
\end{equation}
and note that $\mathcal{U} := \{\mathcal{U}(s)\}_{s \geq 0}$ is a family of linear operators on $C_T(\mathbb{R},X)$. This family is often called the \emph{evolution semigroup} (associated to $U$), but is in the literature always defined on $C_0(\mathbb{R},X)$ or $L^p(\mathbb{R},X)$ for some $1 \leq p \leq \infty$, instead of $C_T(\mathbb{R},X)$. The reason for this is clear as the existence of an evolution semigroup on $C_0(\mathbb{R},X)$ or $L^p(\mathbb{R},X)$ is equivalent to the existence of unique solutions of abstract non-autonomous Cauchy problems, see for example \cite[Section 9]{Engel2000} and the references therein. Therefore, we can not claim directly all well-known facts about evolution semigroups. 

One key element in our upcoming proofs on properties of $\mathcal{U}$ is the fact that we work on $T$-periodic functions. Therefore, we can reformulate our problems on functions defined on the unit circle, represented as the compact Hausdorff quotient space $\mathbb{T} := \mathbb{R} / T \mathbb{Z}$. Indeed, for any (complex) Banach space $E$ and $l \in \mathbb{N}$, the space $C_T^l(\mathbb{R},E) \cong C^l(\mathbb{T},E)$ via the linear isometric isomorphism $f \mapsto \hat{f}$, where $\hat{f}(\tau \mbox{ mod } T) = f(\tau)$ for all $\tau \in \mathbb{R}$. From now on, we will denote for any $\tau \in \mathbb{R}$ the point $\hat{\tau} := \tau \mbox{ mod } T \in \mathbb{T}$. We start with proving the following, rather classical, result.

\begin{lemma} \label{lemma:C0semigroup}
The evolution semigroup $\mathcal{U}$ is a well-defined $\mathcal{C}_0$-semigroup.
\end{lemma}
\begin{proof}
Let us first observe that $\mathcal{U}(s)$ is well-defined for any $s \geq 0$ since $U$ is a strongly continuous forward evolutionary system, and there holds
\begin{align*}
    (\mathcal{U}(s)f)(\tau+T) &= U(\tau+T,\tau+T-s)f(\tau+T-s) \\
    &= U(\tau,\tau-s)f(\tau-s) = (\mathcal{U}(s)f)(\tau), \quad \forall f \in C_T(\mathbb{R},X),
\end{align*}
where the second equality follows from \eqref{eq:Uperiodic}. Since the bounded linear perturbation $B$ from \eqref{eq:perturbationB} is $T$-periodic, there exist a $M \geq 1$ and $\omega \in \mathbb{R}$ such that $\|U(\tau,s)\| \leq Me^{\omega(\tau-s)}$ for all $\tau \geq s$, see \cite[Theorem XII.2.7]{Diekmann1995}. Therefore, one can easily verify from \eqref{eq:evolutionsemi} that $\mathcal{U}$ is a family of bounded linear operators forming a semigroup that satisfies $\| \mathcal{U}(s) \| \leq M e^{\omega s}$ for all $s \geq 0$.

It remains to show that $\mathcal{U}$ is strongly continuous. To do this, let us reformulate \eqref{eq:evolutionsemi} on the space $C(\mathbb{T},X)$. Let $\varepsilon > 0$ and $f \in C(\mathbb{T},X)$ be given. Since $\mathbb{T}$ is compact, we know that $f$ is uniformly continuous and so there exists a $\delta_1 > 0$ such that $\sup_{\hat{\tau} \in \mathbb{T}} \| f(\tau+s) - f(\tau) \|_\infty < \varepsilon/2$ for all $|s| < \delta_1$, where we recall $\hat{\tau} = \tau \mbox{ mod } T$. Again, by compactness of $\mathbb{T}$, there exist $\tau_1,\dots,\tau_l \in [0,T]$ for some $l \in \mathbb{N}$ such that $\mathbb{T}$ is covered by the union of $l$ open balls centered at $\hat{\tau}_i$ and satisfy in addition
\begin{equation*}
    \sup_{\tau \in [\tau_i, \tau_{i+1}]} \|f(\tau) - f(\tau_i)\|_\infty < \frac{\varepsilon}{4(Me^{\omega s} + 1)}, \quad \forall i \in \{1,\dots,l-1\}. 
\end{equation*}
Due to the strong continuity of the forward evolutionary system $U$, there exists a $\delta_2 > 0$ such that $\|U(\tau+s,\tau)f(\tau_i) - f(\tau_i) \| < \varepsilon/4$ for all $s \in (0,\delta_2), \tau \in [\tau_i, \tau_{i+1}]$ and $i \in \{1,\dots,l-1\}$. Hence,
\begin{align*}
    \sup_{\hat{\tau} \in \mathbb{T}}\|U(\tau+s,\tau)f(\tau) - f(\tau) \|_\infty &\leq \max_{1,\dots,l-1} \sup_{\tau \in [\tau_i,\tau_{i+1}]} \bigg( \|U(\tau+s,\tau)[f(\tau) - f(\tau_i)] \|_\infty \\
    &+ \| U(\tau+s,\tau)f(\tau_i) - f(\tau_i) \|_\infty + \|f(\tau) - f(\tau_i) \|_\infty \bigg) < \frac{\varepsilon}{2}.
\end{align*}
If we let $\delta = \min \{\delta_1, \delta_2\} > 0$ and combine the results from above, we obtain for all $s \in (0,\delta)$ that
\begin{align*}
    \|\mathcal{U}(s)f - f \|_\infty \leq \sup_{\hat{\tau} \in \mathbb{T}} \| U(\tau+s,\tau)f(\tau) - f(\tau) \| + \sup_{\hat{\tau} \in \mathbb{T}} \| f(\tau+s) - f(\tau) \|_\infty < \varepsilon,
\end{align*}
which proves the claim.
\end{proof}
Since $\mathcal{U}$ is a $\mathcal{C}_0$-semigroup on $C_T(\mathbb{R},X)$, it has a generator $\hat{\mathcal{A}}$ with closed and densely defined domain $\mathcal{D}(\hat{\mathcal{A}}) \subseteq C_T(\mathbb{R},X)$ consisting of all functions $f \in C_T(\mathbb{R},X)$ such that $\lim_{s \downarrow 0} \frac{1}{s}(\mathcal{U}(s)f - f)$ exists in $C_T(\mathbb{R},X)$. Clearly, if $f \in \mathcal{D}(\mathcal{A})$, i.e. $f \in C_T^1(\mathbb{R},X)$ and $f(\tau) \in \mathcal{D}(A(\tau))$ for all $\tau \in \mathbb{R}$, then a straightforward calculation shows that $(\hat{\mathcal{A}}f)(\tau) = A(\tau)f(\tau) - \dot{f}(\tau)$ and thus $\mathcal{D}(\mathcal{A}) \subseteq \mathcal{D}(\hat{\mathcal{A}})$, i.e. $\hat{\mathcal{A}}$ is a linear extension of $\mathcal{A}$. Since $\mathcal{A}$ is a closable (\cref{prop:closable}), it admits a smallest closed linear extension, also called the \emph{closure} (of $\mathcal{A}$), and will be denoted by $\overline{\mathcal{A}}$, see \cite[Definition B.3]{Engel2000} for more information. To conclude,
\begin{equation} \label{eq:inclusionsA}
    \mathcal{D}(\mathcal{A}) \subseteq \mathcal{D}(\overline{\mathcal{A}}) \subseteq \mathcal{D}(\hat{\mathcal{A}}) \subseteq C_T(\mathbb{R},X),
\end{equation}
where the second inclusion follows from the fact that $\overline{\mathcal{A}}$ is the smallest closed linear extension of $\mathcal{A}$. Our aim is to show that $\mathcal{D}(\overline{\mathcal{A}}) = \mathcal{D}(\hat{\mathcal{A}})$ and that $\mathcal{D}(\mathcal{A})$ is dense in $C_T(\mathbb{R},X)$. The first claim characterizes $\mathcal{A}$ while the second claim will be of great independent importance in \cref{subsec: duality relations} when we will study the spectral structure of $\mathcal{A}$. Note that we can not directly conclude from \eqref{eq:Adecomp} and the density of $\mathcal{D}(M_A)$ and $\mathcal{D}(D)$ (\cref{lemma:MADclosed}) that $\mathcal{A}$ is densely defined. To prove the density, we will construct a dense subset $\mathcal{D}_\mathcal{A} \subseteq \mathcal{D}(\mathcal{A})$ of $C_T(\mathbb{R},X)$ such that $\hat{\mathcal{A}}$ is the closure of $\mathcal{A}$. In particular, the operator $\hat{\mathcal{A}}$ is given by $\mathcal{A}$ on the dense set $\mathcal{D}(\mathcal{A})$. We mention that our construction is closely related to \cite[Theorem 3.12]{Chicone1999} and \cite[Proposition 2.9]{Latushkin1996} but adapted to the $C_T(\mathbb{R},X)$ setting of periodic linear DDEs. To construct the set $\mathcal{D}_\mathcal{A}$, consider for any $s\in \mathbb{R}$ a function $ \alpha \in C_T^1(\mathbb{R},\mathbb{R})$ satisfying $\alpha(s+h) = 0$, an initial condition $\varphi \in \mathcal{D}(A(s))$ and introduce the map $f_{\alpha,s,\varphi} : \mathbb{R} \to X$ by
\begin{equation} \label{eq:falpha}
    f_{\alpha,s,\varphi}(\tau) := \alpha(\tau)U(\tau,s+lT)\varphi, \quad \tau \in [s+lT+h, s+(l+1)T+h), \quad l \in \mathbb{Z},
\end{equation}
and note that this function is well-defined as $A$ is $T$-periodic, see \eqref{eq:D(As)}. Recall that the number $h > 0$ appearing in \eqref{eq:falpha} is the finite delay from \cref{sec:periodic center manifolds}.

\begin{proposition} \label{prop:DAdense}
The set $\mathcal{D}_\mathcal{A} := \spn \{f_{\alpha,s,\varphi} : s \in \mathbb{R}, \varphi \in \mathcal{D}(A(s)), \alpha \in C_T^1(\mathbb{R},\mathbb{R}) \mbox{ with } \alpha(s+h) = 0 \} \subseteq \mathcal{D}(\mathcal{
A})$ is dense in $C_T(\mathbb{R},X)$ and the generator $\hat{\mathcal{A}}$ of the evolution semigroup $\mathcal{U}$ is the closure of $\mathcal{A}$.
\end{proposition}
\begin{proof}
Fix $f = f_{\alpha,s,\varphi}$ as defined in \eqref{eq:falpha}. Let us first prove that $f \in \mathcal{D}(\mathcal{A})$, i.e. $f \in C_T^1(\mathbb{R},X)$ and $f(\tau) \in \mathcal{D}(A(\tau))$ for all $\tau \in \mathbb{R}$. For a given $\tau \in \mathbb{R}$, there exists an integer $l$ such that $\tau \in [s+lT+h,s+(l+1)T+h)$ and one can verify easily from \eqref{eq:falpha} that $f(\tau+T) = f(\tau)$, i.e. $f$ is $T$-periodic. As $\tau \geq s + lT + h$, it follows from \cref{lemma:Dinvariant} that $f(\tau) \in \mathcal{D}(A(\tau))$ and there holds
\begin{equation} \label{eq:fdot}
    \dot{f}(\tau) = \dot{\alpha}(\tau)U(\tau,s+lT)\varphi + \alpha(\tau) A(\tau)U(\tau,s+lT)\varphi = \dot{\alpha}(\tau)U(\tau,s+lT)\varphi + A(\tau)f(\tau).
\end{equation}
Hence, $f \in C_T^1(\mathbb{R},X)$ since $\alpha(s+h) = 0$, the maps $\dot{\alpha},f$ and $\tau \mapsto A(\tau)$ are continuous, and $U$ is strongly continuous. Note that $(\mathcal{U}(t)f)(\tau) = \alpha(\tau-t)U(\tau,s+lT)\varphi$ and so
\begin{equation*}
    (\hat{\mathcal{A}}f)(\tau) = \frac{d}{dt}(\mathcal{U}(t)f)(\tau) \bigg|_{t = 0} = - \dot{\alpha}(\tau)U(\tau,s+lT)\varphi, \quad \tau \in [s + lT +h, s+(l+1)T+h), \quad l \in \mathbb{Z}.
\end{equation*}
Thus, it follows from \eqref{eq:fdot} that $(\hat{\mathcal{A}}f)(\tau) = A(\tau)f(\tau) -\dot{f}(\tau)$, i.e. $\mathcal{A}f = \hat{\mathcal{A}}f$.

To complete the proof, we have to show that $\mathcal{D}_\mathcal{A}$ is dense in $C_T(\mathbb{R},X) \cong C(\mathbb{T},X)$. Let us first observe that $\spn\{ \beta \otimes \varphi \in C(\mathbb{T},X) : \varphi \in X, \ \beta \in C(\mathbb{T},\mathbb{R}) \}$, where $(\beta \otimes \varphi)(\hat{\tau}) := \beta(\hat{\tau})\varphi$ for all $\hat{\tau} \in \mathbb{T}$, is dense in $C(\mathbb{T},X)$. Fix $g = \beta \otimes \varphi$ for some $\varphi \in X$ and $\beta \in C^1(\mathbb{T},\mathbb{R})$, and observe that $\| \beta \|_\infty < \infty$ as $\beta$ is continuous and $\mathbb{T}$ is compact. We will prove that $g$ can be approximated by a sum of functions from $\mathcal{D}_\mathcal{A}$, which can be identified as a subset of $C(\mathbb{T},X)$ as
\begin{equation} \label{eq:DAmod}
    \spn \{\hat{f}_{\alpha,s,\varphi} : s \in \mathbb{R}, \varphi \in \mathcal{D}(A(s)), \alpha \in C^1(\mathbb{T},\mathbb{R}) \mbox{ with } \alpha(\widehat{s+h}) = 0 \},
\end{equation}
where
\begin{equation*}
    \hat{f}_{\alpha,s,\varphi}(\hat{\tau}) := \alpha(\hat{\tau})U(\tau,s+lT)\varphi, \quad \tau \in [s + lT+h, s+(l+1)T+h), \quad l \in \mathbb{Z},
\end{equation*}
Fix $\varepsilon > 0$. Note that for every $\hat{\tau}_0 \in \mathbb{T}$ the map $(\tau,s) \mapsto U(\tau,s)\varphi$ is continuous at $(\tau_0,\tau_0)$. Hence, there are points $s_0 \leq s_0'$ such that $I(\hat{\tau}_0) = (\hat{s}_0,\hat{s}_0') \ni \hat{\tau}_0$ and $\|U(\tau,s_0)\varphi - \varphi \|_\infty \leq \varepsilon/(2\| \beta \|_\infty)$ for all $\tau \in (s_0,s_0')$. Hence, $\{I(\hat{\tau})\}_{\hat{\tau} \in \mathbb{T}}$ is an open covering of $\mathbb{T}$ and thus has a finite subcovering, denoted by $\{I_1,\dots,I_m\}$ for some $m \in \mathbb{N}$ due to compactness of $\mathbb{T}$. For each $I_j := (\hat{s}_j,\hat{s}_j')$ we have that $\|U(\tau,s_j)\varphi - \varphi \|_\infty \leq \varepsilon/(2\| \beta \|_\infty)$ for all $\hat{\tau} \in I_j$ with $j=1,\dots,m$. Choose a smooth partition of unity $\{\gamma_1,\dots,\gamma_m\}$ subordinate to the cover $\{I_1,\dots,I_m\}$. That is, $\gamma_j \in C^\infty(\mathbb{T},[0,1])$ satisfying $\sum_{j=1}^m \gamma_j(\hat{\tau}) = 1$ with $\hat{\tau} \in \mathbb{T}$ and $\mbox{supp } \gamma_j \subseteq I_j$ for all $j=1,\dots,m$. Recall from \cref{sec:periodic center manifolds} that for each $j=1,\dots,m$ the set $\mathcal{D}(A(s_j))$ is dense in $X$ and thus there exists a $\varphi_j \in \mathcal{D}(A(s_j))$ such that
\begin{equation*}
    \| \varphi - \varphi_j \|_\infty \leq \frac{\varepsilon}{2} M^{-1}e^{-\omega (T+h)},
\end{equation*}
where $M$ and $\omega$ are from the proof of \cref{lemma:C0semigroup}. Define the function $h : \mathbb{T} \to X$ by
\begin{equation*}
    h(\hat{\tau}) := \sum_{j=1}^m \beta(\hat{\tau})\gamma_j(\hat{\tau})U(\tau,s_j + lT)\varphi_j, \quad \tau \in [s_j + lT + h, s_j+(l+1)T+h), \quad l \in \mathbb{Z},
\end{equation*}
and since $\mbox{supp } \gamma_j \subseteq I_j$, the function $h$ is in $\mathcal{D}_\mathcal{A}$ as it is a sum of functions from \eqref{eq:DAmod}. Consider $\tau \in [s + lT + h, s+(l+1)T+h)$ for some $l \in \mathbb{Z}$, then
\begin{align*}
    \|g(\hat{\tau}) - h(\hat{\tau}) \|_\infty &\leq \bigg \| \beta(\hat{\tau}) \sum_{j=1}^m \gamma_j(\hat{\tau})  \varphi - \sum_{j=1}^m \beta(\hat{\tau}) \gamma_j(\hat{\tau}) U(\tau,s_j+lT)\varphi_j \bigg \|_\infty \\
    & \leq \|\beta\|_\infty \sum_{j=1}^m \gamma_j(\hat{\tau}) ( \| \varphi - U(\tau,s_j+lT)\varphi \|_\infty + \|U(\tau,s_j+lT)\| \|\varphi-\varphi_j\|_\infty) \leq \varepsilon.
\end{align*}
Taking the supremum over all $\tau \in \mathbb{R}$ yields $\|g-h\|_\infty \leq \varepsilon$, as required.
\end{proof}
Similar to \eqref{eq:Adecomp}, we can decompose the non-closed but closable unbounded linear operator $\mathcal{A}^{\odot \star}$ as 
\begin{equation*}
    \mathcal{D}(\mathcal{A}^{\odot \star}) = \mathcal{D}(M_{A^{\odot \star}}) \cap \mathcal{D}(D^{\odot \star}), \quad \mathcal{A}^{\odot \star} = M_{A^{\odot \star}} - D^{\odot \star},
\end{equation*}
where $M_{A^{\odot \star}} : \mathcal{D}(M_{A^{\odot \star}}) \to C_T(\mathbb{R},X^{\odot \star})$ is the (complexified) multiplication operator (associated to $A^{\odot \star}$) and $D^{\odot \star} : \mathcal{D}(D^{\odot \star}) \to C_T(\mathbb{R},X^{\odot \star})$ is the (complexified) differential operator. Both operators can be defined as in \eqref{eq:Arecht} and \eqref{eq:Drecht}, but on the sun-star level. One can introduce similarly as in \eqref{eq:evolutionsemi} the semigroup $\mathcal{U}^{\odot \star} := \{\mathcal{U}^{\odot \star}(s)\}_{s \geq 0}$ of bounded linear operators on $C_T(\mathbb{R},X^{\odot \star})$, where
\begin{equation*}
    (\mathcal{U}^{\odot \star}(s)f)(\tau) := U^{\odot \star}(\tau,\tau-s)f(\tau-s), \quad \forall f \in C_T(\mathbb{R},X^{\odot \star}), \ s \geq 0.
\end{equation*}
However, as the forward evolutionary system $U^{\odot \star}$ is not strongly continuous (\cref{sec:periodic center manifolds}), we can not expect that $\mathcal{U}^{\odot \star}$ forms a $\mathcal{C}_0$-semigroup on $C_T(\mathbb{R},X^{\odot \star})$.

\begin{remark} \label{remark:ODEclosed}
In the setting of finite-dimensional ODEs, \cref{thm:eigenfunctions} appeared for the first time in \cite[Lemma 2]{Iooss1988} and was later used extensively in \cite{Kuznetsov2005,Witte2013,Witte2014} to study codimension one and two bifurcations of limit cycles of
\begin{equation} \label{eq:ODEremark}
    \dot{x}(t) = f(x(t)),
\end{equation}
where $x(t) \in \mathbb{R}^n$ and the vector field $f : \mathbb{R}^n \to \mathbb{R}^n$ is sufficiently smooth. Here, the (complexified) unbounded linear operator $\mathcal{A} : \mathcal{D}(\mathcal{A}) \to C_T(\mathbb{R},\mathbb{C}^n)$ takes the form
\begin{equation*}
    \mathcal{D}(\mathcal{A}) = C_T^1(\mathbb{R},\mathbb{C}^n) \subseteq  C_T(\mathbb{R},\mathbb{C}^n), \quad \mathcal{A}\varphi = (\tau \mapsto A(\tau)\varphi(\tau) - \dot{\varphi}(\tau)),
\end{equation*}
where $A(\tau) = Df(\gamma(\tau)) \in \mathbb{R}^{n \times n}$, and $\gamma$ denotes a (nonhyperbolic) $T$-periodic solution of \eqref{eq:ODEremark}. We claim that $\mathcal{A}$ is closed. Let us first observe that the linear (complexified) multiplication operator $M_A$ from \eqref{eq:Arecht} in this setting is bounded since $\|M_A\varphi \|_\infty \leq  \| A \|_\infty \| \varphi \|_\infty$ as the real-valued function $\tau \mapsto \| A(\tau) \| $ is continuous, $T$-periodic and thus bounded. It follows from \eqref{eq:Adecomp} that $\mathcal{A}$ is closed since it is the sum of a bounded and closed linear operator. Since $\mathbb{R}^n$ is reflexive and finite-dimensional, the sun-star calculus construction becomes trivial \cite[Section 1]{Neerven1992} and thus $\mathcal{A}^{\odot \star} = \mathcal{A}$ is a closed linear operator. \hfill $\lozenge$
\end{remark}
\cref{prop:closable} and \cref{remark:ODEclosed} show that the underlying (spectral) structure of the operators $\mathcal{A}$ and $\mathcal{A}^{\odot \star}$ heavily depends on the periodic linear evolutionary system of interest. Indeed, in the setting of finite-dimensional ODEs, the operators $\mathcal{A}$ and $\mathcal{A}^{\odot \star}$ are closed while they are not closed but closable in the infinite-dimensional setting of DDEs. Even in the theory of periodic linear PDEs formulated on Bochner and Sobolev spaces, it turns out that $\mathcal{A}$ is frequently a closed unbounded linear operator, see \cite[Corollary 1]{Arendt2009} for more information. This loss of closedness in classical DDEs is a rather surprising result, since such a scenario can not occur when studying autonomous linear systems possessing a strongly continuous semigroup. Indeed, when a strongly continuous (semi)group is available, it is well-known that its associated (infinitesimal) generator, that generates the Jordan chains similar as in \cref{thm:eigenfunctions}, is always closed, see for example \cite[Theorem I.1.4]{Engel2000}.

\begin{remark} \label{remark:notation}
If we write $\mathcal{A}$ as $-\frac{d}{d\tau}+A(\tau)$ and $\mathcal{A}^{\odot \star}$ as $-\frac{d}{d\tau}+A^{\odot \star}(\tau)$, then \eqref{eq:ODEeig2} is equivalent to
\begin{equation*}
    \bigg( -\frac{d}{d\tau} + A(\tau) - \sigma \bigg)\varphi_i(\tau) =
    \begin{cases}
    0, \quad &i = 0,\\
    \varphi_{i-1}(\tau), \quad &i=1,\dots,m_\lambda-1,
    \end{cases}
\end{equation*}
and \eqref{eq:ODEeig} to
\begin{equation*}
    \bigg(-\frac{d}{d\tau} + A^{\odot \star}(\tau) - \sigma \bigg)j\varphi_i(\tau) =
    \begin{cases}
    0, \quad &i = 0,\\
    j\varphi_{i-1}(\tau), \quad &i=1,\dots,m_\lambda-1.
    \end{cases}   
\end{equation*}
This notation is also closely related to the weak$^\star$ formulation of the initial value problems from \eqref{eq:IVPzeta1} and \eqref{eq:IVPzetai}. Moreover, this notation is now in line with the literature \cite{Iooss1988,Iooss1999,Kuznetsov2005,Kuznetsov2023a,Witte2013,Witte2014} on the periodic smooth Jordan chains for finite-dimensional ODEs. Therefore, when we relate results from this paper with the literature, we will frequently use this notation. Strictly speaking, this notation  requires a special interpretation as $\frac{d}{d\tau}$ must act on functions in $C_T^1(\mathbb{R},X) \subseteq C_T(\mathbb{R},X)$ and $A(\tau)$ on functions in $\mathcal{D}(A(\tau)) \subseteq X$. We are confident that this notation would not create any confusion for the reader. \hfill $\lozenge$
\end{remark}

To represent in \cref{sec:characterization} the linear part of our local coordinate system on $\mathcal{W}_{\loc}^c(\Gamma)$ one can introduce the \emph{Floquet operator} (at time $\tau$) associated to the Floquet multiplier $\lambda$, defined as the coordinate map $Q_\lambda(\tau) : \mathbb{C}^{m_\lambda} \to E_\lambda(\tau)$ by
\begin{equation} \label{eq:Floquetmap}
    Q_\lambda(\tau)\xi:=\sum_{i=0}^{m_\lambda-1} \xi_i \varphi_i(\tau), \quad \forall \xi = (\xi_0,\dots,\xi_{m_\lambda-1}) \in \mathbb{C}^{m_\lambda}.
\end{equation}
It is clear from \Cref{thm:eigenfunctions} that the map $\tau \mapsto Q_\lambda(\tau)$ is $T$-periodic, $C^{k+1}$-smooth and takes values in $\mathcal{L}(\mathbb{C}^{m_\lambda}, E_\lambda(\tau))$ for all $\tau \in \mathbb{R}$. Using the notation from \cref{remark:notation}, a direct calculation shows that
\begin{equation} \label{eq:FloquetmapODE}
    \bigg( -\frac{d}{d\tau} +A^{\odot \star}(\tau)\bigg)j(Q_\lambda(\tau)\xi) = j(Q_\lambda(\tau)M_\lambda\xi),
\end{equation}
where $M_\lambda$ is the ${m_\lambda} \times m_\lambda$ Jordan matrix defined by
\begin{equation} \label{eq:matrixM}
    M_\lambda :=
    \begin{pmatrix}
        \sigma & 1 & \cdots & 0 \\
        0 & \sigma & \ddots & \vdots \\
        \vdots  & \ddots  & \ddots & 1  \\
        0 & \dots & 0 & \sigma
    \end{pmatrix}.
\end{equation}
This result is an extension of \cite[Proposition III.3]{Iooss1999} from finite-dimensional ODEs to classical DDEs. Because we are dealing with the real state space $X = C([-h,0],\mathbb{R}^n)$, the linear operator $M_\lambda : \mathbb{C}^{m_\lambda} \to \mathbb{C}^{m_\lambda}$ from \eqref{eq:matrixM} should represent a real operator, recall also \Cref{remark:complexification}. Depending on the location of the Floquet multiplier $\lambda$ in the complex plane, we have three options \cite{Iooss1999}:
\begin{itemize}
    \item If $\lambda \in (0,\infty)$, we choose $\sigma$ and $\varphi_0(\tau),\dots,\varphi_{m_\lambda-1}(\tau)$ real.
    \item If $\lambda \in \mathbb{C} \setminus \mathbb{R}$, then its complex conjugate $\overline{\lambda} \neq \lambda$ is also a Floquet multiplier. Hence, we choose $\sigma$ and $\varphi_0(\tau),\dots,\varphi_{m_\lambda-1}(\tau)$ complex and introduce $\overline{\sigma}$ and $\overline{\varphi_0}(\tau),\dots,\overline{\varphi_{m_\lambda-1}}(\tau)$ for the complex conjugate Jordan block.
    \item If $\lambda \in \mathbb{R}_{-}$, both methods do not succeed. Indeed, the Floquet exponents $\sigma$ are of the form $ \frac{\pi}{T} + \frac{2 \pi i}{T} l$ with $l \in \mathbb{Z}$. The standard way to deal with this situation is to double the period, since if $\lambda \in \sigma(U(s+T,s))$ is a Floquet multiplier then $\lambda^2 \in \sigma(U(s+2T,s))$.
\end{itemize}
To study this last case, one has to require $T$-antiperiodicity of the maps $\varphi_0,\dots,\varphi_{m_\lambda - 1}$ from \Cref{thm:eigenfunctions}. 

\begin{proposition} \label{prop:eigenfunctions2}
Let $\lambda \in \mathbb{R}_{-}$ be a Floquet multiplier of algebraic multiplicity $m_\lambda$ with $\sigma$ its associated Floquet exponent. Then there exist $T$-antiperiodic maps $\varphi_0,\dots,\varphi_{m_\lambda - 1} \in C^{k+1}(\mathbb{R},X)$ satisfying \eqref{eq:ODEeig2} and \eqref{eq:ODEeig}.
Furthermore, there exists a real $T$-periodic projector $P_\lambda: \mathbb{R} \to \mathcal{L}(X)$ such that $\mathcal{R}(P_\lambda(\tau)) = E_\lambda(\tau)$ for all $\tau \in \mathbb{R}$. Moreover, the $T$-antiperiodic Floquet operator $Q_\lambda$ satisfies \eqref{eq:FloquetmapODE}, where $M_\lambda$ represents a linear operator on $\mathbb{R}^{m_\lambda}$.
\end{proposition}
\begin{proof}
To prove the first assertion, we copy the proof of \Cref{thm:eigenfunctions} but in the $T$-antiperiodic setting. The proof goes identical up to \eqref{eq:zeta1}. If we set $\tau = s+T$ in \eqref{eq:zeta1} we get
\begin{equation*}
    |\lambda|\varphi_0(s+T) = e^{\sigma T} \varphi_0(s+T) = U(s+T,s)\varphi_0(s) = \lambda \varphi_0(s)
\end{equation*}
and so $\varphi_0(s + T) = \sign(\lambda) \varphi_0(s) = -\varphi_0(s)$, which proves the $T$-antiperiodicity of $\varphi_0$.

Consider now \eqref{eq:IVPzetai} and suppose that the right-hand side of the differential equation satisfies $\varphi_{i-1}(s+T) = -\varphi_{i-1}(s)$. Our goal is to find $\alpha_{ik}$ such that $\varphi_{i}$ is $T$-antiperiodic. Instead of requiring the $T$-periodicity of $\varphi_i$, we now require $T$-antiperiodicity of $\varphi_i$ and see that this holds if and only if
\begin{equation*}
    (U(s+T,s) - \lambda I)\varphi_i(s) = U(s+T,s) \sum_{k=1}^{i} \frac{(-1)^{k+1} T^k}{k!} \varphi_{i-k}(s) \notag \\ 
\end{equation*}
which is precisely \eqref{eq:thmU(s+T,s)-lambda}. Hence, the same procedure in \Cref{thm:eigenfunctions} can be followed to determine the associated $\alpha_{ik}$'s uniquely. The real spectral projection $P_\lambda(\tau) \in \mathcal{L}(X)$ onto $E_\lambda(\tau)$ for all $\tau \in \mathbb{R}$ is constructed in the same way as the Dunford integral in \cite[Appendix A.2]{Lentjes2023} and so $T$-periodicity follows. For the Floquet operator, it follows from linearity and $\varphi_i(\tau+T) = - \varphi_i(\tau)$ for all $i=0,\dots,m_\lambda-1$ that $Q_\lambda(\tau+T) = - Q_\lambda(\tau)$ for all $\tau \in \mathbb{R}$. 
\end{proof}
Hence, we conclude that the full construction obtained in this subsection also holds for Floquet multipliers $\lambda \in \mathbb{R}_{-}$ when one has to deal with $T$-antiperiodic (generalized) eigenfunctions. Note that all results on the evolution semigroup $\mathcal{U}$ still hold when formulated on a space of continuous $T$-antiperiodic functions, mainly due to \eqref{eq:Uperiodic}.

\subsection{Dual periodic smooth Jordan chains} \label{subsec: dual time periodic smooth Jordan}
As already announced in the \cref{sec:introduction}, we have to construct adjoint (generalized) eigenfunctions for our upcoming article on the study of bifurcations of limit cycles in DDEs. Therefore, our next goal is to repeat the construction from \cref{subsec: time periodic smooth Jordan}, but now for the dual system. 

Let us first observe that $\sigma(U(s+T,s)) = \sigma(U^{\star}(s-T,s))$ because $U(s+T,s)$ is a real operator (\Cref{remark:complexification}) and $U^{\star}(s,s+T) = U^\star(s-T,s)$ due to \eqref{eq:Uperiodic}. Since $U^\odot(s-T,s)$ is the restriction of $U^{\star}(s-T,s)$ to $X^\odot$, it follows from \cite[Proposition IV.2.17]{Engel2000}, that $ \sigma(U^{\star}(s-T,s)) = \sigma(U^{\odot}(s-T,s))$ and taking the adjoint once more yields the spectral equality
\begin{equation} \label{eq:spectraleq}
    \sigma(U(s+T,s)) = \sigma(U^\star(s-T,s)) = \sigma(U^\odot(s-T,s)) = \sigma(U^{\odot\star}(s+T,s)).
\end{equation}
Recall from \cref{sec:periodic center manifolds} that $E_\lambda(s) = \mathcal{N}((\lambda I - U(s+T,s))^{k_\lambda}))$, where $1 \leq k_\lambda < \infty$ is the order of the pole $\lambda = \mu$ of the map $\mu \mapsto (\mu I - U(s+T,s))^{-1}$. Moreover, $k_\lambda$ is the smallest integer such that $ \cup_{j \in \mathbb{N}} \mathcal{N}((\lambda I - U(s+T,s))^j) = \mathcal{N}((\lambda I - U(s+T,s))^{k_\lambda})$. Hence, $(\lambda I - U(s+T,s))^{k_\lambda} \in \mathcal{L}(X)$ is a projection and thus induces a direct sum decomposition of two closed $U(s+T,s)$-invariant subspaces
\begin{equation} \label{eq:directsumdecomp}
    X = E_\lambda(s) \oplus R_\lambda(s),
\end{equation}
where $R_\lambda(s) := \mathcal{R}((\lambda I - U(s+T,s))^{k_\lambda})$, see \cite[Section 3.6]{Lentjes2023} or \cite[Chapter 5]{Taylor1986} for more information. This induces a direct sum decomposition into two closed $U^\star(s-T,s)$-invariant subspaces $X^\star = E_\lambda^\star(s) \oplus R_\lambda^\star(s)$, where $E_\lambda^\star(s) = \mathcal{N}((\lambda I - U^\star(s-T,s))^{k_\lambda})$ and $R_\lambda^\star(s) = \mathcal{R}((\lambda I - U^\star(s-T,s))^{k_\lambda})$. Moreover, we obtain a similar decomposition for the spaces $X^\odot$ and $X^{\odot \star}$, see also \cite[Appendix A]{Lentjes2023}. The (adjoint) (generalized) eigenspaces corresponding to a Floquet multiplier $\lambda$ satisfy in addition 
\begin{equation} \label{eq:adjointeigenspaces}
    j(E_\lambda(s)) = E_\lambda^{\odot \star}(s), \quad E_\lambda^{\star}(s) = E_\lambda^{\odot}(s).
\end{equation}
Let us now prove that the obtained direct sum decompositions can be written in terms of annihilators.

\begin{lemma} \label{lemma:decomposition}
Let $\lambda$ be a Floquet multiplier of algebraic multiplicity $m_\lambda$ and suppose that $k_\lambda$ is the order of the pole $\lambda = \mu$ of $\mu \mapsto (\mu I - U(s+T,s))^{-1}$. Then we have the decompositions
\begin{equation} \label{eq:decomposition}
    X = E_\lambda(s) \oplus E_\lambda^\odot(s)^{\perp}, \quad X^{\odot} = E_\lambda^{\odot}(s) \oplus j(E_\lambda(s))^{\perp}, \quad \forall s \in \mathbb{R},\\
\end{equation}
where $\perp$ denotes the annihilator with respect to the natural dual pairing $\langle \cdot , \cdot \rangle$ from \eqref{eq:Xsunpairing} and \eqref{eq:Xsunstarpairing} in the first and second equality, respectively.
\end{lemma}
\begin{proof}
According to \eqref{eq:directsumdecomp} and the direct sum decomposition of $X^{\odot \star}$, it remains to show that $R_\lambda(s) = E_\lambda^{\odot}(s)^\perp$ and $R_\lambda^{\odot}(s) = E_\lambda(s)^\perp$. It follows immediately from Banach's closed range theorem \cite[Section VII.5]{Yosida1995} that $R_\lambda(s) = E_\lambda^{\star}(s)^\perp$ and $R_\lambda^{\odot}(s) = E_\lambda^{\odot \star}(s)^\perp$ since $X$ and $X^\odot$ are Banach spaces. The claim follows now from using \eqref{eq:adjointeigenspaces}.
\end{proof}
So far, the described construction above is similar as for autonomous DDEs, written in the language of $\mathcal{C}_0$-semigroups, see \cite[Section IV.2]{Diekmann1995} for more information. We are ready now to repeat the construction from \cref{sec:periodic spectral computations} but for the adjoint system. Recall from \cref{sec:periodic center manifolds} that the family $U^\star$ is not strongly continuous. Therefore, it is helpful to study instead some properties of the strongly continuous family $U^\odot$. The following result is a dual version of \cref{lemma:Dinvariant}.
\begin{lemma} \label{lemma:Dsuninvariant}
For any $s \in \mathbb{R}$, the operator $U^\odot(\tau,s)$ maps $X^\odot$ into $\mathcal{D}(A^\odot(\tau))$ for all $\tau \leq s-2h$. Moreover, there holds for any $\varphi^\odot \in X^\odot$ that
\begin{equation} \label{eq:UodotODE}
    \frac{\partial}{\partial \tau} U^\odot(\tau,s)\varphi^\odot = -A^\odot(\tau)U^\odot(\tau,s)\varphi^\odot, \quad \tau \leq s-2h,
\end{equation}
where the partial derivative should be interpreted in the norm topology of $X^\odot$.
\end{lemma}
\begin{proof}
Let $\varphi^\odot \in X^\odot$ and $s\in \mathbb{R}$ be given. It follows from \cite[Theorem 5.8]{Clement1988} that $\tau \mapsto \phi^\odot(\tau) := U^\odot(\tau,s)\varphi^\odot$ satisfies
\begin{equation} \label{eq:weaksol}
     \frac{d}{d\tau} \langle \phi^\odot(\tau),\varphi \rangle =  - \langle A^{\odot \star}(\tau)j\varphi, \phi^\odot(\tau) \rangle, \quad \forall \varphi \in j^{-1}\mathcal{D}(A_0^{\odot \star}), \quad \tau \leq s,
\end{equation}
where the map $\tau \mapsto \langle \phi^\odot(\tau), \varphi \rangle$ is continuously differentiable for all $\tau \leq s$. Because $U^\odot$ leaves $X^\odot$ invariant, we know that $\phi^\odot = (c,g)$ takes values in $X^\odot$. Recalling \eqref{eq:L(t)varphi} and \eqref{eq:Xsunstarpairing} yields
\begin{align*}
    \langle A^{\odot \star}(\tau)j\varphi, \phi^\odot(\tau) \rangle  &= \langle (L(\tau)\varphi,\dot{\varphi}),(c(\tau),g(\tau)) \rangle \\
    &= \langle c(\tau) \zeta(\tau,\cdot),\varphi \rangle + \int_0^h g(\tau)(\theta)\dot{\varphi}(-\theta) d\theta \\
    &= \int_0^h d_\theta[c(\tau) \zeta(\tau,\theta) + g(\tau)(\theta)] \varphi(-\theta),
\end{align*}
where we used integration by parts for Riemann-Stieltjes integrals and the conditions on $g$ coming from \eqref{eq:D(Astar0)} in the third equality. The $d_\theta$ refers to Riemann-Stieltjes integration over the $\theta$-variable. We can also express the left-hand side of \eqref{eq:weaksol} as
\begin{equation*}
    \frac{d}{d\tau} \langle \phi^\odot(\tau),\varphi \rangle =  \dot{c}(\tau) \varphi(0) + \int_0^h \frac{\partial g(\tau)(\theta)}{\partial \tau} \varphi(-\theta) d\theta = \int_0^h d_\theta\bigg[\frac{\partial \phi^\odot(\tau)(\theta)}{\partial \tau} \bigg] \varphi(-\theta).
\end{equation*}
Hence, \eqref{eq:weaksol} is equivalent to
\begin{equation*}
    \int_0^h d_\theta \bigg[c(\tau)\zeta(\tau,\theta) + \bigg(\frac{\partial}{\partial \tau} + \frac{\partial}{\partial \theta} \bigg)  \phi^\odot(\tau)(\theta) \bigg] \varphi(-\theta) = 0, \quad \forall \varphi \in j^{-1}\mathcal{D}(A_0^{\odot \star}),
\end{equation*}
since $g(\tau)$ represents the derivative of $\phi^\odot(\tau)$ with respect to the state. Clearly, if we can show that $\phi^\odot$ satisfies
\begin{equation} \label{eq:PDEadjoint}
    \bigg( \frac{\partial}{\partial \tau} + \frac{\partial}{\partial \theta} \bigg)\phi^\odot(\tau)(\theta) + c(\tau)\zeta(\tau,\theta) = 0, \quad \forall \tau \leq s, \ \theta \in (0,h],
\end{equation}
then \eqref{eq:weaksol} is satisfied. The (inhomogeneous) transport equation \eqref{eq:PDEadjoint} has the unique solution
\begin{equation} \label{eq:phi_isun}
    \phi^\odot(\tau)(\theta) = c(\tau-\theta) - \int_0^\theta c(\tau + v - \theta) \zeta(\tau + v- \theta,v) dv,
\end{equation}
when the initial condition $c$ is specified. To determine $c$ from $\varphi^\odot = \phi^\odot(s)$, let us take a look at the map $\theta \mapsto g(\tau)(\theta) = \frac{\partial}{\partial \theta} \phi^\odot(\tau)(\theta)$ that has to satisfy $g(\tau)(h) = 0$ for all $\tau \leq s$. Differentiating \eqref{eq:phi_isun} with respect to $\theta$ by employing the Leibniz integral rule and integration by parts for Riemann-Stieltjes integrals leads eventually to
\begin{equation*}
    g(\tau)(\theta) = - \dot{c}(\tau - \theta) - \int_0^\theta c(\tau + v - \theta) d_2 \zeta(\tau + v - \theta,v).
\end{equation*}
To simplify the expression even more, recall that $g(\tau - \theta + h)(h) = 0$ for all $\tau \in \mathbb{R}$ and $\theta \in (0,h]$ that satisfy $\tau - \theta + h \leq s$. This is equivalent to
\begin{equation} \label{eq:g2}
    0 = - \dot{c}(\tau - \theta) - \int_0^h c(\tau + v - \theta) d_2 \zeta(\tau + v - \theta,v),
\end{equation}
and so $g$ simplifies to
\begin{equation} \label{eq:g3}
    g(\tau)(\theta) = \int_\theta^h c(\tau - \theta + v)d_2\zeta(\tau - \theta + v,v).
\end{equation}
As $g(\tau + h)(h) = 0$ for all $\tau + h \leq s$, it follows from \eqref{eq:g2} that $c$ satisfies
\begin{equation} \label{eq:cDDE}
    \dot{c}(\tau) = - \int_0^h c(\tau + v) d_2 \zeta(\tau + v,v), \quad \forall \tau \leq s-h,
\end{equation}
which is an advance differential equation that has to be solved backward in time with an initial condition that still has to be specified. To determine the initial condition, recall that $\varphi^\odot$ is known and so it follows from \eqref{eq:phi_isun} that $c$ must satisfy
\begin{equation*}
     c(s-\theta) = \varphi^\odot(\theta) + \int_0^\theta c(s + v - \theta) \zeta(s + v - \theta,v) dv, \quad \forall \theta \in (0,h].
\end{equation*}
Perform the change of variables: $s-\theta = -\xi, z(\xi) = c(-\xi), \psi(\xi) = \varphi^\odot(s+\xi)$ and define the map $K$ as $K(u-w) := K(w-u,w)$. This yields the linear Volterra integral equation of the second kind
\begin{equation} \label{eq:z}
    z(\xi) = \psi(\xi) + \int_{-s}^\xi z(\eta) K(\xi-\eta) d\eta, \quad \forall \xi \in (-s,-s+h].
\end{equation}
Because $\psi \in \AC((-s,-s+h],\mathbb{C}^{n \star})$ and $K \in L^1((0,h],\mathbb{C}^{n \star})$ as the kernel $\zeta$ is $C^{k}$-smooth in the first component and of bounded variation of the second component, it follows from \cite[Theorem 2.3.5]{Gripenberg1990} that \eqref{eq:z} has a unique solution $\tilde{z} = z(-s + h + \cdot) \in \AC([-h,0],\mathbb{C}^{n \star})$, where we define $\tilde{z}(0) := \lim_{t \downarrow 0 } z(-s+ h + t)$. Perform now the change of variables $\tau = -t, z(t) = c(-t)$ in  \eqref{eq:cDDE} and define the kernel $\tilde{\zeta}$ by $\tilde{\zeta}(t,v) := \zeta(v-t,v)$ to obtain
\begin{equation} \label{eq:zDDE}
    \begin{dcases}
       \dot{z}(t) = \int_0^h z(t-v) d_2 \tilde{\zeta}(t,v), \quad &t \geq -s+h, \\
       z(-s+h+\cdot) = \tilde{z}, \quad &\tilde{z} \in C([-h,0],\mathbb{C}^{n \star}),
    \end{dcases}
\end{equation}
which is a periodic linear DDE, since $t \mapsto \tilde{\zeta}(t,\cdot)$ is $T$-periodic, with initial condition $\tilde{z}$ defined on starting time $-s+h$. Recall from the methods of steps (see for example the proof of \cref{lemma:Dinvariant}) that $z|_{[-s+h,\infty)}$ is $C^1$-smooth, and so $c|_{(-\infty,s-h]}$ is $C^1$-smooth. Hence, if we take $\tau \leq s-2h$, then $[\tau,\tau-h] \subset (-\infty,s-h]$ and so $c|_{[\tau,\tau-h]}$ is $C^1$-smooth. We claim that $\tau \mapsto g(\tau)$ is $C^1$-smooth for $\tau \leq s-2h$. Notice from \eqref{eq:g3} that
\begin{equation} \label{eq:gdot}
    \dot{g}(\tau)(\theta) = -\int_\theta^h \dot{c}(\tau-\theta+v) d_2\zeta(\tau+v,v) - \int_\theta^h c(\tau-\theta+v) d_2 D_1\zeta(\tau+v,v),
\end{equation}
and so $\tau \mapsto g(\tau)$ is $C^1$-smooth for $\tau \leq s-2h$ since $\tau \mapsto \dot{c}(\tau)$ is continuous for $\tau \leq s-2h$ and $\tau \mapsto D^1\zeta(\tau,\cdot)$ is $C^{k-1}$ smooth. Hence, we conclude that $\tau \mapsto \phi^\odot(\tau) = (c(\tau),g(\tau))$ is $C^1$-smooth for $\tau \leq s-2h$. Moreover, it follows from \eqref{eq:g3} that the map $\theta \mapsto g(\tau)(\theta)$ has bounded variation for all $\theta \in (0,h]$. As $g(\cdot)(0) = 0$, we have shown that $g(\tau) \in \NBV([0,h],\mathbb{C}^{n \star})$ and so $\phi_i^\odot(\tau) \in \mathcal{D}(A_0^\star)$ for all $\tau \leq s-2h$. According to \eqref{eq:PDEadjoint}, we see that $c(\tau) \zeta(\tau,\cdot) + g(\tau) = -\dot{\phi}^{\odot}(\tau) \in X^\odot$ for $\tau \leq s-2h$ due to \eqref{eq:cDDE} and \eqref{eq:gdot} as $c(\tau) \in \mathbb{C}^{n\star}, \dot{g}(\tau) \in L^1([0,h],\mathbb{C}^{n \star})$ and $\dot{g}(\tau)(h) = 0$ for all $\tau \leq s-2h$. To conclude, we have proven that $\tau \mapsto (c(\tau),g(\tau)) = U^\odot(\tau,s)\varphi^\odot$ is in $\mathcal{D}(A^\odot(\tau))$ for all $\tau \leq s-2h$.

Since $U^\odot$ is a strongly continuous evolutionary system, one can use the same strategy as performed in the proof of \cref{lemma:Dinvariant} to prove \eqref{eq:UodotODE}.
\end{proof}

Following the structure of \cref{subsec: time periodic smooth Jordan}, let us construct now a Jordan chain of $U^{\odot}(s-T,s)$ at a given Floquet multiplier $\lambda \in \sigma(U^\odot(s-T,s))$. Again, by the construction given in \cite[Section IV.4]{Diekmann1995}, there exists an ordered basis $\{\phi_{m_\lambda - 1}^\odot(s),\dots,\phi_{0}^\odot(s)\}$ of $E_\lambda^\odot(s)$, called a \emph{Jordan chain} (of $E_\lambda^\odot(s)$), satisfying
\begin{equation} \label{eq:jordan chain sun}
    (U^\odot(s-T,s) - \lambda I)\phi_i^\odot(s) =
    \begin{dcases}
    0, \quad &i = m_\lambda - 1, \\
    \phi_{i+1}^\odot(s), \quad &i=m_\lambda-2,\dots,0.
    \end{dcases}
\end{equation}
As $U^\odot$ leaves $X^\odot$ invariant, the Jordan chain in \eqref{eq:jordan chain sun} is well-defined. Recall from \Cref{subsec: time periodic smooth Jordan} that the bounded linear operator $U_\lambda(\tau,s) : E_\lambda(s) \to E_\lambda(\tau)$ is a topological isomorphism. Hence, $U_\lambda^{\odot}(\tau,s) := (U_\lambda(s,\tau))^\odot : E_\lambda^\odot(s) \to E_\lambda^\odot(\tau)$ is a topological isomorphism and so $\{\phi_{m_\lambda - 1}^\odot(\tau),\dots,\phi_{0}^\odot(\tau) \}$ is an ordered basis of $E_\lambda^\odot(\tau)$, where $\phi_{i}^\odot(\tau) := U_\lambda^\odot(\tau,s)\phi_i^\odot(s)$ for all $\tau \in \mathbb{R}$ and $i=m_\lambda - 1,\dots,0$. The following lemma, a dual version of \cref{lemma:basisD(A)}, shows that $E_\lambda^\odot(\tau)$ has additional structure.

\begin{lemma} \label{lemma:adjoint D(Astar)}
The ordered basis $\{\phi_{m_\lambda - 1}^\odot(\tau),\dots,\phi_{0}^\odot(\tau) \} \subseteq \mathcal{D}(A^\odot(\tau))$ consists of $C^{k+1}$-smooth functions and forms a Jordan chain for $E_\lambda^\odot(\tau)$ for all $\tau \in \mathbb{R}$.
\end{lemma}
\begin{proof}
The proof of the Jordan chain structure is analogous to that of \Cref{lemma:basisD(A)}. According to \cref{lemma:Dsuninvariant}, we have that $\phi_{i}^\odot(\tau) = U_\lambda^\odot(\tau,s)\phi_i^\odot(s)$ is in $\mathcal{D}(A^\odot(\tau))$ for all $\tau \leq s-2h$. As $\{\phi_{m_\lambda - 1}^\odot(\tau),\dots,\phi_{0}^\odot(\tau) \}$ forms a Jordan chain, one can use the same technique as in \cref{lemma:basisD(A)} to prove that $\phi_{i}^\odot(\tau)$ is in $\mathcal{D}(A^\odot(\tau))$ for all $i=m_\lambda-1,\dots,1$ and $\tau \in \mathbb{R}$. Note that each $\phi_{i}^\odot = (c_i,g_i)$ has an associated $z_i$ that solves the periodic linear DDE \eqref{eq:zDDE}. Using again the methods of step on this DDE and the Jordan chain structure as performed in \Cref{lemma:basisD(A)}, we obtain that $z_i$ is $C^{k+1}$-smooth and thus $\tau \mapsto \phi_{i}^\odot(\tau) = (c_i(\tau),g_i(\tau))$ is $C^{k+1}$-smooth due to the explicit representations \eqref{eq:g3} and \eqref{eq:cDDE}.
\end{proof}

Let us now take a look at the $T$-(anti)periodicity of the Jordan chain for the adjoint system. It is clear from the computation
\begin{equation*}
    \phi_i^\odot(s-T) - \phi_i^\odot(s)  = 
    \begin{dcases}
    (\lambda - 1) \phi_{m_\lambda - 1}^\odot(s), \quad &i= m_\lambda - 1, \\
    (\lambda - 1) \phi_i^\odot(s) + \phi_{i+1}^\odot(s), \quad &i=m_\lambda-2,\dots,0,
    \end{dcases}
\end{equation*}
that $\tau \mapsto \phi_i^\odot(\tau)$ is $T$-periodic if and only if $\lambda = 1$ and $i=m_\lambda - 1$, and $T$-antiperiodic if $\lambda = -1$ and $i=m_\lambda - 1$. However, similarly as for the (generalized) eigenfunctions of $\mathcal{A}$ from \cref{thm:eigenfunctions}, we need a $C^{k+1}$-smooth $T$-(anti)periodic basis of $E_\lambda^\odot(\tau)$ in our upcoming paper on the study of bifurcations of limit cycles in classical DDEs. Therefore, let us introduce the unbounded linear operators $\mathcal{A}^{\odot} : \mathcal{D}(\mathcal{A}^{\odot}) \to C_T(\mathbb{R},X^{\odot})$ and $\mathcal{A}^\star : \mathcal{D}(\mathcal{A}^\star)  \to C_T(\mathbb{R},X^\star)$ by
\begin{align} \label{eq:D(Asun)DDE}
\begin{split} 
    \mathcal{D}(\mathcal{A}^{\odot}) &:= \{ \varphi^\odot \in C_T^1(\mathbb{R},X^{\odot}) : \varphi^\odot(\tau) \in \mathcal{D}(A^{\odot}(\tau)) \mbox{ for all $\tau \in \mathbb{R}$} \} \subseteq C_T(\mathbb{R},X^{\odot}), \\
    \mathcal{D}(\mathcal{A}^\star) &:= \{ \varphi^\star \in C_T^1(\mathbb{R},X^\star) : \varphi^\star(\tau) \in \mathcal{D}(A^\star(\tau)) \mbox{ for all $\tau \in \mathbb{R}$} \} \subseteq C_T(\mathbb{R},X^\star),
\end{split}
\end{align}
with action
\begin{equation*}
    \mathcal{A}^{\odot} \varphi^\odot := ( \tau \mapsto A^{\odot}(\tau)\varphi^\odot(\tau) + \dot{\varphi}^\odot(\tau)), \quad \mathcal{A}^\star\varphi^\star := ( \tau \mapsto A^\star(\tau)\varphi^\star(\tau) + \dot{\varphi}^\star(\tau)).
\end{equation*}
The $T$-antiperiodic setting will be discussed at the end of this subsection. The following result illustrates that the periodic (generalized) eigenfunctions of $\mathcal{A}^\star$ and $\mathcal{A}^\odot$ can be obtained similarly as the periodic (generalized) eigenfunctions of $\mathcal{A}$ and $\mathcal{A}^{\odot \star}$ from \cref{thm:eigenfunctions}.

\begin{theorem} \label{thm:adjoint eigenfunctions}
Let $\lambda \in \mathbb{C} \setminus \mathbb{R}_{-}$ be a Floquet multiplier of algebraic multiplicity $m_\lambda$ with $\sigma$ its associated Floquet exponent. Then there exist $\varphi_{i}^\odot \in C_T^{k+1}(\mathbb{R},X^\odot)$ satisfying
\begin{equation} \label{eq:ODE adjoint2}
    (\mathcal{A}^\odot - \sigma I)\varphi_i^\odot=
    \begin{cases}
    0, \quad &i = m_\lambda-1,\\
    \varphi_{i+1}^\odot, \quad &i=m_\lambda-2,\dots,0,
    \end{cases}
\end{equation}
or equivalently
\begin{equation} \label{eq:ODE adjoint}
    (\mathcal{A}^\star - \sigma I)\varphi_i^\odot=
    \begin{cases}
    0, \quad &i = m_\lambda-1,\\
    \varphi_{i+1}^\odot, \quad &i=m_\lambda-2,\dots,0.
    \end{cases}
\end{equation}
such that the set of functions $\{\varphi_{m_\lambda - 1}^\odot(\tau),\dots,\varphi_{0}^\odot(\tau) \}$ is an ordered basis of $E_\lambda^\odot(\tau)$ for all $\tau \in \mathbb{R}$.
\end{theorem}
\begin{proof}
Let $s \in \mathbb{R}$ be a starting time and consider the ordered basis $\{\phi_{m_\lambda - 1}^\odot(s),\dots,\phi_{0}^\odot(s)\}$ of $E_\lambda^\odot(s)$ from \eqref{eq:jordan chain sun}. We show the claim by induction on $i \in \{m_\lambda-1,\dots, 0 \}$. For the base case ($i=m_\lambda-1)$, consider the initial value problem
\begin{equation} \label{eq:IVPpsi1}
    \begin{dcases}
    ( d^\star + A^\star(\tau) - \sigma)\varphi_{m_\lambda-1}^\odot(\tau) = 0, \quad \tau \leq s,\\
    \varphi_{m_\lambda-1}^\odot(s) = \phi_{m_\lambda-1}^\odot(s),
    \end{dcases}
\end{equation}
where $\phi_{m_\lambda-1}^\odot(s)$ is the first basis vector of $E_\lambda^\odot(s)$. It follows from \eqref{eq:IVPpsi1} that
\begin{align*}
    d^\star (e^{\sigma(s-\cdot)}\varphi_{m_\lambda-1}^\odot)(\tau) &= -\sigma e^{\sigma(s-\tau)}\varphi_{m_\lambda-1}^\odot(\tau) + e^{\sigma(s-\tau)} d^\star \varphi_{m_\lambda-1}^\odot(\tau) \\
    &= e^{\sigma(s-\tau)}(d^\star -\sigma)\varphi_{m_\lambda-1}^\odot(\tau) = -A^\star(\tau)[e^{\sigma(s-\tau)}\varphi_{m_\lambda-1}^\odot(\tau)].
\end{align*}
This differential equation is of the form \cite[Equation (5.8)]{Clement1988} and so \eqref{eq:IVPpsi1} admits a unique solution \cite[Theorem 5.8]{Clement1988} on $(-\infty,s]$ given by
\begin{equation} \label{eq:adjoint_eigUsun}
    e^{\sigma(s-\tau)}\varphi_{m_\lambda-1}^\odot(\tau) = U^{\odot}(\tau,s)\varphi_{m_\lambda-1}^\odot(s), \quad \forall \tau \leq s,
\end{equation}
whenever $\varphi_{m_\lambda-1}^\odot(s) \in \mathcal{D}(A_0^\star)$. Since $\varphi_{m_\lambda-1}^\odot(s) = \phi_{m_\lambda-1}^\odot(s)$, the claim follows from \Cref{lemma:adjoint D(Astar)} because $\phi_{m_\lambda-1}^\odot(s) \in \mathcal{D}(A^{\odot}(s)) \subseteq \mathcal{D}(A_0^{\star})$, recall \cref{lemma:Astar}. $T$-periodicity of $\varphi_{m_\lambda-1}^\odot$ is proven as in \eqref{eq:Tperiodiceig} and so this map extends to $\mathbb{R}$. Hence, $\varphi_{m_\lambda-1}^\odot(\tau) = e^{-\sigma(s-\tau)}\phi_{m_\lambda - 1}^\odot(\tau)$ for all $\tau \in \mathbb{R}$, and so it follows from \cref{lemma:adjoint D(Astar)} that $\varphi_{m_\lambda-1}^\odot$ is $C^{k+1}$-smooth and takes values in $\mathcal{D}(A^{\odot}(\tau)) \subseteq X^{\odot}$, i.e. $\varphi_{m_\lambda-1}^\odot \in \mathcal{D}(\mathcal{A}^\odot) \subseteq \mathcal{D}(\mathcal{A}^\star)$. This proves the base case for \eqref{eq:ODE adjoint} and \eqref{eq:ODE adjoint2}.

Assume now that the maps $\varphi_{m_\lambda - 1}^\odot,\dots,\varphi_{i+1}^\odot \in C_T^{k+1}(\mathbb{R},X^\odot)$ are constructed for some $ i \in \{m_\lambda - 2,\dots,0 \}$ and consider the initial value problem
\begin{equation} \label{eq:IVPpsii}
    \begin{dcases}
    ( d^\star + A^\star(\tau) - \sigma)\varphi_i^\odot(\tau) = \varphi_{i+1}^\odot(\tau), \quad \tau \leq s, \\
    \varphi_i^\odot(s) = \sum_{k = i}^{m_\lambda - 1} \alpha_{ik} \phi_k^\odot(s).    \end{dcases}
\end{equation}
The goal is to find scalars $\alpha_{ik}$ such that $\varphi_i^\odot$ becomes $T$-periodic. Notice from \eqref{eq:IVPpsii} that
\begin{equation} \label{eq:ODEvarphiodot}
    d^\star (e^{\sigma(s-\cdot)} \varphi_i^\odot)(\tau) =  e^{\sigma(s-\tau)} ( d^\star  - \sigma) \varphi_i^\odot(\tau) = -A^\star(\tau)[e^{\sigma(s-\tau)}\varphi_i^\odot(\tau)] + e^{\sigma (s-\tau)}\varphi_{i+1}^\odot(\tau).
\end{equation}
Hence, it suffices to prove that \eqref{eq:ODEvarphiodot} together with the initial condition specified in \eqref{eq:IVPpsii} admits a unique solution on $(-\infty,s]$. Consider the function $w_i^\odot : (-\infty,s] \to X^{\odot}$ defined by
\begin{equation} \label{eq:varphi_i^sun_uitdrukking}
    w_i^\odot(\tau) := U^\odot(\tau,s) \sum_{l=i}^{m_\lambda - 1} \frac{(\tau-s)^{l-i}}{(l-i)!}\varphi_l^\odot(s), \quad \forall \tau \in (-\infty,s].
\end{equation}
Since $\phi_{m_\lambda - 1}^\odot(s),\dots,\phi_{i}^\odot(s) \in \mathcal{D}(A^\odot(s))$, we get from \eqref{eq:IVPpsii} that $\varphi_{m_\lambda - 1}^\odot(s),\dots,\varphi_{i}^\odot(s) \in \mathcal{D}(A^\odot(s))$. It is clear that
\begin{equation*}
    U^\odot(\tau,s) \varphi_j^\odot(s) = \sum_{k=j}^{m_\lambda - 1} \alpha_{ik} \phi_k^\odot(\tau) \in \mathcal{D}(A^\odot(\tau)), \quad \forall \tau \in (-\infty,s], \ j \in \{m_\lambda-1,\dots,i\},
\end{equation*}
and so it follows from \Cref{lemma:adjoint D(Astar)} that $\tau \mapsto w_i^\odot(\tau)$ takes values in $\mathcal{D}(A^\odot(\tau)) $ and is $C^{k+1}$-smooth, which implies weak$^\star$ differentiability of $w_i^\odot$. Clearly $w_i^\odot(s) = \varphi_i^\odot(s)$ and notice that
\begin{align*}
    d^\star w_i^\odot(\tau) &= -A^\star(\tau) U^\odot(\tau,s) \sum_{l=i}^{m_\lambda - 1} \frac{(\tau-s)^{l-i}}{(l-i)!}\varphi_l^\odot(s) + U^{\odot}(\tau,s)\sum_{l = i+1}^{m_\lambda - 1} \frac{(\tau - s)^{l-i-1}}{(l-i-1)!}\varphi_l^\odot(s) \\
    &= - A^\star(\tau) w_i^\odot(\tau) + w_{i+1}^\odot(\tau),
\end{align*}
and so $w_i^\odot$ is a solution on $(-\infty,s]$ of \eqref{eq:ODEvarphiodot}. Since $w_{i+1}^\odot$ is at least continuous, it follows from \Cref{prop:solutionadjoint} and by the construction of $w_i^\odot$ that \eqref{eq:ODEvarphiodot} admits a unique solution $w_i^\odot$ on $(-\infty,s]$ where $w_i^\odot = e^{\sigma(s-\cdot)} \varphi_i^\odot$. As a consequence, $\varphi_i^\odot = e^{-\sigma(s-\cdot)}w_i^{\odot}$ is $C^{k+1}$-smooth. To prove the $T$-periodicity and linear independence at time $\tau \leq s$ of the adjoint (generalized) eigenfunctions, one uses the same technique as the proof on the $T$-periodicity and linear independence of the (generalized) eigenfunctions from \cref{thm:eigenfunctions}. We conclude that $\varphi_i^\odot \in \mathcal{D}(\mathcal{A}^\odot) \subseteq \mathcal{D}(\mathcal{A}^{\star})$, as required.
\end{proof}
\cref{thm:adjoint eigenfunctions} shows us that $\{ \sigma \in \mathbb{C} : \sigma \mbox{ is a Floquet exponent} \}$ is a subset of the point spectrum of the unbounded linear operators $\mathcal{A}^\odot$ and $\mathcal{A}^{\star}$. To study the spectral structure of these operators in \cref{subsec: duality relations}, we have to prove similar to \cref{prop:closable} that $\mathcal{A}^{\odot}$ and $\mathcal{A}^\star$ are at least closable linear operators.
\begin{proposition} \label{prop:closable2}
    The unbounded linear operators $\mathcal{A}^{\odot}$ and $\mathcal{A}^\star$ are not closed but closable.
\end{proposition}
\begin{proof}
Let us first prove that $\mathcal{A}^\odot$ is closable. Let $(\varphi_m^{\odot})_m$ be a sequence in $\mathcal{D}(\mathcal{A}^\odot)$ converging in norm to zero, and assume that the sequence $(\mathcal{A}^\odot \varphi_m^\odot)_m$ converges in norm to some $\psi^\odot \in C_T(\mathbb{R},X^{\odot})$. Set $\psi_m^\odot = \mathcal{A}^\odot \varphi_m^\odot$ and notice that this is equivalent to
\begin{equation} \label{eq:Asunclosable}
    \bigg(\frac{\partial}{\partial \theta} + \frac{\partial}{\partial \tau} \bigg) \varphi_m^\odot(\tau)(\theta) = -\varphi_m^\odot(\tau)(0^+) \zeta(\tau,\theta) + \psi_m^\odot(\tau)(\theta),  \quad \forall \tau \in \mathbb{R}, \ \theta \in (0,h].
\end{equation}
Analogous to the computations in the proof of \cref{lemma:Dsuninvariant}, this (inhomogeneous) transport equation has a solution $ \varphi_m^\odot = (c_m,g_m)$ given by
\begin{align} \label{eq:gncnAsun}
\begin{split}
    g_m(\tau)(\theta) &= \int_\theta^h c_m(\tau-\theta+s)d_2\zeta(\tau-\theta+s,s) - \int_\theta^h k_m(\tau-\theta+s)(s) ds, \\ 
    \dot{c}_m(\tau) &= - \int_0^h c_m(\tau+s) d_2\zeta(\tau+s,s) + d_m(\tau) + \int_0^h k_m(\tau+s)(s) ds,
\end{split}
\end{align}
where $\psi_m^\odot = (d_m,k_m)$. Since $\varphi_m^\odot \to 0$ in norm, we have that $c_m \to 0$ uniformly, and $g_m \to 0$ uniformly in the first component and in $L^1$ in the second component as $m \to \infty$. Hence, it follows from the first equation in \eqref{eq:gncnAsun} that
\begin{equation*}
    \int_\theta^h k_m(\tau-\theta+s)(s) ds \to 0, \quad m \to \infty, \quad \forall \tau \in \mathbb{R}, \ \theta \in (0,h].
\end{equation*}
Similar to the proof of \cref{prop:closable}, one obtains
\begin{equation} \label{eq:integralkn}
    \int_{\theta_1}^{\theta_2} k_m(\tau-\theta_1+s)(s) ds \to 0, \quad m \to \infty, \quad \forall \tau \in \mathbb{R}, \ \theta_1,\theta_2 \in (0,h].
\end{equation}
Since $\psi_m^\odot \to \psi^\odot$, we have that $d_m \to d$ and $k_m \to k$ as $m \to \infty$ for some $d \in C_T(\mathbb{R},\mathbb{C}^{n\star})$ and $k \in C_T(\mathbb{R},L^1([0,h],\mathbb{C}^{n\star}))$. The same strategy as in the proof of \cref{prop:closable} can be applied to \eqref{eq:integralkn} to conclude that $k = 0$ almost everywhere. The only minor difference is by the $L^1$-convergence of $(r_m)_m$, now defined by $s \mapsto k_m(\tau - \theta_1 + s)(s)$, that one has to (possibly) extract a subsequence $(r_{m_j})_{m_j}$ of $(r_m)_m$ that converges almost everywhere pointwise towards $r$, defined by $s \mapsto k(\tau - \theta_1 + s)(s)$, as $j \to \infty$. It remains to show that $d = 0$. Because $c_m \to 0, d_m \to d$ and $k_m \to k$ in norm as $m \to \infty$, we deduce from the second equation in \eqref{eq:gncnAsun} that $\dot{c}_m \to \dot{c}$ in uniformly as $m \to \infty$. Hence, it follows from taking the limit as $m \to \infty$ from the second equation of \eqref{eq:gncnAsun} that $\dot{c} = d$ and so $c(\tau) = c(s) + \int_s^\tau d(t) dt$ for all $\tau \in \mathbb{R}$. But $c$ is the zero function and so $d=0$ by continuity of $d$. Thus $\psi^\odot = 0$ almost everywhere, which proves that $\mathcal{A}^\odot$ is closable. To prove closability of $\mathcal{A}^\star$, notice that the sequence of functions $(g_m)_m$ and $(k_m)_m$ also belong to $\NBV([0,h],\mathbb{C}^{n\star})$ due to \eqref{eq:D(Astar0)} and \cref{lemma:adjoint D(Astar)}. Since $\NBV([0,h],\mathbb{R}^{\star}) \subseteq L^1([0,h],\mathbb{C}^{n\star})$, the same strategy can be applied as above to prove that $\mathcal{A}^\star$ is closable.

Let us now prove that $\mathcal{A}^{\odot}$ is not closed. Let $(f_m)_m$ be a sequence in $C_T^1(\mathbb{R},\mathbb{C}^{n\star})$ converging uniformly to some $f \in C_T(\mathbb{R},\mathbb{C}^{n\star}) \setminus C_T^1(\mathbb{R},\mathbb{C}^{n\star})$. Introduce the sequences $(F_m)_m$ and $(\varphi_m^{\odot})_m$ by $F_m(\tau,\theta) = \int_0^\theta f_m(\tau-s) ds$ and $\varphi_m^\odot(\tau)(\theta) = (\theta-h)F_m(\tau,\theta) - \int_0^\theta F_m(\tau,s) ds$, and introduce the maps $F$ and $\varphi^\odot$ by $F(\tau,\theta) = \int_0^\theta f(\tau-s) ds$ and $\varphi^\odot(\tau)(\theta) = (\theta-h)F(\tau,\theta) - \int_0^\theta F(\tau,s) ds$, for all $\theta \in (0,h]$ and set $\varphi_m^\odot(\tau)(0) = \varphi^\odot(\tau)(0) = 0$ for all $\tau \in \mathbb{R}$ and $m \in \mathbb{N}$. Then, we can write $\varphi_m^\odot = (0,g_m)$ and $\varphi^\odot = (0,g)$, where $g_m(\tau)(\theta) = (\theta-h)f_m(\tau-\theta)$ and $g(\tau)(\theta) = (\theta-h)f(\tau-\theta)$ for all $\theta \in (0,h]$, and we set $g_m(\tau)(0) = g(\tau)(0) = 0$ for all $\tau \in \mathbb{R}$ and $m \in \mathbb{N}$. Let us now verify that $(\varphi_m^\odot)_m$ is a sequence in $\mathcal{D}(\mathcal{A}^{\odot})$. Clearly, $\varphi_m^\odot(\tau)(0) = \varphi_m^\odot(\tau)(0^+) = 0$, the map $\tau \mapsto g_m(\tau)$ is continuously differentiable and $T$-periodic, and $\theta \mapsto g_m(\tau)(\theta)$ is in $L^1([0,h],\mathbb{C}^{n\star})$, which proves that in $(\varphi_m^\odot)_m$ is a sequence in $C_T^1(\mathbb{R},X^{\odot})$ as $g_m(\tau)(h) = 0$. Moreover, $g_m(\tau) \in \NBV([0,h],\mathbb{C}^{n \star})$ which shows that $\varphi_m^\odot(\tau) \in \mathcal{D}(A^{\star}(\tau)) = \mathcal{D}(A_0^\star)$, where this equality of the domains was obtained in \cref{lemma:adjoint D(Astar)}. It remains to show that $g_m(\tau) \in X^{\odot}$. Notice that $g_m(\tau)(0) = 0$ and $g_m(\tau) = (-hf_m(\tau), \theta \mapsto f_m(\tau-\theta) +(\theta-h)f_m'(\tau-\theta))$ where the map $\theta \mapsto f_m(\tau-\theta) +(\theta-h)f_m'(\tau-\theta)$ is clearly in $L^{1}([0,h],\mathbb{C}^{n \star})$ since $f_m$ is $C^1$-smooth. Hence, $(\varphi_m^{\odot})_m$ is a sequence in $\mathcal{D}(\mathcal{A}^{\odot})$. We claim that $\varphi_m^\odot \to \varphi^\odot$ in norm as $m \to \infty$. To see this, notice that
\begin{equation*}
    \| \varphi_m^\odot - \varphi^\odot \| = \sup_{\tau \in \mathbb{R}}\| g_m(\tau) - g(\tau) \|_{L^1} \leq 2h^2 \| f_m - f \|_{\infty} \to 0, \quad m \to \infty,
\end{equation*}
due to the uniform convergence of $(f_m)_m$ towards $f$. Using \eqref{eq:Asunclosable}, a direct calculation shows that $(\mathcal{A}^{\odot} \varphi_m^{\odot})(\tau)(\theta) = -hf_m(\tau) + F_m(\tau,\theta)$. To show that $\mathcal{A}^{\odot} \varphi_m^{\odot}$ converges in norm towards $\psi^\odot$, where $\psi^\odot(\tau)(\theta) = -hf(\tau) + F(\tau,\theta)$, notice that
\begin{equation*}
    \| \mathcal{A}^{\odot} \varphi_m^{\odot} - \psi^\odot \| \leq h \| f_m - f \|_\infty + \sup_{\tau \in \mathbb{R}} \| F_m(\tau,\cdot)-F(\tau,\cdot)\|_{\AC} \leq 2h \| f_m-f\|_{\infty} \to 0, \quad m \to \infty,
\end{equation*}
which proves the claim. Since $\tau \mapsto \dot{\varphi}^\odot(\tau)(h) = -hf(\tau) + F(\tau,h)$, is the sum of a non-continuously differentiable function $-hf$ and a continuously differentiable function $\tau \mapsto F(\tau,h)$, we have that $\tau \mapsto \varphi^\odot(\tau)(h)$ is not $C^1$-smooth and so $\varphi^\odot$ is not in $\mathcal{D}(\mathcal{A}^{\odot})$, i.e. $\mathcal{A}^\odot$ is not closed.

Since this counterexample is constructed via a sequence of $C^1$-smooth functions, it can also be used to prove analogously that $\mathcal{A}^{\star}$ is not closed.
\end{proof}
Note that we can decompose $\mathcal{A}^\odot$ as the difference of linear operators 
\begin{equation} \label{eq:Adecomp2}
    \mathcal{D}(\mathcal{A}^\odot) = \mathcal{D}(M_{A^\odot}) \cap \mathcal{D}(D^\odot), \quad \mathcal{A}^\odot = M_{A^\odot} + D^\odot,
\end{equation}
where $M_{A^\odot} : \mathcal{D}(M_{A^\odot}) \to C_T(\mathbb{R},X^\odot)$ is the (complexified) multiplication operator (associated to $A^\odot$) and $D^\odot : \mathcal{D}(D^\odot) \to C_T(\mathbb{R},X^\odot)$ is the (complexified) differential operator. Both operators can be defined as in \eqref{eq:Arecht} and \eqref{eq:Drecht}, but on the sun level. To put \eqref{eq:Adecomp2} in the light of evolution semigroups, let us introduce for any $s \geq 0$ the operator $\mathcal{U}^\odot(s) : C_T(\mathbb{R},X^\odot) \to C_T(\mathbb{R},X^\odot)$ by
\begin{equation} \label{eq:evolutionsemisun}
    (\mathcal{U}^\odot(s)f)(\tau) := U^\odot(\tau-s,\tau)f(\tau-s), \quad \forall f \in C_T(\mathbb{R},X),
\end{equation}
and note that $\mathcal{U}^\odot := \{\mathcal{U}^\odot(s)\}_{s \geq 0}$ is a family of bounded linear operators on $C_T(\mathbb{R},X^\odot)$ forming a semigroup. This family is called the \emph{sun evolution semigroup} (associated to $U^\odot$) and a similar proof as performed in \cref{lemma:C0semigroup} shows that $\mathcal{U}^\odot $ is a well-defined $\mathcal{C}_0$-semigroup. Hence, it has a generator $\hat{\mathcal{A}}^\odot$ with closed and densely defined domain $\mathcal{D}(\hat{\mathcal{A}}^\odot) \subseteq C_T(\mathbb{R},X^\odot)$. Analogous to \eqref{eq:inclusionsA}, we have the inclusions
\begin{equation} \label{eq:inclusionsAodot}
    \mathcal{D}(\mathcal{A}^\odot) \subseteq \mathcal{D}(\overline{\mathcal{A}^\odot}) \subseteq \mathcal{D}(\hat{\mathcal{A}}^\odot) \subseteq C_T(\mathbb{R},X^\odot),
\end{equation}
where $\overline{\mathcal{A}^\odot}$ denotes the closure of $\mathcal{A}^\odot$, which exists due to \cref{prop:closable2}. Again, we want to show that $\mathcal{D}(\overline{\mathcal{A}^\odot}) = \mathcal{D}(\hat{\mathcal{A}}^\odot)$ and that $\mathcal{D}(\mathcal{A}^\odot)$ is dense in $C_T(\mathbb{R},X^\odot)$ as this will be important in \cref{subsec: duality relations}. To proceed, we use a similar construction as performed in \cref{subsec: time periodic smooth Jordan} and therefore, we will construct a set $\mathcal{D}_{\mathcal{A}^\odot} \subseteq \mathcal{D}(\mathcal{A}^\odot)$ that is dense in $C_T(\mathbb{R},X^\odot)$. To construct this set, consider for any $s\in \mathbb{R}$ the function $\alpha \in C_T^1(\mathbb{R},\mathbb{R})$ satisfying $\alpha(s-2h) = 0$, the initial condition $\psi \in \mathcal{D}(A^\odot(s))$ and introduce the map $f_{\alpha,s,\psi}^\odot : \mathbb{R} \to X^\odot$ by
\begin{equation} \label{eq:falphasun}
    f_{\alpha,s,\psi}^\odot(\tau) := \alpha(\tau) U^\odot(\tau,s+lT)\psi, \quad \tau \in (s+(l-1)T-2h,s+lT-2h], \quad l \in \mathbb{Z}, 
\end{equation}
and note that this function is well-defined as $A^\odot$ is $T$-periodic, see \eqref{eq:D(Asuns)}.

\begin{proposition} \label{prop:DAsundense}
The set $\mathcal{D}_{\mathcal{A}^\odot} := \spn \{f_{\alpha,s,\psi}^\odot : s \in \mathbb{R}, \varphi \in \mathcal{D}(A^\odot(s)), \alpha \in C_T^1(\mathbb{R},\mathbb{R}) \mbox{ with } \alpha(s-2h) = 0 \} \subseteq \mathcal{D}(\mathcal{
A}^\odot)$ is dense in $C_T(\mathbb{R},X^\odot)$ and the generator $\hat{\mathcal{A}}^\odot$ of the sun evolution semigroup $\mathcal{U}^\odot$ is the closure of the linear operator $\mathcal{A}^\odot$.
\end{proposition}
\begin{proof}
The proof is essentially identical to the proof of \cref{prop:DAdense} but in the sun setting. One uses \cref{lemma:Dsuninvariant} instead of \cref{lemma:Dinvariant} and therefore the $+h$ in the action of $f_{\alpha,s,\varphi}$ gets replaced by a $-2h$ in the action of $f_{\alpha,s,\psi}^\odot$, see \eqref{eq:falpha} and \eqref{eq:falphasun} respectively. Moreover, as $U^\odot$ is strongly continuous, $\mathcal{D}(A^\odot(\tau))$ is dense in $X^\odot$ for all $\tau \in \mathbb{R}$ and $\|U^\odot(s,\tau) \| = \|U(\tau,s)^\odot \| \leq \|U(\tau,s)^\star \| = \|U(\tau,s)\|$ for all $\tau \geq s$, all claims in the proof of \cref{prop:DAdense} can be lifted to the sun level.
\end{proof}

Similar to \eqref{eq:Adecomp2}, we can decompose the non-closed but closable unbounded linear operator $\mathcal{A}^{\star}$ as 
\begin{equation*}
    \mathcal{D}(\mathcal{A}^\star) = \mathcal{D}(M_{A^\star}) \cap \mathcal{D}(D^\star), \quad \mathcal{A}^\star = M_{A^\star} + D^\star,
\end{equation*}
where $M_{A^{\star}} : \mathcal{D}(M_{A^\star}) \to C_T(\mathbb{R},X^\star)$ is the (complexified) multiplication operator (associated to $A^{\star}$) and $D^{\star} : \mathcal{D}(D^\star) \to C_T(\mathbb{R},X^\star)$ denotes the (complexified) differential operator. Both operators can be defined as in \eqref{eq:Arecht} and \eqref{eq:Drecht}, but on the star level. One can introduce similarly as in \eqref{eq:evolutionsemisun} the semigroup $\mathcal{U}^{\star} := \{\mathcal{U}^{\star}(s)\}_{s \geq 0}$ of bounded linear operators on $C_T(\mathbb{R},X^{\star})$, where
\begin{equation*}
    (\mathcal{U}^{\star}(s)f)(\tau) := U^{\star}(\tau,\tau-s)f(\tau-s), \quad \forall f \in C_T(\mathbb{R},X^{\star}), \ s \geq 0.
\end{equation*}
However, as the backward evolutionary system $U^{\star}$ is not strongly continuous (\cref{sec:periodic center manifolds}), we can not expect that $\mathcal{U}^{\star}$ forms a $\mathcal{C}_0$-semigroup on $C_T(\mathbb{R},X^{\star})$.

\begin{remark} \label{remark:curlyAadjoint}
    Consider for a moment the finite-dimensional ODE setting from \cref{remark:ODEclosed}. Here, the (complexified) unbounded linear operator $\mathcal{A}^\star : \mathcal{D}(\mathcal{A}^\star) \to C_T(\mathbb{R},\mathbb{C}^{n \star})$ takes the form
\begin{equation*}
    \mathcal{D}(\mathcal{A}^\star) = C_T^1(\mathbb{R},\mathbb{C}^{n \star}) \subseteq  C_T(\mathbb{R},\mathbb{C}^{n \star}), \quad \mathcal{A}^\star\varphi^\star = (\tau \mapsto \varphi^\star(\tau) A^\star(\tau) - \dot{\varphi}^\star(\tau)),
\end{equation*}
where $A^\star(\tau) = A(\tau)$ is the matrix $Df(\gamma(\tau)) \in \mathbb{R}^{n \times n}$. Note that $\varphi^\star(\tau)$ is now a row vector, and therefore $\varphi^\star(\tau)$ is placed in front of $A^\star(\tau)$. In the articles on the periodic normalization method \cite{Kuznetsov2005,Witte2013,Witte2014}, the authors worked for the adjoint problem on the space $C_T(\mathbb{R},\mathbb{C}^n)$. In that case, $A^\star(\tau)$ denotes the transpose of the matrix $Df(\gamma(\tau))$ and $\mathcal{A}^\star \varphi^\star = (\tau \mapsto A^\star(\tau) \varphi^\star(\tau) - \dot{\varphi}^\star(\tau))$, where $\varphi^\star(\tau)$ is a column vector. Moreover, $\mathcal{A}^\star$ is closed due to a similar proof as performed in \cref{remark:ODEclosed}. Recall from \cref{remark:ODEclosed} that the sun-star calculus construction in finite-dimensional ODEs becomes trivial and thus $\mathcal{A}^{\odot} = \mathcal{A}^\star$ is a closed linear operator. \hfill $\lozenge$
\end{remark}

Similar as for the (generalized) eigenfunctions associated to a Floquet multiplier $\lambda \in \mathbb{R}_{-}$, we need a $T$-antiperiodic dual version of \Cref{prop:eigenfunctions2}. The proof will be omitted as it is analogous to the proof of \Cref{prop:eigenfunctions2} but one relies on \Cref{thm:adjoint eigenfunctions}, the dual version of \cref{thm:eigenfunctions}.

\begin{proposition} \label{prop:adjoint2}
Let $\lambda \in \mathbb{R}_{-}$ be a Floquet multiplier of algebraic multiplicity $m_\lambda$ with $\sigma$ its associated Floquet exponent. Then there exist $T$-antiperiodic maps $\varphi_{m_\lambda - 1}^\odot,\dots,\varphi_{0}^\odot \in C^{k+1}(\mathbb{R},X^\odot)$ satisfying \eqref{eq:ODE adjoint}.
\end{proposition}
We conclude that, similar as at the end of \cref{subsec: time periodic smooth Jordan}, that the whole construction from this subsection also holds when the Floquet multiplier $\lambda \in \mathbb{R}_{-}$, i.e. the adjoint (generalized) eigenfunctions are $T$-antiperiodic.

\subsection{Duality and spectral relations} \label{subsec: duality relations}
Let us now take the time to study the duality and spectral relations between the non-closed but closable unbounded linear operators $\mathcal{A}, \mathcal{A}^{\star}, \mathcal{A}^{\odot}$ and $\mathcal{A}^{\odot \star}$ from \cref{subsec: time periodic smooth Jordan} and \cref{subsec: dual time periodic smooth Jordan}. Recall from the definition of the (complexified) natural dual pairing $\langle \cdot , \cdot \rangle : X^\star \times X \to \mathbb{C}$ that
\begin{alignat*}{2}
    \langle U^\star(s-T,s) \varphi^\star, \varphi \rangle &= \langle \varphi^\star, U(s+T,s)\varphi \rangle, \quad &&\forall (\varphi^\star,\varphi) \in X^\star \times X, \\
    \langle A^\star(s) \varphi^\star, \varphi \rangle &= \langle \varphi^\star, A(s) \varphi \rangle, \quad &&\forall (\varphi^\star,\varphi) \in \mathcal{D}(A^\star(s)) \times \mathcal{D}(A(s)).
\end{alignat*}
By construction, $U^\star(s-T,s)$ and $A^\star(s)$ are the unique linear operators satisfying these relations since $\mathcal{D}(A(s))$ is norm dense in $X$, recall \cref{sec:periodic center manifolds}. Since $X^{\odot}$ is again a Banach space, one can play this game once more, and obtain
\begin{alignat*}{2}
    \langle U^{\odot \star}(s+T,s) \varphi^{\odot \star}, \varphi^\star \rangle &= \langle \varphi^{\odot\star}, U^\odot(s-T,s)\varphi^\odot \rangle, \quad &&\forall (\varphi^{\odot \star},\varphi^\odot) \in X^{\odot \star} \times X^\odot, \\
    \langle A^{\odot\star}(s) \varphi^{\odot \star}, \varphi^\odot \rangle &= \langle \varphi^{\odot \star}, A^\odot(s) \varphi^\odot \rangle, \quad &&\forall (\varphi^{\odot \star},\varphi^\odot) \in \mathcal{D}(A^{ \odot \star}(s)) \times \mathcal{D}(A^\odot(s)),
\end{alignat*}
where $\langle \cdot, \cdot \rangle : X^{\odot \star} \times X^\odot \to \mathbb{C}$ is now the (complexified) natural dual pairing between the (complexified) Banach spaces $X^\odot$ and $X^{\odot \star}$. Again, $U^{\odot \star}(s+T,s)$ and $A^{\odot \star}(s)$ are the unique linear operators satisfying these relations since $\mathcal{D}(A^{\odot}(s))$ is norm dense in $X^{\odot}$, recall \cref{sec:periodic center manifolds}.

Let us now turn our attention to the study of the (point) spectrum of the four mentioned operators of interest. As sufficiently large iterates of $U(s+T,s)$ are compact, see \cite[Corollary XII.3.4]{Diekmann1995} and \eqref{eq:Uperiodic}, $\sigma(U(s+T,s))$ is a countable set in $\mathbb{C}$, independent of $s \in \mathbb{R}$, consisting of $0$ and isolated eigenvalues of finite type (Floquet multipliers) that can possibly accumulate to $0$, recall \cref{sec:periodic center manifolds}. Since $U(s+T,s) \in \mathcal{L}(X)$, it follows from Schrauder's theorem that sufficiently large iterates of $U^\star(s-T,s) \in \mathcal{L}(X^\star)$ are also compact and thus the spectrum of $U^\star(s-T,s)$ is a countable set in $\mathbb{C}$, independent of $s \in \mathbb{R}$, consisting of $0$ and isolated eigenvalues of finite type that can possibly accumulate to $0$. Recalling \eqref{eq:spectraleq}, we conclude
\begin{align}\label{eq:spectraequalU}
\begin{split}
    \sigma_p(U(s+T,s)) &= \sigma_p(U^{\star}(s-T,s)) \\
    &= \sigma_p(U^\odot(s-T,s)) \\
    &= \sigma_p(U^{\odot \star}(s+T,s)) = \{ \lambda \in \mathbb{C} : \mbox{$\lambda$ is a Floquet multiplier} \},  
\end{split}
\end{align}
where $\sigma_p(M)$ denotes the point spectrum of a bounded linear operator $M$ defined on a complex Banach space $E$. The aim of this subsection is to show that a similar duality relation, spectral equality \eqref{eq:spectraequalU} and decomposition \eqref{eq:decomposition} hold for the operators $\mathcal{A}, \mathcal{A}^\star, \mathcal{A^\odot}$ and $\mathcal{A}^{\odot \star}$. The results in this subsection will be helpful for our upcoming article on the study of bifurcations of limit cycles in classical DDEs, but are interesting on its own.

To construct a duality relation for the unbounded linear operators $\mathcal{A}, \mathcal{A}^\star, \mathcal{A}^{\odot}$ and $\mathcal{A}^{\odot \star}$, it is clear that we first need bilinear maps defined on the (complexified) spaces $C_T(\mathbb{R},X^\star) \times C_T(\mathbb{R},X)$ and $C_T(\mathbb{R},X^{\odot \star}) \times C_T(\mathbb{R},X^{\odot})$. To do this, let us introduce the (complexified) bilinear maps $\langle \cdot, \cdot \rangle_T : C_T(\mathbb{R},X^\star) \times C_T(\mathbb{R},X) \to \mathbb{C}$ and $\langle \cdot, \cdot \rangle_T : C_T(\mathbb{R},X^{ \odot \star}) \times C_T(\mathbb{R},X^\odot) \to \mathbb{C}$, also called \emph{pairings}, by
\begin{align} \label{eq:pairingT}
\begin{split}
    \langle \varphi^\star, \varphi \rangle_T &:= \int_0^T \langle \varphi^\star(\tau), \varphi(\tau) \rangle d\tau, \quad \forall (\varphi^\star,\varphi) \in  C_T(\mathbb{R},X^\star) \times  C_T(\mathbb{R},X), \\
    \langle \varphi^{\odot\star}, \varphi^\odot \rangle_T &:= \int_0^T \langle \varphi^{\odot \star}(\tau), \varphi^{\odot}(\tau) \rangle d\tau, \quad \forall (\varphi^{\odot\star},\varphi^\odot) \in  C_T(\mathbb{R},X^{\odot \star}) \times  C_T(\mathbb{R},X^\odot),
\end{split}
\end{align}
where the natural dual pairing $\langle \cdot, \cdot \rangle$ is between $X^\star$ and $X$ in the first integral, and between $X^{\odot \star}$ and $X^\odot$ in the second integral, recall \cref{sec:periodic center manifolds}. An important preliminary result, that we need in several upcoming results, is the following.

\begin{lemma} \label{lemma:pairingnondeg}
The bilinear maps $\langle \cdot, \cdot \rangle_T$ from \eqref{eq:pairingT} are nondegenerate.
\end{lemma}
\begin{proof}
The main idea of the proof is to extend the nondegeneracy of the natural duality pairings, appearing in the integrands of \eqref{eq:pairingT}, towards the pairings $\langle \cdot, \cdot \rangle_T$. Recall from the definition of nondegeneracy of bilinear maps that one must prove that the linear maps $ \varphi^\star \mapsto \langle \varphi^\star, \cdot \rangle_T$ and $ \varphi \mapsto \langle \cdot, \varphi \rangle_T$ are injective. We first prove the claim for the first pairing appearing in \eqref{eq:pairingT}. 

Let us start by proving that the map $ \varphi^\star \mapsto \langle \varphi^\star, \cdot\rangle_T$ is injective. Suppose that $\varphi^\star \neq 0$, then there exists a $\tau_0 \in (0,T)$ such that $\varphi^\star(\tau_0) \neq 0$. Hence, there exists by linearity a $x \in X$ such that $\operatorname{Re} \langle \varphi^\star(\tau_0), x \rangle_T = 2 > 0$. Because the map $\tau \mapsto \operatorname{Re} \langle \varphi^\star(\tau), x \rangle$ is continuous, there exists a $\delta \in (0,\min \{\tau_0,T-\tau_0\})$ such that $\operatorname{Re} \langle \varphi^\star(\tau), x \rangle \geq 1$ for all $\tau \in [\tau_0-\delta,\tau_0+\delta]$. Let $\varepsilon \in (0,\delta)$ be given and consider the $X$-valued map $\varphi_\varepsilon$ defined on $[0,T]$ by
\begin{equation*}
\varphi_\varepsilon(\tau) :=
\begin{dcases}
    \bigg(1-\frac{|\tau - \tau_0|}{\varepsilon}\bigg)x , \quad &\tau \in (\tau_0 - \varepsilon, \tau_0 + \varepsilon), \\
    0, \quad &\tau \in [0,\tau_0-\varepsilon] \cup [\tau_0+\varepsilon,T].
\end{dcases}    
\end{equation*}
Extending this function $T$-periodically on $\mathbb{R}$ shows that $\varphi_\varepsilon \in C_T(\mathbb{R},X)$. Hence, a straightforward calculation proves that $\operatorname{Re} \langle \varphi^\star, \varphi_\varepsilon \rangle_T \geq \varepsilon > 0$ and thus $\langle \varphi^\star, \varphi_\varepsilon \rangle_T \neq 0$. We conclude that the map $ \varphi^\star \mapsto \langle \varphi^\star, \cdot\rangle_T$ is injective.

To prove that the map $ \varphi \mapsto \langle \cdot, \varphi \rangle_T$ is injective, we assume that $\varphi(\tau_0) \neq 0$ for some $\tau_0 \in (0,T)$. By a corollary of the (geometric) Hahn-Banach theorem, there exists a $x^\star \in X^\star$ such that $\operatorname{Re} \langle x^\star, \varphi(\tau_0) \rangle_T > 0$. The same strategy as above can be used now to prove that $\langle \varphi_\varepsilon^\star, \varphi \rangle_T \neq 0$ for some appropriate $\varphi_\varepsilon^\star \in C_T(\mathbb{R},X^\star)$ with sufficiently small $\varepsilon > 0$.

The proof on the nondegeneracy of the second pairing in \eqref{eq:pairingT} is analogous because $X$ and $X^\odot$ are both Banach spaces.
\end{proof}

\begin{proposition} \label{prop:pairingcurlyA}
For any $\varphi \in \mathcal{D}(\mathcal{A}), \varphi^\star \in \mathcal{D}(\mathcal{A}^\star), \varphi^\odot \in \mathcal{D}(\mathcal{A}^\odot)$ and $\varphi^{\odot \star} \in \mathcal{D}(\mathcal{A}^{\odot \star})$, there holds
\begin{equation} \label{eq:curlyadjoint}
    \langle \mathcal{A}^\star \varphi^\star, \varphi \rangle_T = \langle \varphi^\star,\mathcal{A} \varphi \rangle_T, \quad \langle \mathcal{A}^{\odot \star} \varphi^{\odot \star}, \varphi^\odot \rangle_T = \langle \varphi^{\odot \star},\mathcal{A}^\odot \varphi^\odot \rangle_T.
\end{equation}
Moreover $\mathcal{A}^\star$ and $\mathcal{A}^{\odot \star}$ are the unique linear operators satisfying \eqref{eq:curlyadjoint}.
\end{proposition}
\begin{proof}
Let $\varphi \in \mathcal{D}(\mathcal{A})$ and $\varphi^\star \in \mathcal{D}(\mathcal{A}^\star)$ be given. Then,
\begin{align*}
    \langle \varphi^\star,\mathcal{A} \varphi \rangle_T &= \int_0^T \langle \varphi^\star(\tau), A(\tau)\varphi(\tau) - \dot{\varphi}(\tau) \rangle d\tau
    = - \int_0^T \langle \varphi^\star(\tau),\dot{\varphi}(\tau) \rangle d\tau + \int_0^T \langle \varphi^\star(\tau),A(\tau)\varphi(\tau) \rangle d\tau.
\end{align*}
Performing integration by parts on the first integral shows that $\langle \varphi^\star, \dot{\varphi} \rangle_T = -\langle \dot{\varphi}^\star, \varphi \rangle_T$ as $\varphi^\star$ and $\varphi$ are $T$-periodic. Hence,
\begin{equation*}
    \langle \varphi^\star,\mathcal{A} \varphi \rangle_T = \int_0^T \langle \dot{\varphi}^\star(\tau), \varphi(\tau) \rangle d\tau + \int_0^T \langle A^\star(\tau) \varphi^\star(\tau), \varphi(\tau) \rangle d\tau = \langle \mathcal{A}^\star \varphi^\star, \varphi \rangle_T,
\end{equation*}
which proves the claim. Let us now prove that $\mathcal{A}^\star$ is the unique operator satisfying the first equation in \eqref{eq:curlyadjoint}. Let $\mathcal{E}^\star : \mathcal{D}(\mathcal{A}^\star) \to C_T(\mathbb{R},X^\star)$ be a linear operator satisfying $\langle \mathcal{E}^\star \varphi^\star, \varphi \rangle_T = \langle \varphi^\star,\mathcal{A} \varphi \rangle_T$, or equivalently $\langle (\mathcal{E}^\star - \mathcal{A}^\star)\varphi^\star,\varphi \rangle = 0$, for all $(\varphi^\star,\varphi) \in \mathcal{D}(\mathcal{A}^\star)\times \mathcal{D}(\mathcal{A})$. We claim that $\langle (\mathcal{E}^\star - \mathcal{A}^\star)\varphi^\star,\varphi \rangle = 0$ for all $\varphi \in C_T(\mathbb{R},X)$. Recall from \cref{prop:DAdense} that $\mathcal{D}(\mathcal{A})$ is norm dense in $C_T(\mathbb{R},X)$ and so there exists for every $\varphi \in C_T(\mathbb{R},X)$ a sequence $(\varphi_m)_m$ in $\mathcal{D}(\mathcal{A})$ such that $\varphi_m \to \varphi$ in norm as $m \to \infty$. Since $(\mathcal{E}^\star - \mathcal{A}^\star)\varphi^\star$ takes values in $C_T(\mathbb{R},X^\star)$, there exists a $\psi^\star \in C_T(\mathbb{R},X^\star)$ such that $(\mathcal{E}^\star - \mathcal{A}^\star)\varphi^\star = \psi^\star$. As $\langle (\mathcal{E}^\star - \mathcal{A}^\star)\varphi^\star,\varphi_m \rangle = 0$ for all $m \in \mathbb{N}$, we obtain
\begin{equation*}
    |\langle (\mathcal{E}^\star - \mathcal{A}^\star)\varphi^\star,\varphi \rangle_T| = |\langle \psi^\star,\varphi_m - \varphi \rangle_T| \leq T \| \psi^\star \|_\infty \| \varphi - \varphi_m \|_\infty, \quad m \to \infty,
\end{equation*}
which proves that $\langle (\mathcal{E}^\star - \mathcal{A}^\star)\varphi^\star,\varphi \rangle = 0$ for all $\varphi \in C_T(\mathbb{R},X)$. As $\langle \cdot , \cdot \rangle_T$ is nondegenerate by \cref{lemma:pairingnondeg}, we obtain $(\mathcal{E}^\star - \mathcal{A}^\star)\varphi^\star = 0$ for all $\varphi^\star \in \mathcal{D}(\mathcal{A}^\star)$ and so $\mathcal{E}^\star = \mathcal{A}^\star$ on $\mathcal{D}(\mathcal{A}^\star)$, i.e. $\mathcal{A}^\star$ is uniquely determined. The proof for the second identity in \eqref{eq:curlyadjoint} and uniqueness of $\mathcal{A}^{\odot \star}$ is analogous, but one uses \cref{prop:DAsundense} to guarantee that $\mathcal{D}(\mathcal{A}^\odot)$ is norm dense in $C_T(\mathbb{R},X^\odot)$.
\end{proof}
The previous result shows that there exists a duality relation between, on the one hand $\mathcal{A}$ and $\mathcal{A}^\star$, and on the other hand $\mathcal{A^\odot}$ and $\mathcal{A}^{\odot \star}$. However, the bilinear maps $\langle \cdot, \cdot \rangle_T$ are not the natural dual pairings since $C_T(\mathbb{R},X)^\star \neq C_T(\mathbb{R},X^\star)$ and $C_T(\mathbb{R},X^\odot)^\star \neq C_T(\mathbb{R},X^{\odot \star})$. Hence, we can not expect in general that the standard results from adjoint theory for closable linear operators hold, see \cref{remark:mazur} and in particular \cite[Appendix B]{Engel2000} for more information.

Let us now study the spectra of the non-closed but closable unbounded linear operators $\mathcal{A}, \mathcal{A}^\star, \mathcal{A}^{\odot}$ and $\mathcal{A}^{\odot \star}$. To do this, we first recall general spectral theory for closable, not necessarily closed, unbounded linear operators, see for example \cite{Engel2000,Hille1957,Goldberg1966,Gohberg1990,Taylor1986}. Let $E$ be a complex Banach space and consider a closable (unbounded) linear operator $A: \mathcal{D}(A) \to E$ with \emph{domain} $\mathcal{D}(A)$, a linear subspace of $E$. A complex number $z$ belongs to the \emph{resolvent set} $\rho(A)$ of $A$ if the operator $zI - A$ is injective, has dense range $\mathcal{R}(zI-A)$, and the \emph{resolvent} of $A$ at $z$ defined by $(z I - A)^{-1}$ is a densely defined bounded linear operator from $\mathcal{R}(zI-A)$ to $\mathcal{D}(A)$. The \emph{spectrum} $\sigma(A)$ of $A$ is defined to be the complement of $\rho(A)$ in $\mathbb{C}$, and the \emph{point spectrum} $\sigma_p(A)$ of $A$ is the set of those $\mu \in \mathbb{C}$ such that $\mu I - A$ is not injective, i.e. $A \varphi = \mu \varphi$ for some nonzero \emph{eigenvector} $\varphi \in \mathcal{D}(A)$ corresponding to the \emph{eigenvalue} $\mu$. 

As we are only interested in \eqref{eq:spectraequalU} in the point spectra, we will only look into the point spectra of the operators $\mathcal{A}, \mathcal{A}^\star, \mathcal{A}^\odot$ and $\mathcal{A}^{\odot \star}$. A full detailed analysis on the relation between the spectra of $\mathcal{A}, \mathcal{A}^\star, \mathcal{A}^{\odot}$ and $\mathcal{A}^{\odot \star}$, and their point spectra will be provided in our upcoming paper on the study of bifurcations of limit cycles in classical DDEs.
\begin{proposition}
The point spectra of $\mathcal{A}, \mathcal{A}^\star, \mathcal{A}^\odot$ and $\mathcal{A}^{\odot \star}$ satisfy
\begin{equation} \label{eq:pointspectra}
    \sigma_p(\mathcal{A}) = \sigma_p(\mathcal{A}^\star) = \sigma_p(\mathcal{A}^\odot) = \sigma_p(\mathcal{A}^{\odot \star}) =\{ \sigma \in \mathbb{C} : \sigma \mbox{ is a Floquet exponent} \}.
\end{equation}
\end{proposition}
\begin{proof}
It follows from \cref{thm:eigenfunctions} and \cref{thm:adjoint eigenfunctions} that $\{ \sigma \in \mathbb{C} : \sigma \mbox{ is a Floquet exponent} \}$ is a subset of $\sigma_p(\mathcal{A}), \sigma_p(\mathcal{A}^\star), \sigma_p(\mathcal{A}^\odot)$ and $\sigma_p(\mathcal{A}^{\odot \star})$. To prove the converse, assume that $z \in \sigma_p(\mathcal{A})$, then $\mathcal{A}\varphi = z \varphi$ for some nonzero $\varphi \in \mathcal{D}(\mathcal{A})$. This equation is of the form \eqref{eq:zeta1} and so $\varphi(\tau) = e^{-z(\tau-s)}U(\tau,s)\phi(s)$ 
for some nonzero initial condition $\phi(s)$, defined at starting time $s \in \mathbb{R}$. As $\varphi$ must be $T$-periodic, there holds $U(s+T,s) \phi(s) = e^{z T} \phi(s)$, i.e. $\lambda = e^{zT}$ is a Floquet multiplier meaning that $z$ is a Floquet exponent. Hence, there is an equality between the first and last set of \eqref{eq:pointspectra}. 

Let us now prove an equality between the second and fourth set of \eqref{eq:pointspectra}. Let $z \in \sigma_p(\mathcal{A}^\star)$, then $\mathcal{A}^\star\varphi^\odot = z \varphi^\odot$ for some nonzero $\varphi^\odot \in \mathcal{D}(\mathcal{A}^\star)$. This is an equation of the form \eqref{eq:adjoint_eigUsun} and so $\varphi^\odot(\tau) = e^{-z(s-\tau)}U^\odot(\tau,s)\phi^\odot(s)$ for some nonzero initial condition $\phi^\odot(s)$, defined at starting time $s \in \mathbb{R}$. As $\varphi^\odot$ must be $T$-periodic, there holds $U^\odot(s-T,s)\phi^\odot(s) = e^{zT} \phi^\odot(s)$, meaning that $e^{zT} \in \sigma_p(U^\odot(s-T,s))$. Due to the discussion below \eqref{eq:spectraequalU}, $\lambda = e^{zT}$ is a Floquet multiplier and so $z$ is a Floquet exponent.

Analogous proofs to compute the point spectrum of $\mathcal{A}^\odot$ and $\mathcal{A}^{\odot \star}$ using \cref{thm:eigenfunctions} and \cref{thm:adjoint eigenfunctions} in combination with the discussion below \eqref{eq:spectraequalU} proves the claim. The $T$-antiperiodic case is proven similarly.
\end{proof}

Our next aim is to construct a direct sum decomposition as in \eqref{eq:decomposition}, but now for the spaces $C_T(\mathbb{R},X)$ and $C_T(\mathbb{R},X^{\odot \star})$. To do this, let $\lambda \in \mathbb{C} \setminus \mathbb{R}_{-}$ be Floquet multiplier of algebraic multiplicity $m_\lambda$ with $\sigma$ its associated Floquet exponent and let $k_\lambda$ denote the order of a pole of $\lambda = \mu$ of the map $\mu \mapsto (\mu I - U(\tau+T,\tau))^{-1}$ for all $\tau \in \mathbb{R}$. Recall from \cref{thm:eigenfunctions} that we constructed a  $T$-periodic $C^{k+1}$-smooth basis of the (generalized) eigenspace $E_\lambda(\tau) = \spn \{ \varphi_0(\tau),\dots,\varphi_{m_\lambda - 1}(\tau) \}$. Therefore, it is clear that $E_\sigma := \spn \{\varphi_0,\dots,\varphi_{m_\sigma - 1} \}$ is a $m_\sigma$-dimensional subspace of $C_T(\mathbb{R},X)$, where $m_\sigma := m_\lambda$ is called \emph{algebraic multiplicity} (of $\sigma$). Let us also introduce $1 \leq k_\sigma < \infty$ as the order of the pole $\sigma = \mu$ of the map $\mu \mapsto (\mu I - \mathcal{A})^{-1}$. With this notation, it follows directly from \cref{thm:eigenfunctions} that $E_\sigma = \mathcal{N}((\sigma I - \mathcal{A})^{k_\sigma})$ and let us also introduce $R_\sigma := \mathcal{R}((\sigma I - \mathcal{A})^{k_\sigma})$. The following result is a $C_T(\mathbb{R},X)$-version of \eqref{eq:directsumdecomp}.

\begin{lemma} \label{lemma:CTdecomp}
Let $\lambda \in \mathbb{C} \setminus \mathbb{R}_{-}$ be a Floquet multiplier with associated Floquet exponent of algebraic multiplicity $m_\sigma$ and let $k_\sigma$ be the order of of the pole $\sigma = \mu$ of the map $\mu \mapsto (\mu I - \mathcal{A})^{-1}$. Then, we have the direct sum decomposition into two closed $\mathcal{A}$-invariant subspaces
\begin{equation*}
    C_T(\mathbb{R},X) = E_\sigma \oplus R_\sigma.
\end{equation*}
\end{lemma}
\begin{proof}
Recall from general spectral theory that $k_\sigma$ is the smallest integer such that $\cup_{j \in \mathbb{N}} \mathcal{N}((\sigma I - \mathcal{A})^{j}) = \mathcal{N}((\sigma I - \mathcal{A})^{k_\sigma})$. According to \cite[Theorem 8.41.F]{Taylor1986}, we have that $\mathcal{D}(\mathcal{A}) = E_\sigma \oplus R_\sigma$, where both $\mathcal{A}$-invariant subspaces in the direct sum are norm closed. But recall from \cref{prop:DAdense} that $\mathcal{A}$ is densely defined and thus taking the norm closure of this direct sum equality proves the claim.    
\end{proof}
This allows us to introduce in addition the closed $\mathcal{A}^{\star}$-invariant subspaces $E_\sigma^\star := \mathcal{N}((\sigma I - \mathcal{A}^\star)^{k_\sigma})$ and $R_\sigma^\star := \mathcal{R}((\sigma I - \mathcal{A}^\star)^{k_\sigma})$, and its sun and sun-star counterparts denoted by $E_\sigma^\odot, R_\sigma^\odot, E_\sigma^{\odot \star}$ and $R_{\sigma}^{\odot \star}$. Moreover, it follows from \cref{thm:eigenfunctions} and \cref{thm:adjoint eigenfunctions} that 
\begin{equation} \label{eq:perpequalities}
    \iota(E_\sigma) = E_\sigma^{\odot \star}, \quad E_\sigma^\star = E_\sigma^\odot,
\end{equation}
which is an extension of \eqref{eq:adjointeigenspaces}. This allows us to prove the following result, which is a generalization of \cref{lemma:decomposition}.

\begin{proposition} \label{prop:decomposition}
Let $\lambda \in \mathbb{C} \setminus \mathbb{R}_{-}$ be a Floquet multiplier with associated Floquet exponent $\sigma$ of algebraic multiplicity $m_\sigma$ and let $k_\sigma$ be the order of of the pole $\sigma = \mu$ of the map $\mu \mapsto (\mu I - \mathcal{A})^{-1}$. Then,
\begin{equation} \label{eq:perp}
    R_\sigma^\perp = E_\sigma^\star, \quad R_\sigma = E_\sigma^{\star \perp}, \quad R_\sigma^{\odot \perp} = E_\sigma^{\odot \star}, \quad R_\sigma^\odot = E_\sigma^{\odot \star \perp},
\end{equation}
where the annihilator $\perp$ is defined in terms of the pairings $\langle \cdot,\cdot \rangle_T$ from \eqref{eq:pairingT}. In particular, we have the direct sum decompositions
\begin{equation} \label{eq:decompCT2}
    C_T(\mathbb{R},X) = E_\sigma \oplus E_\sigma^{\odot \perp}, \quad C_T(\mathbb{R},X^\odot) = E_\sigma^\odot \oplus \iota(E_\sigma)^{\perp}.
\end{equation}
\end{proposition}
\begin{proof}
Let us prove the first equality of \eqref{eq:perp}. If $\varphi^\star \in E_\sigma^\star$, then $\langle \varphi^\star, (\sigma I - \mathcal{A})^{k_\sigma} \varphi \rangle_T = \langle (\sigma I - \mathcal{A}^\star)^{k_\sigma}\varphi^\star,\varphi \rangle_T = 0$ for all $\varphi \in \mathcal{D}(\mathcal{A})$ and thus $\varphi^\star \in R_\sigma^\perp$. Here, we used the first equation of \eqref{eq:curlyadjoint} $m_\lambda$ times. The proof for the other inclusion is analogous. Taking the annihilator of the first equality in \eqref{eq:perp} yields $R_\sigma \subseteq R_\sigma^{\perp \perp} = E_\sigma^{\star \perp} \subseteq C_T(\mathbb{R},X)$, and so the first inclusion of the second equality in \eqref{eq:perp} is proven. Here we used \cite[Theorem 8.9.2]{Narici2010} in the first inclusion and noticed from the same reference that $R_\sigma^{\perp \perp}$ is a linear subspace of $C_T(\mathbb{R},X)$.

To prove the other inclusion, assume by contraction that there exists a $\varphi_0 \in E_\sigma^{\star \perp}$ such that $\varphi_0 \notin R_\sigma$. For a fixed $\tau_0 \in (0,T)$, we obtain $\varphi_0(\tau_0) \notin R_\lambda(\tau_0)$, where we recall from \cref{subsec: dual time periodic smooth Jordan} that $R_\lambda(\tau_0)$ is a norm closed subspace of $X$. A corollary of the (geometric) Hahn Banach theorem guarantees the existence of a $x^\star \in X^\star$ such that $\operatorname{Re} \langle x^\star, \varphi_0(\tau_0) \rangle = 2 > 0$ and $\langle x^\star, \psi \rangle = 0$ for all $\psi \in R_\lambda(\tau_0)$. The aim is to show that this construction extends over the whole period, i.e. we will show that there exists a $\varphi_\varepsilon^\star \in C_T(\mathbb{R},X^\star)$ for some sufficiently small $\varepsilon > 0$ such that $\operatorname{Re} \langle \varphi_\varepsilon^\star,\varphi \rangle_T > 0$ and $ \langle \varphi_\varepsilon^\star, (\sigma I - \mathcal{A})^{k_\sigma} \varphi \rangle_T = 0$ for all $\varphi \in \mathcal{D}(\mathcal{A})$. Since $x^\star \in R_\lambda(\tau_0)^\perp = E_\lambda^\odot(\tau_0)$, recall \cref{lemma:decomposition}, we know from \cref{thm:adjoint eigenfunctions} that $x^\star = \sum_{j=0}^{m_\lambda - 1} \mu_j \varphi_j^\odot(\tau_0)$ for some $\mu_{m_\lambda - 1},\dots,\mu_{0} \in \mathbb{C}$. Hence, we introduce the $T$-periodic map $\hat{x}^\star : \mathbb{R} \to X^\star$ by
\begin{equation*}
    \hat{x}^\star(\tau) := \sum_{j=0}^{m_\lambda - 1} \mu_j \varphi_j^\odot(\tau), \quad \forall \tau \in \mathbb{R},
\end{equation*}
and note that $\hat{x}^\star(\tau_0) = x^\star$. Recall from \cref{thm:adjoint eigenfunctions} that the adjoint (generalized) eigenfunctions $\varphi_{m_\lambda-1}^\odot,\dots,\varphi_0^\odot$ are at least continuous and so $\hat{x}^\star \in C_T(\mathbb{R},X^\star)$ by linearity. Because the map $\tau \mapsto \langle \hat{x}^\star(\tau), \varphi_0(\tau) \rangle$ is continuous and $\operatorname{Re} \langle \hat{x}^\star(\tau_0), \varphi_0(\tau_0) \rangle = 2$, there exists a $\delta \in (0,\min\{\tau_0,T-\tau_0\})$ such $\operatorname{Re} \langle \hat{x}^\star(\tau), \varphi_0(\tau) \rangle \geq 1$ for all $\tau \in [\tau_0-\delta,\tau_0+\delta]$. Introduce now for any $\varepsilon \in (0,\delta)$ the $X^\star$-valued map $\varphi_\varepsilon^\star$ on $[0,T]$ by 
\begin{equation*}
\varphi_\varepsilon^\star(\tau) :=
\begin{dcases}
    \bigg(1-\frac{|\tau - \tau_0|}{\varepsilon}\bigg)\hat{x}^\star(\tau) , \quad &\tau \in (\tau_0 - \varepsilon, \tau_0 + \varepsilon), \\
    0, \quad &\tau \in [0,\tau_0-\varepsilon] \cup [\tau_0+\varepsilon,T].
\end{dcases}    
\end{equation*}
Extending this function $T$-periodically on $\mathbb{R}$ proves that $\varphi_\varepsilon^\star \in C_T(\mathbb{R},X^\star)$ and a straightforward calculation shows that $\operatorname{Re} \langle \varphi_\varepsilon^\star, \varphi_0 \rangle_T \geq \varepsilon > 0$. Let us now prove that $ \langle \varphi_\varepsilon^\star, (\sigma I - \mathcal{A})^{k_\sigma} \varphi \rangle_T = 0$ for all $\varphi \in \mathcal{D}(\mathcal{A})$. Let $\varphi \in \mathcal{D}(\mathcal{A})$ be given and set $\psi = (\sigma I - \mathcal{A})^{k_\sigma} \varphi$. Then,
\begin{align*}
    \langle \varphi_\varepsilon^\star, (\sigma I - \mathcal{A})^{k_\sigma} \varphi \rangle_T = \sum_{j=0}^{m_\lambda - 1} \mu_j \int_{\tau_0-\varepsilon}^{\tau_0+\varepsilon} \bigg(1-\frac{|\tau - \tau_0|}{\varepsilon} \bigg) \langle \varphi_j^\odot(\tau), \psi(\tau) \rangle d\tau.
\end{align*}
Since $\varphi_j^\odot(\tau) \in E_\lambda^\odot(\tau)$ and $\psi(\tau) \in R_\lambda(\tau) = E_\lambda^\star(\tau)^\perp = E_\lambda^\odot(\tau)^\perp$, recall \cref{thm:adjoint eigenfunctions} and \eqref{eq:decomposition}, we have that the natural dual pairings $\langle \varphi_j^\odot(\tau), \psi(\tau) \rangle$ vanish for all $\tau \in \mathbb{R}$ and $j=m_\lambda-1,\dots,0$. Thus, $\langle \varphi_\varepsilon^\star, (\sigma I - \mathcal{A})^{k_\sigma} \varphi \rangle_T = 0$ for all $\varphi \in \mathcal{D}(\mathcal{A})$, which proves the claim. To finalize the proof, notice that $\langle (\sigma I - \mathcal{A}^\star)^{k_\sigma}\varphi_\varepsilon^\star, \varphi \rangle_T = 0$ for all $\varphi \in \mathcal{D}(\mathcal{A})$, but then $(\sigma I - \mathcal{A}^\star)^{k_\sigma}\varphi_\varepsilon^\star = 0$, or equivalently $\varphi_\varepsilon^\star \in E_\sigma^\star$, by norm density of $\mathcal{D}(\mathcal{A})$ (\cref{prop:DAdense}) and nondegeneracy of the pairing $\langle \cdot, \cdot \rangle_T$ (\cref{lemma:pairingnondeg}). But then in particular $\langle \varphi_\varepsilon^\star, \varphi_0 \rangle_T = 0$ because $\varphi_0 \in E_\sigma^{\star \perp}$, which is a contradiction. Hence, we conclude that $E_\sigma^{\star \perp} \subseteq R_\sigma$, which proves the second equality of \eqref{eq:perp}. The remaining two equalities from \eqref{eq:perp} are proven analogously. 

The first equality in \eqref{eq:decompCT2} follows from \eqref{eq:perp} in combination with \eqref{eq:perpequalities}. To prove the second equality in \eqref{eq:decompCT2}, recall from \cref{prop:DAsundense} that $\mathcal{A}^\odot$ is also densely defined and thus a similar proof as for \cref{lemma:CTdecomp} can be used to conclude, in combination with \eqref{eq:perpequalities}, that the second direct sum decomposition in \eqref{eq:decompCT2} holds.
\end{proof}

\begin{remark} \label{remark:mazur}
Recall from previous results that we essentially lifted, for a specific example, the Hahn-Banach theorem on $X \times X^\star$ towards $C_T(\mathbb{R},X) \times C_T(\mathbb{R},X^\star)$, where $\langle \cdot , \cdot \rangle_T$ is not a (complexified) natural dual pairing but just a nondegenerate bilinear map. This puts, according to \cite[Chapter 8]{Narici2010}, the spaces $C_T(\mathbb{R},X)$ and $C_T(\mathbb{R},X^\star)$ in \emph{duality}. Let us study for a moment this topic from a more topological point of view. According to \cite[Theorem 8.9.2]{Narici2010}, we have the characterization 
\begin{equation} \label{eq:weakclosure}
    R_\sigma^{\perp \perp} = \overline{R_\sigma}^{\weak},
\end{equation}
where the closure is with respect to the weak topology on $C_T(\mathbb{R},X)$ induced by $C_T(\mathbb{R},X^\star)$ in terms of the first pairing $\langle \cdot , \cdot \rangle_T$ from \eqref{eq:pairingT}. This topology will be denoted by $\sigma(C_T(\mathbb{R},X),C_T(\mathbb{R},X^\star))$, see \cite{Narici2010,Engel2000} for other examples. This topology is the weakest TVS topology on $C_T(\mathbb{R},X)$ making all the (complexified) linear maps $\langle \varphi^\star , \cdot \rangle_T : C_T(\mathbb{R},X) \to \mathbb{C}$ continuous as $\varphi^\star$ ranges over $C_T(\mathbb{R},X^\star)$. Recall from \cref{lemma:pairingnondeg} that $C_T(\mathbb{R},X)$ separates the points of $C_T(\mathbb{R},X^\star)$ and therefore $\sigma(C_T(\mathbb{R},X),C_T(\mathbb{R},X^\star))$ is Hausdorff, recall \cite[Theorem 8.2.2]{Narici2010}. If the annihilator $\perp$ in \cref{prop:decomposition} would be in terms of the natural dual pairing, this is a well-known result, see for example \cite[Theorem 8.41.F]{Taylor1986} or \cite[Theorem IV.2.5]{Diekmann1995}. However, all equivalent proofs of this result in the literature  rely somehow on the Hahn-Banach theorem. For example, \eqref{eq:weakclosure} holds where the weak closure would be defined in terms of the standard weak topology induced via the dual space in terms of the natural dual pairing. As every linear subspace is convex, Mazur's theorem \cite[Theorem V.1.4]{Conway1985}, a corollary of the (geometric) Hahn Banach theorem, guarantees that \eqref{eq:weakclosure} also equals its norm closure, but recall from \cref{prop:decomposition} that $R_\sigma$ is norm closed. It is unclear to the authors if the $\sigma(C_T(\mathbb{R},X),C_T(\mathbb{R},X^\star))$ topology is a so-called \emph{topology of the pair} \cite[Definition 8.7.1]{Narici2010}, but we doubt it. If this was true, Mazur's theorem would hold immediately \cite[Theorem 8.8.1]{Narici2010} and the proof of \cref{prop:decomposition} could be simplified. \hfill $\lozenge$
\end{remark}

The following result will be helpful in our upcoming article, when we will study codimension one bifurcations of limit cycles in classical DDEs. However, this result can also be used to study generic bifurcations of limit cycles of higher codimension.

\begin{corollary} \label{cor:decompsimple}
Let $\lambda \in \mathbb{C} \setminus \mathbb{R}_{-}$ be a Floquet multiplier with associated Floquet exponent $\sigma$ of algebraic multiplicity $m_\sigma$ and geometric multiplicity $1$. Then, we have the direct sum decompositions
\begin{align*}
    C_T(\mathbb{R},X) &= \mathcal{N}((\sigma I -\mathcal{A})^{m_\sigma}) \oplus \mathcal{N}((\sigma I -\mathcal{A}^\odot)^{m_\sigma})^{\perp}, \\ 
    C_T(\mathbb{R},X^\odot) &= \mathcal{N}((\sigma I -\mathcal{A}^\odot)^{m_\sigma}) \oplus \mathcal{N}((\sigma I -\mathcal{A}^{\odot \star})^{m_\sigma})^{\perp}.
\end{align*}
\end{corollary}
\begin{proof}
If the geometric multiplicity of the Floquet exponent is $1$, then $k_\sigma = m_\sigma$, recall for example \cite[Section IV.1.17]{Engel2000}, and so the result follows from \cref{prop:decomposition}.
\end{proof}

When the Floquet multiplier $\lambda \in \mathbb{R}_{-}$, we know from \cref{prop:eigenfunctions2} and \cref{prop:adjoint2} that we have to work with $T$-antiperiodic (adjoint) (generalized) eigenfunctions. All the results from this subsection also hold when we replace $T$-periodicity with $T$-antiperiodicity. Only some parts of several proofs have to be modified a little bit to ensure $T$-antiperiodicity of certain functions.

\begin{remark}
Consider the finite-dimensional ODE setting from \cref{remark:ODEclosed} and \cref{remark:curlyAadjoint}. Clearly, $\mathcal{A} = \mathcal{A}^{\odot \star}$ and $\mathcal{A^\star} = \mathcal{A}^\odot$ are densely defined linear operators as $C_T^1(\mathbb{R},\mathbb{C}^{n})$ and $C_T^1(\mathbb{R},\mathbb{C}^{n \star})$ are norm dense in $C_T(\mathbb{R},\mathbb{C}^{n})$ and $C_T(\mathbb{R},\mathbb{C}^{n \star})$, respectively. We claim that the pairing $\langle \cdot , \cdot \rangle_T : C_T(\mathbb{R},\mathbb{C}^{n \star}) \times C_T(\mathbb{R},\mathbb{C}^{n}) \to \mathbb{C}$ given by
\begin{equation*}
    \langle \varphi^\star,\varphi \rangle_T = \int_0^T \varphi^\star(\tau) \varphi(\tau) d\tau,
\end{equation*}
is nondegenerate. If $\langle \cdot , \varphi \rangle_T = 0$, choose $\varphi_0^\star$ as the Hermitian adjoint of $\varphi$ since then $\langle \varphi_0^\star, \varphi \rangle_T = \int_0^T | \varphi(\tau)|^2 d\tau = 0$ and so $\varphi = 0$ by continuity of $\varphi$. If $\langle \varphi^\star , \cdot \rangle_T = 0$, choose $\varphi_0$ as the Hermitian adjoint of $\varphi^\star$ to conclude that $\varphi^\star =0$. Moreover, for all $\varphi \in \mathcal{D}(\mathcal{A})$ and $\varphi \in \mathcal{D}(\mathcal{A}^\star)$, we get
\begin{align} \label{eq:dualityODE}
\begin{split}
    \langle \varphi^\star, \mathcal{A}\varphi \rangle_T &= \int_0^T \varphi^\star(\tau) A(\tau) \varphi(\tau) - \varphi^\star(\tau) \dot{\varphi}(\tau) d\tau \\
    &= \int_0^T \varphi^\star(\tau) A^\star(\tau) \varphi(\tau) + \dot{\varphi}^\star(\tau) {\varphi}(\tau) d\tau = \langle \mathcal{A}^\star \varphi^\star, \varphi \rangle_T,
\end{split}
\end{align}
where we used integration by parts and the $T$-periodicity of $\varphi$ and $\varphi^\star$. As $\mathcal{A}$ is densely defined, $\mathcal{A}^\star$ is the unique operator satisfying \eqref{eq:dualityODE}. The equality \eqref{eq:pointspectra} follows from a similar proof, where \cref{thm:eigenfunctions} must be replaced by its ODE-variant \cite[Proposition III.1]{Iooss1999} or \cite[Lemma 2]{Iooss1988}. \cref{prop:decomposition} holds as well and the proof simplifies as a Hahn-Banach argument is not necessary. \cref{cor:decompsimple} in the light of finite-dimensional ODEs justifies the fact that particular pairings between adjoint (generalized) eigenfunctions and (generalized) eigenfunctions from the periodic normalization method \cite{Kuznetsov2005,Witte2013,Witte2014} do not vanish, indeed. \hfill $\lozenge$
\end{remark}

\section{Characterization of the center manifold and normal forms} \label{sec:characterization}
In this section, we study the dynamics of \eqref{eq:DDEphi} on the center manifold $\mathcal{W}_{\loc}^c(\Gamma)$ near the \emph{nonhyperbolic} cycle $\Gamma$, meaning that there are, except of the trivial Floquet multiplier, other Floquet multipliers present on the unit circle in the complex plane, or the trivial Floquet multiplier has an algebraic multiplicity higher than one. Recall from \Cref{sec:introduction} that there are three generic codimension one bifurcation of limit cycles: the fold bifurcation, when the trivial Floquet multiplier has algebraic multiplicity $2$ and geometric multiplicity $1$, the period-doubling bifurcation, when a simple Floquet multiplier is located at $-1$, and the Neimark-Sacker bifurcation when a simple complex conjugate pair of Floquet multipliers is located on the unit circle. 

To study these bifurcations, we first separate the dynamics generated by the trivial Floquet multiplier from the rest of the dynamics in \Cref{subsec: separating}. Afterwards, in \Cref{subsec: normal form}, we prove the existence of a special coordinate system on the center manifold and provide in addition the periodic (critical) normal forms of the mentioned three codimension one bifurcations of limit cycles. These results are an extension of the work by Iooss (and Adelmeyer) \cite{Iooss1988,Iooss1999} from finite-dimensional ODEs to infinite-dimensional DDEs. As a consequence, the periodic normal forms obtained in \cref{thm:normalformI}, \cref{thm:normalformII} and \cref{thm:normalformIII} for classical DDEs are exactly the same as obtained by Iooss for ODEs, see \cite[Theorem III.7, Theorem III.10 and Theorem III.13]{Iooss1999}. Instead of only being interested in codim 1 bifurcations of limit cycles, the provided framework is also suited to study bifurcations of limit cycles of higher codimension, see for example \cite[Table 1]{Witte2013} and \cite{Witte2014} for the periodic normal forms of all generic codim 2 bifurcations of limit cycles in ODEs, and hence DDEs. It is nevertheless helpful to keep, throughout the upcoming construction, these three codimension one bifurcations in mind.

Before we start proving the characterization and periodic normal form theorems, let us first recall some useful results from \cite{Lentjes2023} on the topological direct sum decomposition of $X$ and $X^{\odot \star}$. It turns out from \cite[Proposition 10]{Lentjes2023} that we can decompose $X$ as
\begin{equation} \label{eq:decomposition X hyp}
        X = X_{-}(\tau) \oplus X_0(\tau) \oplus X_{+}(\tau), \quad \forall \tau \in \mathbb{R},
\end{equation}
where $X_{-}(\tau)$ and $X_{+}(\tau)$ denote the \emph{stable eigenspace} and \emph{unstable eigenspace} (at time $\tau$) respectively, see \cite[Section 3.6]{Lentjes2023} for their definitions. Furthermore, it turns out from \cite[Appendix A.1]{Lentjes2023} that we can lift the decomposition \eqref{eq:decomposition X hyp} towards a decomposition in $X^{\odot \star}$ as
\begin{equation} \label{eq:decomposition Xsunstar hyp}
        X^{\odot \star} = X^{\odot \star}_{-}(\tau) \oplus X^{\odot \star}_{0}(\tau) \oplus X^{\odot \star}_{+}(\tau), \quad \forall \tau \in \mathbb{R},
\end{equation}
where $X^{\odot \star}_{0}(\tau) = j(X_0(\tau))$ and $X^{\odot \star}_{+}(\tau) = j(X_+(\tau))$, see \cite[Appendix A.2]{Lentjes2023} for more information.

\subsection{Separating the dynamics of the periodic orbit} \label{subsec: separating} \label{subsec:separating}
The coordinate system and normal forms we will present consist of two parts and is inspired by \cite{Iooss1988,Iooss1999}. The first part expresses the dynamics along $\Gamma$ by a time-dependent phase, and the other part expresses the dynamics transverse to $\Gamma$ in terms of this phase. The normal forms depend on the location and multiplicities of the Floquet multipliers on the unit circle. Let us first separate the dynamics of the periodic orbit in terms of coordinates along this phase and transverse to this phase. 

Recall that $X_0(\tau)$ is a $(n_0 + 1)$-dimensional subspace of $X$ for all $\tau \in \mathbb{R}$. For each $\lambda \in \Lambda_0$, we know that the (generalized) eigenspace $E_\lambda(\tau)$ has a ($T$ or 2$T$)-periodic basis that satisfies the conditions from \Cref{thm:eigenfunctions} or \Cref{prop:eigenfunctions2}, depending on the location of $\lambda$ on the unit circle. Recall that the trivial Floquet multiplier is always present on the unit circle and $\dot{\gamma}_\tau$ is the associated $T$-periodic eigenfunction of $U(\tau+T,\tau)$. We choose $\varphi_0(\tau)$ to be $\dot{\gamma}_\tau$ and denote $\tilde{X}_{0}(\tau)$ as the space spanned by $\varphi_1(\tau),\dots,\varphi_{n_0}(\tau)$ that forms a ($T$ or $2T$)-periodic $C^{k+1}$-smooth basis as presented in \Cref{thm:eigenfunctions} or \Cref{prop:eigenfunctions2}. Define for any $\tau \in \mathbb{R}$ the bounded linear operator $\tilde{Q}_0(\tau): \mathbb{R}^{n_0} \to \tilde{X}_0(\tau)$ as
\begin{equation} \label{eq:Q0tilde}
    \tilde{Q}_0(\tau)\xi := \sum_{i=1}^{n_0} \xi_i\varphi_i(\tau), \quad \forall \xi = (\xi_1,\dots,\xi_{n_0}) \in \mathbb{R}^{n_0}.
\end{equation}
With this notation, it is clear that the Floquet operator (at time $\tau$) associated to $\Lambda_0$, denoted by $Q_0(\tau) : \mathbb{R} \times \mathbb{R}^{n_0} \to X_0(\tau)$ has action
\begin{equation*}
    Q_0(\tau)(\xi_0,\xi) := \xi_0 \dot{\gamma}_\tau + \tilde{Q}_0(\tau)\xi, \quad \forall (\xi_0,\xi) \in \mathbb{R} \times \mathbb{R}^{n_0}.
\end{equation*}
The $({n_0}+1) \times ({n_0}+1)$ matrix $M_0$ from takes the form 
\begin{equation} \label{eq:M0}
    M_0 =  \left(
    \begin{array}{c|c}
      0 & \star \ \cdots \ 0\\
      \hline
      \vspace{-4pt}
      0\\
      \vdots & \tilde{M}_0\\
      0
    \end{array}
    \right),
\end{equation}
where $\star \in \{0,1\}$ depends on the algebraic multiplicity of the trivial Floquet multiplier.

\subsection{Characterization and normal form theorems} \label{subsec: normal form}
Depending on the algebraic multiplicity of the trivial Floquet multiplier and the location of the other Floquet multipliers on the unit circle, the periodic normal forms will have a different shape and therefore three different normal form theorems will be presented.

\begin{figure}[ht]
    \centering
    \includegraphics[width = 11cm]{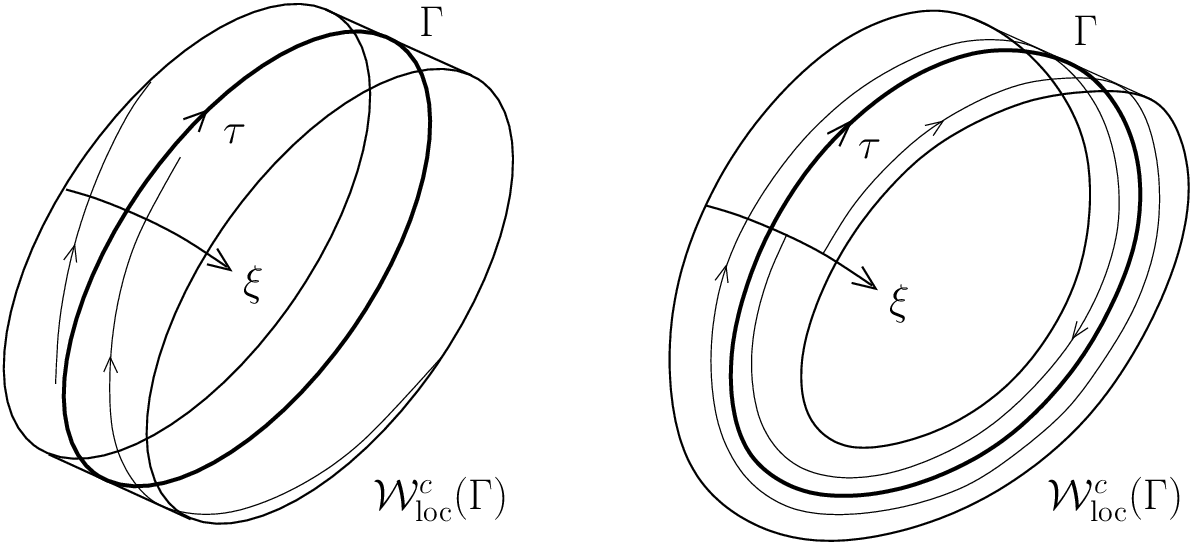}
    \caption{Illustration of two-dimensional center manifolds $\mathcal{W}_{\loc}^{c}(\Gamma)$ together with the coordinate system $(\tau,\xi)$. The left figure represents the case when $-1 \not \in \Lambda_0$ and then $\mathcal{W}_{\loc}^{c}(\Gamma)$ is locally diffeomorphic to a cylinder in a neighborhood of $\Gamma$, see \Cref{thm:normalformII}. The right figure represents the case when $-1 \in \Lambda_0$ and then $\mathcal{W}_{\loc}^{c}(\Gamma)$ is locally diffeomorphic to a Möbius band in a neighborhood of $\Gamma$, see \Cref{thm:normalformIII}.}
    \label{fig:CM}
\end{figure}

The main idea to prove the existence of suitable coordinates on  $\mathcal{W}_{\loc}^{c}(\Gamma)$ is to use the invariance property of $\mathcal{W}_{\loc}^{c}(\Gamma)$ to the fullest. We try to parametrize a solution $x_t$ of \eqref{eq:DDEphi} on $\mathcal{W}_{\loc}^{c}(\Gamma)$ in the vicinity of the periodic orbit $\Gamma$ as
\begin{equation} \label{eq:x_t manifold}
    x_{t}= \gamma_\tau + \tilde{Q}_0(\tau)\xi + H(\tau,\xi).
\end{equation}
Here, $\tau$ is a function of $t$ that expresses the dynamics along $\Gamma$ by a time-dependent phase while $\xi$ is a function of $\tau$ that expresses the dynamics transverse to $\Gamma$ in terms of this phase. Such a coordinate system is visualized for a two-dimensional local center manifold around $\Gamma$ in \Cref{fig:CM}.

The only unknown in \eqref{eq:x_t manifold} is the nonlinear operator $H : \mathbb{R} \times \mathbb{R}^{n_0} \to X$. We will show in the upcoming normal form theorems that this operator exists and is sufficiently smooth. This will be done by using the invariance property of $\mathcal{W}_{\loc}^{c}(\Gamma)$. To be more precise, if we take an initial condition $x_s = \varphi \in \mathcal{W}_{\loc}^{c}(\Gamma)$, then the solution $x_t$ of \eqref{eq:DDEphi} must lie on $\mathcal{W}_{\loc}^{c}(\Gamma)$ for all $t$ in the time domain of definition, say $I \subseteq \mathbb{R}$ with $s \in I$. Hence, \cite[Theorem 3.6]{Clement1989} tells us that $x_t$ satisfies the abstract ODE
\begin{equation} \label{eq:invariance equation}
    \frac{d}{dt}j(x_t) = A_0^{\odot \star}j(x_t) + G(x_t), \quad t \in I,
\end{equation}
where the map $G \in C^{k+1}(X,X^{\odot \star})$ is defined by $G(\psi) := F(\psi)r^{\odot \star}$ for all $\psi \in X$. Here, recall that $F \in C^{k+1}(X,\mathbb{R}^n)$ denotes the right-hand side of \eqref{eq:DDEphi} for a given $k \geq 1$.

First, we consider the case where the trivial Floquet multiplier has algebraic multiplicity $1$ and there is no Floquet multiplier located at $-1$. Such a normal form theorem would be helpful to study the Neimark-Sacker bifurcation. 
\begin{theorem}[Normal Form I] \label{thm:normalformI}
Assume that the algebraic multiplicity of the trivial Floquet multiplier is one and that  $-1$ is not a Floquet multiplier. Then there exist $C^k$-smooth maps $H : \mathbb{R} \times \mathbb{R}^{n_0} \to X, p : \mathbb{R} \times \mathbb{R}^{n_0} \to \mathbb{R}$ and $P:\mathbb{R} \times \mathbb{R}^{n_0} \to \mathbb{R}^{n_0}$ such that any solution $x_t$ of \eqref{eq:DDEphi} on $\mathcal{W}_{\loc}^{c}(\Gamma)$ may be represented as
\begin{equation*}
    x_{t}= \gamma_\tau + \tilde{Q}_0(\tau)\xi + H(\tau,\xi), \quad t \in I,
\end{equation*}
where the time dependence of the coordinates $(\tau,\xi)$ describing
the dynamics of \eqref{eq:DDEphi} on $\mathcal{W}_{\loc}^{c}(\Gamma)$ is defined by the normal form
\begin{equation*}
    \begin{cases}
    \dfrac{d \tau}{dt} = 1+p(\tau,\xi) + \mathcal{O}(|\xi|^{k+1}),\\
    \\[-8pt]
    \dfrac{d \xi}{d\tau} = \tilde{M}_0 \xi + P(\tau,\xi) + \mathcal{O}(|\xi|^{k+1}).
    \end{cases}
\end{equation*}
Here the functions $H,p$, and $P$ are $T$-periodic in $\tau$ and at least quadratic in $\xi$. The $\mathcal{O}$-terms are also $T$-periodic in $\tau$.  Moreover, $p$ and $P$ are polynomials in $\xi$ of degree less than or equal to $k$ satisfying
\begin{equation*}
    \frac{d}{d\tau}p(\tau,e^{-\tau \tilde{M}_0^{\star}} \xi) = 0, \quad \frac{d}{d\tau}\bigg( e^{\tau \tilde{M}_0^{\star}} P(\tau,e^{-\tau \tilde{M}_0^{\star}} \xi)\bigg) = 0,
\end{equation*}
for all $\tau \in \mathbb{R}$ and $\xi \in \mathbb{R}^{n_0}$.
 \end{theorem}
The proof of this theorem is quite long and technical. Essentially, it is a careful generalization of \cite[Theorem III.7]{Iooss1999} and \cite[Theorem 1]{Iooss1988}. Therefore, we first sketch the idea of the proof and break it up into several steps. In \textbf{Step 1}, we prove that $\mathcal{W}_{\loc}^c(\Gamma)$ can be parametrized in terms of $(\tau,\xi)$-coordinates. In the remaining steps, we will show that any solution of \eqref{eq:DDEphi} on $\mathcal{W}_{\loc}^c(\Gamma)$, written in the original time $t$, can be expressed in $(\tau,\xi)$-coordinates by the normal form in \cref{thm:normalformI}. This will be done by characterizing the map $H$ by its Taylor expansion up to a certain order. In \textbf{Step 2}, we write down this Taylor expansion in suitable form and start in \textbf{Step 3} with collecting terms in powers of $\xi^q$ for $q=0,\dots,k$ obtained from the invariance equation of the center manifold. We get for $q=2,\dots,k$ an equation for the coefficient $H_q$ and we must show that this can be uniquely solved. This will be accomplished by decomposing $H_q$ into the decomposition provided in \eqref{eq:decomposition X hyp} together with the separation of the dynamics along $\Gamma$ made in \Cref{subsec: separating}, see \textbf{Step 4}. Hence, we get for each $q=2,\dots,k$ the terms $H_q^{00},\tilde{H}_q^0, H_q^{-}$ and $H_q^{+}$ and then prove the existence and continuity of each of these terms in \textbf{Step 5} $(H_q^{+})$, \textbf{Step 6} $(H_q^{-})$ and \textbf{Step 7} ($H_q^{00}$ and $\tilde{H}_q^0$). This proves that the map $H$ is $C^k$-smooth. The provided normal forms are partially derived in step 6 of the proof in combination with \cite[Theorem III.7]{Iooss1999}.

\begin{proof}[Proof of \Cref{thm:normalformI}]
\textbf{Step 1: Existence.} Consider a point $\Psi_0 \in \mathcal{W}_{\loc}^c(\Gamma)$ and consider an associated fixed pair $(\tau_0,\xi_0) \in \mathbb{R} \times \mathbb{R}^{n_0}$. Results from \cref{subsec:separating} in combination with \eqref{eq:decomposition X hyp} allows us to write
\begin{equation} \label{eq:Psi0}
    \Psi_0 = \gamma_{\tau_0} + \Tilde{Q}_0(\tau_0)\xi_0 + \hat{H}(\tau_0,\xi_0),
\end{equation}
where $\gamma_{\tau_0} + \Tilde{Q}_0(\tau_0)\xi_0 \in \mathbb{R}\dot{\gamma}_{\tau_0} \oplus \tilde{X}_0(\tau_0) = X_0(\tau_0)$ and $\hat{H}(\tau_0,\xi_0) \in {X}_{-}(\tau_0) \oplus {X}_{+}(\tau_0)$. Assume for a moment that $\hat{H}$ is continuously differentiable. This will be proven later, independent on the upcoming construction. Define then the continuously differentiable operator $\mathcal{K} : \mathbb{R} \times \mathbb{R}^{n_0} \times X \to X$ by
\begin{equation*}
    \mathcal{K}((\tau,\xi),\Psi) := \Psi - (\gamma_\tau + \Tilde{Q}_0(\tau)\xi + \hat{H}(\tau,\xi)),
\end{equation*}
and notice from \eqref{eq:Psi0} that $\mathcal{K}((\tau_0,\xi_0),\Psi_0) = 0$. Since the partial derivative $D_\Psi \mathcal{K}((\tau_0,\xi_0),\Psi_0)$ is the identity operator on $X$, we deduce from the implicit function theorem in Banach spaces that there exists a continuously differentiable function $\Psi : \mathbb{R} \times \mathbb{R}^{n_0} \to X$, defined in a sufficiently small neighborhood of $\Gamma$, that satisfies $\mathcal{K}((\tau,\xi),\Psi(\tau,\xi)) = 0$ or equivalently
\begin{equation} \label{eq:Psitau}
    \Psi(\tau,\xi) = \gamma_\tau + \Tilde{Q}_0(\tau)\xi + \hat{H}(\tau,\xi).
\end{equation}
This proves that $\Psi$ is a parametrization of the center manifold $\mathcal{W}_{\loc}^c(\Gamma)$ in $(\tau,\xi)$-coordinates. Similar to the proof of \cite[Theorem 1]{Iooss1988} and \cite[Theorem III.7]{Iooss1999}, we allow $\hat{H}(\tau,\cdot)$ taking values in $X_0(\tau)$. The map $\hat{H}$ with this additional property will be denoted by $H$. Hence, we incorporate an eventual nonlinear change of coordinates on $\tau$ and $\xi$ in \eqref{eq:Psitau} in such a way that \eqref{eq:DDEphi} written on $\mathcal{W}_{\loc}^c(\Gamma)$ is already in normal form. It remains to show that the map $H$ is $C^k$-smooth and that $x_t = \Psi(\tau,\xi)$, where the relation between the two (different) time scales $t$ and $\tau$ is defined by the normal form in \cref{thm:normalformI}. This will be proven in the remaining steps.

\textbf{Step 2: Taylor expansion.} Let us write \eqref{eq:DDEphi} in the form of \eqref{eq:invariance equation} and notice that
\begin{equation} \label{eq:Gexpansion}
   G(\gamma_\tau + \varphi) = G(\gamma_\tau) + B(\tau)\varphi + \sum_{q=2}^k G_q(\tau,\varphi^{(q)}) + \mathcal{O}(\|\varphi\|_{\infty}^{k+1}), \quad \forall \varphi \in X,
\end{equation}
where $B(\tau)\varphi = [DF(\gamma_\tau)\varphi] r^{\odot \star}$ is the time-dependent bounded linear perturbation and the nonlinear terms are given by $G_q(\tau,\varphi^{(q)}) := \frac{1}{q!}D^{q}F(\gamma_\tau)(\varphi^{(q)})r^{\odot \star},$ where $D^{q}F(\gamma_\tau) : X^{q} \to \mathbb{R}^n$ is the $q$th order Fr\'echet derivative evaluated at $\gamma_\tau$ and $\varphi^{(q)} := (\varphi,\dots,\varphi) \in X^{q} := X \times \dots \times  X$ for $q=2,\dots,k$. We also expand the maps $H,p$ and $P$ as
\begin{align*}
    H(\tau,\xi) &= \sum_{q=2}^k H_q(\tau,\xi^{(q)}) + \mathcal{O}(|\xi|^{k+1}), \quad p(\tau,\xi) = \sum_{q=2}^k p_q(\tau,\xi^{(q)}), \quad  P(\tau,\xi) = \sum_{q=2}^k P_q(\tau,\xi^{(q)}),
\end{align*}
with coefficients  $H_q(\tau,\xi^{(q)}) \in X, p_q(\tau,\xi^{(q)}) \in \mathbb{R}$ and $P_q(\tau,\xi^{(q)}) \in \mathbb{R}^{n_0}$, where $\xi^{(q)} := (\xi,\dots,\xi) \in [\mathbb{R}^{n_0}]^{q}$. As already announced, we will use the invariance property of $\mathcal{W}_{\loc}^{c}(\Gamma)$ to obtain a representation of the coefficients $H_q(\tau,\xi^{(q)})$ for all $q = 2,\dots,k$. Hence, we compare the expansions of 
\begin{equation*}
    \frac{d}{dt}j(\gamma_\tau + \tilde{Q}_0(\tau)\xi + H(\tau,\xi)) = j\bigg(\dot{\gamma}_\tau + \frac{\partial \tilde{Q}_0(\tau)}{\partial \tau}\xi + \frac{\partial H(\tau,\xi)}{\partial \tau} +\bigg(\tilde{Q}_0(\tau) + D_\xi H(\tau,\xi) \bigg) \frac{d\xi}{d\tau} \bigg) \frac{d\tau}{dt}
\end{equation*}
and
\begin{equation*}
     A_0^{\odot \star}j(x_t)+G(x_t) = A_0^{\odot \star}j(\gamma_\tau + \tilde{Q}_0(\tau)\xi + H(\tau,\xi)) + G(\gamma_\tau + \tilde{Q}_0(\tau)\xi + H(\tau,\xi))
\end{equation*}
by substituting 
\begin{equation*}
    \frac{d \tau}{dt} = 1+p(\tau,\xi) + \mathcal{O}(|\xi|^{k+1}), \quad \frac{d \xi}{d\tau} = \Tilde{M}_0 \xi + P(\tau,\xi) + \mathcal{O}(|\xi|^{k+1}).
\end{equation*}
Using the expansions of $H,p$ and $P$ together with \eqref{eq:Gexpansion}, where $\varphi$ must be substituted by $\tilde{Q}_0(\tau)\xi + H(\tau,\xi)$, we get
\begin{align*}
    &j \bigg[ \dot{\gamma}_\tau + \frac{\partial \tilde{Q}_0(\tau)}{\partial \tau} \xi + \sum_{q=2}^k \frac{\partial H_q(\tau,\xi^{(q)})}{\partial \tau} + \bigg(\tilde{Q}_0(\tau) + \sum_{q=2}^k D_\xi H_q(\tau,\xi^{(q)}) \bigg) \bigg(\tilde{M}_0 + \sum_{q=2}^{k} P_q(\tau,\xi^{(q)}) \bigg) \bigg] \\
    &\bigg(1+  \sum_{q=2}^k p_q(\tau,\xi^{(q)}) \bigg) + \mathcal{O}(|\xi|^{k+1}) \\
   &= A_0^{\odot \star}j(\gamma_\tau) + G(\gamma_\tau) + A^{\odot \star}(\tau)j\bigg(\tilde{Q}_0(\tau)\xi + \sum_{q=2}^k H_q(\tau,\xi^{(q)}) \bigg) \\
   &+ \sum_{q=2}^k G_q \bigg(\tau,\bigg[\tilde{Q}_0(\tau)\xi + \sum_{p=2}^k H_q(\tau,\xi^{(p)}) \bigg]^{(q)} \bigg) + \mathcal{O}(|\xi|^{k+1}).
\end{align*}

\textbf{Step 3: Collecting terms.} Let us now compare the terms in powers of $\xi$ on both side of this equation. Collecting the $\xi^0$-terms give us
\begin{equation*}
   \frac{d}{d\tau} j({\gamma}_\tau) = A_0^{\odot \star} j(\gamma_\tau) + G(\gamma_\tau),
\end{equation*}
which means that $\tau \mapsto \gamma$ is a solution \eqref{eq:invariance equation}. This was already known since $\gamma$ is a periodic solution of \eqref{eq:DDEphi}. The $\xi^1$-terms give us
\begin{equation} \label{eq:linear part CMT}
    \bigg(-\frac{\partial}{\partial \tau} + A^{\odot \star}(\tau) \bigg) j(\tilde{Q}_0(\tau)\xi) = j(\tilde{Q}_0(\tau)\tilde{M}_0 \xi),
\end{equation}
which is exactly the result established in \eqref{eq:FloquetmapODE}, but now for all Floquet multipliers on the unit circle and this characterizes the linear part. After collecting the $\xi^{(2)}$-terms, we get the expression
\begin{align*}
    &\bigg( -\frac{\partial}{\partial \tau} + A^{\odot \star}(\tau) \bigg) j(H_2(\tau,\xi^{(2)}))\\
    &= j(D_\xi H_2(\tau,\xi^{(2)}) \tilde{M}_0 \xi + p_2(\tau,\xi^{(2)})\dot{\gamma}_\tau + \tilde{Q}_0(\tau) P_2(\tau,\xi^{(2)})) - R_2(\tau,\xi^{(2)}),
\end{align*}
where $R_2(\tau,\xi^{(2)}) = G_2(\tau,(\tilde{Q}_0(\tau) \xi)^{(2)})$. Finally, after collecting the $\xi^{(q)}$-terms for $q = 3,\dots,k$ one obtains
\begin{align} \label{eq:diffeqpsi}
        &\bigg( -\frac{\partial}{\partial \tau} + A^{\odot \star}(\tau) \bigg) j(H_q(\tau,\xi^{(q)})) \nonumber \\
        &= j(D_\xi H_q(\tau,\xi^{(q)}) \tilde{M}_0 \xi + p_q(\tau,\xi^{(q)})\dot{\gamma}_\tau + \tilde{Q}_0(\tau) P_q(\tau,\xi^{(q)})) - R_q(\tau,\xi^{(q)}),
\end{align}
where $R_q(\tau,\xi^{(q)})$ depends on $G_{q'}(\tau,\xi^{(q')})$ for $2\leq q' \leq q$ and $j(H_{q'}(\tau,\xi^{(q')}))$, $j(p_{q'}(\tau,\xi^{(q')})\dot{\gamma}_\tau)$ and $j( \tilde{Q}_0(\tau)P_{q'}(\tau,\xi^{(q')}))$ for $q' = 2,\dots,q-1$.

\textbf{Step 4: Projecting on subspaces.} We want to project \eqref{eq:diffeqpsi} onto the spaces $\mathbb{R}j\dot{\gamma}_\tau, j(\tilde{X}_0(\tau))$ and $ {X}_{\pm}^{\odot \star}(\tau)$ to obtain a representation of $H_q$ separately on each individual space. Because $X$ can be decomposed as in $\eqref{eq:decomposition X hyp}$ where $X_0(\tau) = \mathbb{R}\dot{\gamma}_\tau \oplus \tilde{X}_0(\tau)$ for any $\tau \in \mathbb{R},$ we can decompose for any $q=2,\dots,k$ the function $H_q$ as
\begin{align*}
    H_q(\tau,\xi^{(q)}) = H_q^{00}(\tau,\xi^{(q)})\dot{\gamma}_\tau + \tilde{Q}_0(\tau)\tilde{H}_q^0 (\tau,\xi^{(q)}) + H_q^{-}(\tau,\xi^{(q)}) + H_q^{+}(\tau,\xi^{(q)}),
\end{align*}
where $H_q^{\pm}(\tau,\xi^{(q)}) = P_{\pm} (\tau) H_q (\tau,\xi^{(q)}) \in X_{\pm}(\tau)$ together with $H_q^{00}(\tau,\xi^{(q)}) \in \mathbb{R}$ and $\tilde{H}_q^0 (\tau,\xi^{(q)}) \in \mathbb{R}^{n_0}$ for all $\tau \in \mathbb{R}$ and $\xi \in \mathbb{R}^{n_0}$. It follows from \eqref{eq:decomposition Xsunstar hyp} that
\begin{equation*}
    R_q(\tau,\xi^{(q)}) = R_q^{00}(\tau,\xi^{(q)})j\dot{\gamma}_\tau + j(\tilde{Q}_0(\tau)\tilde{R}_q^0 (\tau,\xi^{(q)})) + R_q^{-}(\tau,\xi^{(q)}) + R_q^{+}(\tau,\xi^{(q)}),
\end{equation*}
where $R_q^{\pm}(\tau,\xi^{(q)}) = P_{\pm}^{\odot \star} (\tau) R_q (\tau,\xi^{(q)}) \in X_{\pm}^{\odot \star}(\tau), R_q^{00}(\tau,\xi^{(q)}) \in \mathbb{R}$ and $\tilde{R}_q^0 (\tau,\xi^{(q)}) \in \mathbb{R}^{n_0}$ for all $\tau \in \mathbb{R}$ and $\xi \in \mathbb{R}^{n_0}$. The precise definition of the spectral projector $P_{\pm}^{\odot \star}$ can be found in \cite[Proposition 8]{Lentjes2023} but it is not necessary to know this form explicitly for the upcoming construction. Substituting these decompositions into \eqref{eq:diffeqpsi} yields for the left-hand side of this equation
\begin{align*}
    \bigg( -\frac{\partial}{\partial \tau} + A^{\odot \star}(\tau) \bigg) j(H_q(\tau,\xi^{(q)}))
    &= -j\bigg(\frac{\partial H_q^{00}(\tau,\xi^{(q)})}{\partial \tau} \dot{\gamma}_\tau + H_q^{00}(\tau,\xi^{(q)})\ddot{\gamma}_\tau \bigg) \\
    &+ A^{\odot \star}(\tau)j(H_q^{00}(\tau,\xi^{(q)})\dot{\gamma}_\tau) \\
    &
    -j\bigg(\frac{\partial \tilde{Q}_0(\tau)}{\partial \tau} \tilde{H}_q^0(\tau,\xi^{(q)}) + \tilde{Q}_0(\tau) \frac{\partial \tilde{H}_q^0(\tau,\xi^{(q)})}{\partial \tau}\bigg) \\
    &+ A^{\odot \star}(\tau)j(\tilde{Q}_0(\tau)\tilde{H}_q^0(\tau,\xi^{(q)}))\\
    &
    + \bigg( -\frac{\partial}{\partial \tau} + A^{\odot \star}(\tau) \bigg) j(H_q^{-}(\tau,\xi^{(q)}) + H_q^{+}(\tau,\xi^{(q)})),
\end{align*}
where we used twice the product rule for differentiation. Recall from \Cref{thm:eigenfunctions} that $(-\frac{d}{d\tau} + A^{\odot \star}(\tau))j \dot{\gamma}_\tau = 0$ and using this in combination with \eqref{eq:linear part CMT} yields
\begin{align} \label{eq:jHq}
\bigg( -\frac{\partial}{\partial \tau} + A^{\odot \star}(\tau) \bigg) j(H_q(\tau,\xi^{(q)})) &=j\bigg(- \frac{\partial H_q^{00}(\tau,\xi^{(q)})}{\partial \tau} \dot{\gamma}_\tau \bigg) \\
&+ j\bigg(\tilde{Q}_0(\tau) \bigg(- \frac{\partial}{\partial \tau} + \tilde{M}_0 \bigg) \tilde{H}_q^0(\tau,\xi^{(q)}) \bigg) \nonumber \\
&
+ \bigg( -\frac{\partial}{\partial \tau} + A^{\odot \star}(\tau) \bigg) j(H_q^{-}(\tau,\xi^{(q)}) + H_q^{+}(\tau,\xi^{(q)})) \nonumber
\end{align}
and this must be equal to the right-hand side of \eqref{eq:diffeqpsi}. Let us first obtain a representation of $H_q^{\pm}$ by projecting it onto the spaces $X_{\pm}^{\odot \star}(\tau)$. On these subspaces, we get the equation
\begin{equation*}
    \bigg( -\frac{\partial}{\partial \tau} + A^{\odot \star}(\tau) \bigg) j(H_q^{\pm}(\tau,\xi^{(q)})) = j(D_\xi H_q^{\pm}(\tau,\xi^{(q)})\tilde{M}_0 \xi) - R_q^{\pm}(\tau,\xi^{(q)}).
\end{equation*}
Substituting $\tau = \theta$ and $\xi = e^{(\theta - \tau)\tilde{M}_0} \xi =: \tilde{\xi}$ leads to
\begin{equation*}
    -\frac{\partial}{\partial \theta} j(H_q^{\pm}(\theta,\tilde{\xi}^{(q)})) + A^{\odot \star}(\theta)j(H_q^{\pm}(\theta,\tilde{\xi}^{(q)})) - j(D_{\tilde{\xi}} H_q^{\pm}(\theta,\tilde{\xi}^{(q)}) \tilde{M}_0 \xi) = - R_q^{\pm}(\theta,\tilde{\xi}^{(q)}).
\end{equation*}
When the operator $-U^{\odot \star}(\tau,\theta)$ acts on both side of the equation, we obtain
\begin{align} \label{eq:CMTcalculation}
     &-U^{\odot \star}(\tau,\theta)\bigg[-\frac{\partial}{\partial \theta} j(H_q^{\pm}(\theta,\tilde{\xi}^{(q)})) + A^{\odot \star}(\theta)j(H_q^{\pm}(\theta,\tilde{\xi}^{(q)})) - j(D_{\tilde{\xi}} H_q^{\pm}(\theta,\tilde{\xi}^{(q)}) \tilde{M}_0 \xi) \bigg] \nonumber \\
     &= U^{\odot \star}(\tau,\theta)R_q^{\pm}(\theta,\tilde{\xi}^{(q)}). 
\end{align}
Let us focus on the left-hand-side of this equation. It follows from \cite[Theorem 5.5]{Clement1988} that
\begin{equation*}
    -U^{\odot \star}(\tau,\theta)A^{\odot \star}(\theta) j(H_q^{\pm}(\theta,\tilde{\xi}^{(q)})) = - [{\partial_\theta^\star}U^{\odot \star}(\tau,\theta)] j(H_q^{\pm}(\theta,\tilde{\xi}^{(q)})).
\end{equation*}
Filling this result back into \eqref{eq:CMTcalculation} and using the partial weak$^\star$ derivative operator yields
\begin{equation*}
    U^{\odot \star}(\tau,\theta) [ {\partial_\theta^\star} j(H_q^{\pm}(\theta,\tilde{\xi}^{(q)})) ] + [{\partial_\theta^\star}U^{\odot \star}(\tau,\theta)]j(H_q^{\pm}(\theta,\tilde{\xi}^{(q)})) = U^{\odot \star}(\tau,\theta)R_q^{\pm}(\theta,\tilde{\xi}^{(q)}),
\end{equation*}
where we have used the product rule for differentiation, but essentially in the natural dual pairings due to the partial weak$^\star$ derivative. Using the product rule again and recalling that $ \tilde{\xi} = e^{(\theta - \tau)\tilde{M}_0} \xi$, we get the identity
\begin{equation*}
    \partial_\theta^\star[U^{\odot \star}(\tau,\theta) j(H_q^{\pm}(\theta,(e^{(\theta - \tau)\tilde{M}_0})^{(q)}))] = U^{\odot \star}(\tau,\theta)R_q^{\pm}(\theta,(e^{(\theta - \tau)\tilde{M}_0})^{(q)}).
\end{equation*}
Using the definition of the weak$^\star$ derivative, we get for all $x^\odot \in X^\odot$ that
\begin{equation} \label{eq:wkstarexpressionPsi}
    \frac{\partial}{ \partial \theta} \langle j(U(\tau,\theta)H_q^{\pm}(\theta,(e^{(\theta - \tau)\tilde{M}_0} \xi)^{(q)})), x^\odot \rangle = \langle U^{\odot \star}(\tau,\theta)R_q^{\pm}(\theta,(e^{(\theta - \tau)\tilde{M}_0} \xi)^{(q)}) ,x^\odot \rangle. \\
\end{equation}

\textbf{Step 5: Existence of $H_q^{+}$.} Let us first find an expression for $H_q^{+}(\tau,\xi^{(q)})$. As $X_+(s)$ is finite-dimensional, $U(\tau,s)$ extends to all $\tau,s \in \mathbb{R}$ on the subspace $X_+(s)$, recall \cite[Proposition 10]{Lentjes2023}. So let $s \geq \tau$ be given and integrate \eqref{eq:wkstarexpressionPsi} over the interval $[\tau,s]$ to obtain
\begin{align} \label{eq:wkstarintegralPsi}
    \langle j(H_q^{+}(\tau,\xi^{(q)})), x^\odot \rangle &= \langle j(U(\tau,s)H_q^{+}(s,(e^{(s - \tau)\tilde{M}_0} \xi)^{(q)})), x^\odot \rangle \nonumber \\
   & - \int_\tau^{s} \langle U^{\odot \star}(\tau,\theta)R_q^{+}(\theta,(e^{(\theta - \tau)\tilde{M}_0} \xi)^{(q)}) ,x^\odot \rangle d\theta.
\end{align}
Let us focus on the first term of the right-hand side. Notice that
\begin{equation*}
    H_q^{+}(s,(e^{(s - \tau)\tilde{M}_0} \xi)^{(q)}) = \sum_{|\alpha| = q} \frac{1}{\alpha!} P_+(s)H_{s}^\alpha ((e^{(s - \tau)\tilde{M}_0} \xi)^\alpha)
\end{equation*}
where $H_{s}^\alpha ((e^{(s - \tau)\tilde{M}_0} \xi)^\alpha) \in X$. Then, we get
\begin{equation*}
    U(\tau,s)H_q^{+}(s,(e^{(s - \tau)\tilde{M}_0} \xi)^{(q)}) = \sum_{|\alpha| = q} \frac{1}{\alpha!} U(\tau,s) P_{+}(s)H_{s}^\alpha ((e^{(s - \tau)\tilde{M}_0} \xi)^\alpha)
\end{equation*}
and using the exponential trichotomy property of the forward evolutionary system \cite[Proposition 10]{Lentjes2023}, there is a $b > 0$ such that for a given $\varepsilon > 0$ there exists a $K_\varepsilon > 0$ with the property
\begin{equation*}
   \|U(\tau,s)H_q^{+}(s,(e^{(s - \tau)\tilde{M}_0} \xi)^{(q)})\|_{\infty} \leq
   K_\varepsilon e^{b(\tau - s)} \sum_{|\alpha| = q} \frac{1}{\alpha!} \|H_{s}^\alpha ((e^{(s - \tau)\tilde{M}_0} \xi)^\alpha)\|_{\infty},
\end{equation*}
where the number $N$ from \cite[Proposition 10]{Lentjes2023} is absorbed in the $K_\varepsilon$ constant. Since the diagonal elements of the matrix $\tilde{M}_0$ have real part zero, $e^{(s - \tau)\tilde{M}_0} \xi$ is a polynomial in $\xi$ and so $\|H_{s}^\alpha ((e^{(s - \tau)\tilde{M}_0} \xi)^\alpha)\|_{\infty}$ can grow at most polynomially for $s \to \pm \infty$. With this in mind, we get
\begin{align*}
    |\langle j(U(\tau,s)H_q^{+}(s,(e^{(s - \tau)\tilde{M}_0} \xi)^{(q)})), x^\odot \rangle| & \leq K_\varepsilon e^{b(\tau - s)} \|x^\odot \| \sum_{|\alpha| = q}  \frac{1}{\alpha!} \|H_{s}^\alpha ((e^{(s - \tau)\tilde{M}_0} \xi)^\alpha)\|_{\infty} \\
    &\leq M_\varepsilon e^{b(\tau - s)} \max_{|\alpha| = q} \|H_{s}^\alpha ((e^{(s - \tau)\tilde{M}_0} \xi)^\alpha)\|_{\infty}\\
    & \to 0, \quad s \to \infty.
\end{align*}
Using this convergence, taking the limit in \eqref{eq:wkstarintegralPsi} yields
\begin{align}
    \langle j(H_q^{+}(\tau,\xi^{(q)})), x^\odot \rangle & = \langle \int_\tau^{\infty} -U^{\odot \star}(\tau,\theta) R_q^{+}(\theta,(e^{(\theta - \tau)\tilde{M}_0} \xi)^{(q)}) d\theta ,x^\odot \rangle, \label{eq:jHq+}
\end{align}
if we can show that for any fixed $x^\odot \in X^\odot$, $\tau \in \mathbb{R}$ and $\xi \in \mathbb{R}^{n_0}$ the map $g_{q,\tau,\xi}^{+} : [\tau,\infty) \to \mathbb{R} $ defined by $g_{q,\tau,\xi}^{+}(\theta) = \langle -U^{\odot \star}(\tau,\theta) R_q^{+}(\theta,(e^{(\theta - \tau)\tilde{M}_0} \xi)^{(q)}) ,x^\odot \rangle$ is in $L^1([\tau,\infty),\mathbb{R})$. From \cite[Proposition 10]{Lentjes2023} we get
\begin{align*}
    \int_\tau^\infty |g_{q,\tau,\xi}^{+}(\theta) | d\theta &\leq K_\varepsilon N \|x^\odot \| e^{b\tau} \int_\tau^\infty e^{-b \theta} \|R_q(\theta,(e^{(\theta - \tau)\tilde{M}_0} \xi)^{(q)}) \| d\theta.
\end{align*}
Recall that $e^{(\theta - \tau)\tilde{M}_0} \xi$ is a polynomial in $\xi$ and that $R_q(\tau,\xi^{(q)})$ depends on $G_{q'}(\tau,\xi^{(q')})$ for $2\leq q' \leq q$ and $j(H_{q'}(\tau,\xi^{(q')}))$, $j(p_{q'}(\tau,\xi^{(q')})\dot{\gamma}_\tau)$ and $j( \tilde{Q}_0(\tau)P_{q'}(\tau,\xi^{(q')}))$ for $q' = 2,\dots,q-1$. Since $G_{q'}$ is periodic in the first variable and evaluated at a polynomial, $G_{q'}$ grows at most polynomially for $2\leq q' \leq q$. As we can assume that $H_{q'}$ is $T$-periodic in the first variable for $q' = 2,\dots,q-1$ (we will show this later for $q' = q$) and evaluated at a polynomial, it follows that $j(H_{q'}(\tau,\xi^{(q')}))$ grows at most polynomially for $q' = 2,\dots,q-1$. As $p_{q'}$ and $P_{q'}$ are $T$-periodic in the first variable for $q' = 2,\dots,q-1$ (we will show this later for $q' = q$), it follows that $j(p_{q'}(\tau,\xi^{(q')})\dot{\gamma}_\tau)$ and $j( \tilde{Q}_0(\tau)P_{q'}(\tau,\xi^{(q')}))$ can grow at most polynomially. To conclude, there exists a polynomial $r_{q,\tau,\xi}^{+}: \mathbb{R} \to \mathbb{R}$ such that $\|R_q^+(\theta,(e^{(\theta - \tau)\tilde{M}_0} \xi)^{(q)}) \| \leq r_{q,\tau,\xi}^{+}(\theta)$ for all $\theta \geq \tau$. Hence, 
\begin{equation} \label{eq:Hq+ bounded}
    \int_\tau^\infty |g_{q,\tau,\xi}^{+}(\theta) | d\theta \leq K_\varepsilon N \|x^\odot \| e^{b\tau} \int_\tau^\infty e^{-b \theta} r_{q,\tau,\xi}^{+}(\theta) d\theta < \infty,
\end{equation}
because the map $ [\tau,\infty) \ni \theta \mapsto e^{-b \theta} r_{q,\tau,\xi}^{+}(\theta) \in \mathbb{R}$ decays fast enough to zero as $\theta \to \infty$. We have proven that the weak$^\star$ integral in \eqref{eq:jHq+} exists. Because $R_q^{+}(\theta,(e^{(\theta - \tau)\tilde{M}_0} \xi)^{(q)}) \in j(X_{+}(\theta))$ and \eqref{eq:jHq+} holds for any $x^\odot \in X^\odot$, we obtain
\begin{align*}
    j(H_q^{+}(\tau,\xi^{(q)})) &= j\int_\tau^{\infty} -U(\tau,\theta) j^{-1} R_q^{+}(\theta,(e^{(\theta - \tau)\tilde{M}_0} \xi)^{(q)}) d\theta.
\end{align*}
As $X$ is $\odot$-reflexive with respect to the shift semigroup $T_0$, we have that $j$ is an isomorphism on $j(X) = X^{\odot \odot}$ and so
\begin{equation} \label{eq:psi_q^+}
    H_q^{+}(\tau,\xi^{(q)}) = -\int_\tau^{\infty} U(\tau,\theta) j^{-1} R_q^{+}(\theta,(e^{(\theta - \tau)\tilde{M}_0} \xi)^{(q)}) d\theta
\end{equation}
can be evaluated as a standard Riemann integral. It can easily be checked that $H_q^{+}$ is $T$-periodic in the first variable because $P_{+}^{\odot \star}$ is $T$-periodic and $R_q$ is $T$-periodic in the first variable. Let us now prove the continuity of the map $H_q^{+}$. As $U^{\odot \star}(t,\tau)$ restricted to $j(X_{+}(\tau))$ is invertible, we can adjust the proof from \cite[Lemma 1]{Lentjes2023} to prove continuity of the limiting function $v(\tau,\infty,\tau,g)$ (notation from \cite[Lemma 1]{Lentjes2023})  for a continuous function $g : [\tau,\infty) \to X^{\odot \star}$ under the assumption that $H_q^{+}$ is bounded in norm. As it is proved in \eqref{eq:Hq+ bounded} that $H_q^+$ is bounded in norm and noticing that $P_+^{\odot \star}$ and $R_q$ are continuous for all $q \in \{1,\dots,k\}$, the result follows. 

\textbf{Step 6: Existence of $H_q^{-}$.} Now, we can look for an explicit expression of $H_q^{-}(\tau,\xi^{(q)})$. Integrating \eqref{eq:wkstarexpressionPsi} over $[s,\tau]$ for a fixed $s\in \mathbb{R}$ shows that 
\begin{align*}
    \langle j(H_q^{-}(\tau,\xi^{(q)})), x^\odot \rangle &= \langle j(U(\tau,s)H_q^{-}(s,(e^{(s - \tau)\tilde{M}_0} \xi)^{(q)})), x^\odot \rangle \nonumber \\
    &+ \int_\tau^{s} \langle U^{\odot \star}(\tau,\theta)R_q^{-}(\theta,(e^{(\theta - \tau)\tilde{M}_0} \xi)^{(q)}) ,x^\odot \rangle d\theta, \
\end{align*}
for all $x^\odot \in X^\odot$. Similar to the $H_q^{+}$-case, we want to show that the first term goes to zero, but now as $s \to -\infty$. Recall that $\|H_{s}^\alpha ((e^{(s - \tau)\tilde{M}_0} \xi)^\alpha)\|_{\infty}$ can grow at most polynomially for $s \to \pm \infty$ and so due to the exponential trichotomy of the forward evolutionary system \cite[Proposition 10]{Lentjes2023}, there exists an $a < 0$ such that for a given $\varepsilon > 0$ there is a $M_\varepsilon >0$ with the property
\begin{align*}
    |\langle j(U(\tau,s)H_q^{-}(s,(e^{(s - \tau)\tilde{M}_0} \xi)^{(q)})), x^\odot \rangle| & \leq  M_{\varepsilon}e^{a (\tau - s)} \max_{|\alpha| = q} \|H_s^{\alpha}(s,(e^{(s-\tau) \tilde{M}_0}\xi)^{(q)})\|_{\infty}\\
    & \to 0, \quad s \to -\infty,
\end{align*}
where the other constants are already absorbed in $M_\varepsilon$. We conclude that
\begin{equation} \label{eq:Hq-}
    \langle j(H_q^{-}(\tau,\xi^{(q)})), x^\odot \rangle = \langle \int_{- \infty}^\tau U^{\odot \star}(\tau,\theta)R_q^{-}(\theta,(e^{(\theta - \tau)\tilde{M}_0} \xi)^{(q)}) d\theta, x^\odot \rangle,  
\end{equation}
if we are able to show that for any $x^\odot \in X^{\odot}$ and fixed $\tau \in \mathbb{R}$ and $\xi \in \mathbb{R}^{n_0}$ that the map $g_{q,\tau,\xi}^{-} : (-\infty,\tau] \to \mathbb{R}$ defined by $g_{q,\tau,\xi}^{-}(\theta) = \langle U^{\odot \star}(\tau,\theta)R_q^{-}(\theta,(e^{(\theta - \tau)\tilde{M}_0} \xi)^{(q)}), x^\odot \rangle$ is in $L^1((-\infty,\tau],\mathbb{R})$. The exponential trichotomy implies that for a given $\varepsilon > 0$ one can find a $K_\varepsilon > 0$ such that 
\begin{align*}
    \int_\tau^\infty |g_{q,\tau,\xi}^{-}(\theta) | d\theta &\leq K_\varepsilon N \|x^\odot \| e^{a\tau} \int_{-\infty}^\tau e^{-a\theta} \|R_q(\theta,(e^{(\theta - \tau)\tilde{M}_0} \xi)^{(q)}) \| d\theta.
\end{align*}
From the same reasoning as in the $H_q^{+}$-case, there exists a polynomial $r_{q,\tau,\xi}^{-} : \mathbb{R} \to \mathbb{R}$ that satisfies the estimate $\|R_q(\theta,(e^{(\theta - \tau)\tilde{M}_0} \xi)^{(q)}) \| \leq r_{q,\tau,\xi}^{-}(\theta)$ for all $\theta \leq \tau$. Hence,
\begin{equation} \label{eq:Hq- bounded}
    \int_\tau^\infty |g_{q,\tau,\xi}^{-}(\theta) | d\theta \leq K_\varepsilon N \|x^\odot \| e^{a\tau} \int_{-\infty}^\tau e^{-a\theta} r_{q,\tau,\xi}^{-}(\theta) d\theta < \infty,
\end{equation}
because the map $\theta \mapsto e^{-a \theta} r_{q,\tau,\xi}^{-}(\theta)$ decays fast enough to zero as $\theta \to - \infty$. Hence, $g_{q,\tau,\xi}^{-} \in L^1((-\infty,\tau],\mathbb{R})$ and so the weak$^\star$ integral exists. Since \eqref{eq:Hq-} holds for all $x^\odot \in X^\odot$ we get
\begin{equation} \label{eq:psi_q^-}
    H_q^{-}(\tau,\xi^{(q)}) = j^{-1}\int_{- \infty}^\tau U^{\odot \star}(\tau,\theta)R_q^{-}(\theta,(e^{(\theta - \tau)\tilde{M}_0} \xi)^{(q)}) d\theta,
\end{equation}
if we can prove that the weak$^\star$ integral takes values in $j(X)$. Notice that we proved in \eqref{eq:Hq- bounded} that $H_q^{-}$ is bounded in norm. Using the notation from \cite[Lemma 1]{Lentjes2023}, we have that $j(H_q^{-}(\tau,\xi^{(q)})) = v(\tau,\tau,-\infty,g)$ where the continuous map $g$ is defined by $g(\theta) = P_{-}^{\odot \star}(\theta) R_q(\theta,(e^{(\theta - \tau)\tilde{M}_0} \xi)^{(q)})$ for all $\theta \in (-\infty,\tau]$ since $P_{-}^{\odot \star}$ and $R_q$ are continuous for all $q \in \{1,\dots,k\}$. It follows from \cite[Lemma 1]{Lentjes2023} that $H_q$ takes values in $j(X)$ and is continuous. It is not difficult to show that $H_q^-$ is $T$-periodic in the first variable because $P_{-}^{\odot \star}$ is $T$-periodic and $R_q$ is $T$-periodic in the first variable. 

\textbf{Step 7: Existence of $H_q^{00}$ and $\tilde{H}_q^{0}$.} To obtain $H_q^{00}(\tau,\xi^{(q)})$ and $H_q^{0}(\tau,\xi^{(q)})$, we project \eqref{eq:diffeqpsi} onto $\mathbb{R}j\dot{\gamma}_\tau$ and $j(\tilde{X}_0(\tau))$. Since $j$ is an isomorphism on $j(X) = X^{\odot \odot}$,  we get from combining \eqref{eq:diffeqpsi} and \eqref{eq:jHq} that the coefficients must satisfy
\begin{equation*}
\begin{aligned}
    -\frac{\partial H_q^{00}(\tau,\xi^{(q)})}{\partial \tau} - D_\xi H_q^{00}(\tau,\xi^{(q)}) \tilde{M}_0 \xi &= p_q(\tau,\xi^{(q)}) - R_q^{00}(\tau,\xi^{(q)}), \\
    -\frac{\partial \tilde{H}_q^0(\tau,\xi^{(q)})}{\partial \tau} + \tilde{M}_0 \tilde{H}_q^0(\tau,\xi^{(q)}) - D_\xi \tilde{H}_q^0(\tau,\xi^{(q)})\tilde{M}_0 \xi &= P_q(\tau,\xi^{(q)}) - {R}_q^0(\tau,\xi^{(q)}).
\end{aligned}
\end{equation*}
These are precisely the equations obtained in \cite[Theorem III.7]{Iooss1999} and so the same proof can be followed to obtain the results on periodic normal forms of \cite[Theorem III.7]{Iooss1999}. In addition, it is proven in \cite[Theorem III.7]{Iooss1999} that $H_q^{00},\tilde{H}_q^0,p_q$ and $P_q$ are continuous, and so we conclude that $H,p$ and $P$ are $C^k$-smooth maps.
\end{proof}

Recall from \cref{thm:eigenfunctions} that the map $\tau \mapsto \gamma_\tau$ is $T$-periodic and $C^k$-smooth. Furthermore, from \eqref{eq:Floquetmap} in combination with \eqref{eq:Q0tilde} we have that $\tau \mapsto \tilde{Q}_0(\tau)$ is $T$-periodic and $C^k$-smooth. It also follows from previous theorem that $(\tau,\xi) \mapsto H(\tau,\xi)$ is $T$-periodic in the first component and $C^k$-smooth. Hence, $\mathcal{W}_{\loc}^{c}(\Gamma)$ can be written locally around $\Gamma$ as
\begin{equation} \label{eq:Wloc2}
    \mathcal{W}_{\loc}^{c}(\Gamma) = \{ \gamma_\tau  + \tilde{Q}_0(\tau)\xi + H(\tau,\xi) \in X : \tau \in \mathbb{R}  \mbox{ and } \xi \in \mathbb{R}^{n_0} \},
\end{equation}
and has exactly the same properties as the description of $\mathcal{W}_{\loc}^{c}(\Gamma)$ given in \eqref{eq:Wloc}. Hence, $\mathcal{W}_{\loc}^{c}(\Gamma)$ is the center manifold for \eqref{eq:DDEphi} around the periodic orbit $\Gamma$ whenever the Floquet multiplier $\lambda$ fulfills the requirements of \Cref{thm:normalformI}. 

Recall from \cite[Corollary 2]{Lentjes2023} that the center manifold $\mathcal{W}_{\loc}^{c}(\Gamma)$ is $T$-periodic, i.e. the map $\mathcal{C}$ from \eqref{eq:Wloc} is $T$-periodic in the first component, independent of the location and multiplicities of the Floquet multipliers.  Moreover, \eqref{eq:Wloc2} shows that $\mathcal{W}_{\loc}^{c}(\Gamma)$ has a $T$-periodic \emph{parametrization} since for any $\xi \in \mathbb{R}^{n_0}$ the map $\tau \mapsto \gamma_\tau  + \tilde{Q}_0(\tau)\xi + H(\tau,\xi)$ is $T$-periodic. A similar representation of $\mathcal{W}_{\loc}^{c}(\Gamma)$ as shown in \eqref{eq:Wloc2} was also obtained in \cite[Section III.2.1]{Iooss1999} for finite-dimensional ODEs.

Next, we consider the case when the trivial Floquet multiplier has an algebraic multiplicity larger than $1$ and there is no Floquet multiplier located at $-1$. Such a normal form theorem would be helpful to study the fold bifurcation. 

\begin{theorem}[Normal Form II] \label{thm:normalformII}
Assume that the algebraic multiplicity of the trivial Floquet multiplier is larger than one and that $-1$ is not a Floquet multiplier. Then there exist $C^k$-smooth maps $H : \mathbb{R} \times \mathbb{R}^{n_0} \to X, p : \mathbb{R} \times \mathbb{R}^{n_0} \to \mathbb{R}$ and $P:\mathbb{R} \times \mathbb{R}^{n_0} \to \mathbb{R}^{n_0}$ such that any solution $x_t$ of \eqref{eq:DDEphi} on $\mathcal{W}_{\loc}^{c}(\Gamma)$ may be represented as
\begin{equation*}
    x_{t}= \gamma_\tau + \tilde{Q}_0(\tau)\xi + H(\tau,\xi), \quad t \in I,
\end{equation*}
where the time dependence of the coordinates $(\tau,\xi)$ describing
the dynamics of \eqref{eq:DDEphi} on $\mathcal{W}_{\loc}^{c}(\Gamma)$ is defined by the normal form
\begin{equation*}
    \begin{cases}
    \dfrac{d \tau}{dt} = 1+ \xi_1 + p(\tau,\xi) + \mathcal{O}(|\xi|^{k+1}),\\
    \\[-8pt]
    \dfrac{d \xi}{d\tau} = \tilde{M}_0 \xi + P(\tau,\xi) + \mathcal{O}(|\xi|^{k+1}).
    \end{cases}
\end{equation*}
Here the functions $H,p$ and $P$ are $T$-periodic in $\tau$ and at least quadratic in $\xi$. The $\mathcal{O}$-terms are also $T$-periodic in $\tau$.  Moreover, $p$ and $P$ are polynomials in $\xi$ of degree less than or equal to $k$ satisfying
\begin{equation*}
    \frac{d}{d\tau}p(\tau,e^{-\tau \tilde{M}_0^{\star}} \xi) = 0, \quad \frac{d}{d\tau}\bigg( e^{\tau \tilde{M}_0^{\star}} P(\tau,e^{-\tau \tilde{M}_0^{\star}} \xi)\bigg) = 0,
\end{equation*}
for all $\tau \in \mathbb{R}$ and $\xi \in \mathbb{R}^{n_0}$.
\end{theorem}
Notice the appearance of the $\xi_1$-term in the normal form description. This is due to the fact that the $\star$ in \eqref{eq:M0} is now replaced by $1$ instead of $0$ compared to \Cref{thm:normalformI}.

\begin{proof}[Proof of \Cref{thm:normalformII}]
We proceed in the same way as the proof of \Cref{thm:normalformI} and start eventually by comparing the expansions of
\begin{equation*}
    \frac{d}{dt}j(\gamma_\tau + \tilde{Q}_0(\tau)\xi + H(\tau,\xi)) = j\bigg(\dot{\gamma}_\tau + \frac{\partial \tilde{Q}_0(\tau)}{\partial \tau}\xi + \frac{\partial H(\tau,\xi)}{\partial \tau} +\bigg(\tilde{Q}_0(\tau) + D_\xi H(\tau,\xi) \bigg) \frac{d\xi}{d\tau} \bigg) \frac{d\tau}{dt}
\end{equation*}
and
\begin{equation*}
    A_0^{\odot \star}j(\gamma_\tau + \tilde{Q}_0(\tau)\xi + H(\tau,\xi)) + G(\gamma_\tau + \tilde{Q}_0(\tau)\xi + H(\tau,\xi))
\end{equation*}
by substituting
\begin{equation*}
    \frac{d \tau}{dt} = 1+ \xi_1 + p(\tau,\xi) + \mathcal{O}(|\xi|^{k+1}) \ \mbox{ and } \frac{d \xi}{d\tau} = \Tilde{M}_0 \xi + P(\tau,\xi) + \mathcal{O}(|\xi|^{k+1}).
\end{equation*}
We copy the same notation from the proof of \Cref{thm:normalformI} and use the expansions of $H,p$ and $P$ together with \eqref{eq:Gexpansion} where now $\varphi$ must be substituted by $\tilde{Q}_0(\tau)\xi + H(\tau,\xi)$. Eventually,
\begin{align*}
    &j \bigg[ \dot{\gamma}_\tau + \frac{\partial \tilde{Q}_0(\tau)}{\partial \tau} \xi + \sum_{q=2}^k \frac{\partial H_q(\tau,\xi^{(q)}))}{\partial \tau} + \bigg(\tilde{Q}_0(\tau) + \sum_{q=2}^k D_\xi H_q(\tau,\xi^{(q)}) \bigg) \bigg(\tilde{M}_0 + \sum_{q=2}^{k} P_q(\tau,\xi^{(q)}) \bigg) \bigg]\\
    & \bigg(1+ \xi_1+ \sum_{q=2}^k p_q(\tau,\xi^{(q)}) \bigg) + \mathcal{O}(|\xi|^{k+1})\\
   &= A_0^{\odot \star}j(\gamma_\tau) + G(\gamma_\tau) + A^{\odot \star}(\tau)j\bigg(\tilde{Q}_0(\tau)\xi + \sum_{q=2}^k H_q(\tau,\xi^{(q)}) \bigg) \\
   &+ \sum_{q=2}^k G_q \bigg(\tau,\bigg[\tilde{Q}_0(\tau)\xi + \sum_{p=2}^k H_q(\tau,\xi^{(p)}) \bigg]^{(q)} \bigg) + \mathcal{O}(|\xi|^{k+1}) .
\end{align*}
Let us now compare the terms in powers of $\xi$ on both side of the equation. The $\xi^0$-terms give us
\begin{equation*}
   \frac{d}{d\tau} j({\gamma}_\tau) = A_0^{\odot \star} j(\gamma_\tau) + G(\gamma_\tau),
\end{equation*}
which means that $\tau \mapsto \gamma_\tau$ is a solution \eqref{eq:invariance equation}. The $\xi^1$-terms tell us
\begin{equation} \label{eq:linear part CMT2}
    \bigg(-\frac{\partial}{\partial \tau} + A^{\odot \star}(\tau) \bigg) j(\tilde{Q}_0(\tau)\xi) = j((\tilde{Q}_0(\tau)\tilde{M}_0 + \gamma_\tau \Pi_1) \xi),
\end{equation}
which is exactly the result established in \eqref{eq:Floquetmap}, but now for all Floquet multipliers on the unit circle and characterizes the linear part. Here $\Pi_1 : \mathbb{R}^{n_0} \to \mathbb{R}$ is the projection on the first component, defined by $\Pi_1 \xi := \xi_1$, where $\xi = (\xi_1,\dots,\xi_{n_0})$. After collecting the $\xi^{(2)}$-terms, we get 
\begin{align*}
    &\bigg( -\frac{\partial}{\partial \tau} + A^{\odot \star}(\tau) \bigg) j(H_2(\tau,\xi^{(2)})) \\
    &= j(D_\xi H_2(\tau,\xi^{(2)}) \tilde{M}_0 \xi + p_2(\tau,\xi^{(2)})\dot{\gamma}_\tau + \tilde{Q}_0(\tau) P_2(\tau,\xi^{(2)})) - R_2(\tau,\xi^{(2)}),
\end{align*}
where $R_2(\tau,\xi^{(2)}) = G_2(\tau,(\tilde{Q}_0(\tau) \xi)^2) - \xi_1(\frac{\partial \tilde{Q}_0(\tau)}{\partial \tau} \xi + \tilde{Q}_0(\tau)\tilde{M}_0 \xi)$. Finally, after collecting the $\xi^{(q)}$-terms for $q = 3,\dots,k$, we get
\begin{align} \label{eq:diffeqpsi2}
\begin{split}
        &\bigg( -\frac{\partial}{\partial \tau} + A^{\odot \star}(\tau) \bigg) j(H_q(\tau,\xi^{(q)})) \\
        &= j(D_\xi H_q(\tau,\xi^{(q)}) \tilde{M}_0 \xi + p_q(\tau,\xi^{(q)})\dot{\gamma}_\tau + \tilde{Q}_0(\tau) P_q(\tau,\xi^{(q)})) - R_q(\tau,\xi^{(q)}),
\end{split}
\end{align}
where $R_q(\tau,\xi^{(q)})$ depends on $G_{q'}(\tau,\xi^{(q')})$ for $2\leq q' \leq q$ and $j(H_{q'}(\tau,\xi^{(q')}))$, $j(p_{q'}(\tau,\xi^{(q')})\dot{\gamma}_\tau)$ and $j( \tilde{Q}_0(\tau)P_{q'}(\tau,\xi^{(q')}))$ for $q' = 2,\dots,q-1$.

We want to project \eqref{eq:diffeqpsi2} onto the spaces $\mathbb{R}\dot{\gamma}_\tau, \tilde{X}_0(\tau)$ and $ {X}_{\pm}(\tau)$ to obtain a representation of $H_q$ separately on each individual space. We decompose for any $q=2,\dots,k$ the functions $H_q$ and $R_q$ as
\begin{align*}
    H_q(\tau,\xi^{(q)}) &= H_q^{00}(\tau,\xi^{(q)})\dot{\gamma}_\tau + \tilde{Q}_0(\tau)\tilde{H}_q^0 (\tau,\xi^{(q)}) + H_q^{-}(\tau,\xi^{(q)}) + H_q^{+}(\tau,\xi^{(q)}) \\
    R_q(\tau,\xi^{(q)}) &= R_q^{00}(\tau,\xi^{(q)})\dot{\gamma}_\tau + \tilde{Q}_0(\tau)\tilde{R}_q^0 (\tau,\xi^{(q)}) + R_q^{-}(\tau,\xi^{(q)}) + R_q^{+}(\tau,\xi^{(q)}),
\end{align*}
where $H_q^{\pm}(\tau,\xi^{(q)}) = P_{\pm} (\tau) H_q (\tau,\xi^{(q)}) \in X_{\pm}(\tau)$ and $R_q^{\pm}(\tau,\xi^{(q)}) = P_{\pm}^{\odot \star} (\tau) R_q (\tau,\xi^{(q)}) \in X_{\pm}^{\odot \star}(\tau)$ with coordinates $H_q^{00}(\tau,\xi^{(q)}), R_q^{00}(\tau,\xi^{(q)}) \in \mathbb{R}$ and $\tilde{H}_q^0 (\tau,\xi^{(q)}), \tilde{R}_q^0 (\tau,\xi^{(q)}) \in \mathbb{R}^{n_0}$ for all $\tau \in \mathbb{R}$ and $\xi \in \mathbb{R}^{n_0}$.

Carrying out the calculations in the same way as the proof of \Cref{thm:normalformII} and recalling that $(-\frac{d}{d\tau} + A^{\odot \star}(\tau))j\dot{\gamma}_\tau = 0$ together with \eqref{eq:linear part CMT2}, we obtain
\begin{align*}
\bigg( -\frac{\partial}{\partial \tau} + A^{\odot \star}(\tau) \bigg) j(H_q(\tau,\xi^{(q)})) &=j\bigg(- \frac{\partial H_q^{00}(\tau,\xi^{(q)})}{\partial \tau} \dot{\gamma}_\tau + \Pi_1 \tilde{H}_q^0(\tau,\xi^{(q)})\dot{\gamma}_\tau \bigg) \\
&
+ j\bigg(\tilde{Q}_0(\tau) \bigg(- \frac{\partial}{\partial \tau} + \tilde{M}_0 \bigg) \tilde{H}_q^0(\tau,\xi^{(q)}) \bigg) \nonumber \\
&
+ \bigg( -\frac{\partial}{\partial \tau} + A^{\odot \star}(\tau) \bigg) j(H_q^{-}(\tau,\xi^{(q)}) + H_q^{+}(\tau,\xi^{(q)}))
\end{align*}
and this must be equal to the right-hand side of \eqref{eq:diffeqpsi2}. To obtain the coefficients, we project onto the spaces $\mathbb{R}j\dot{\gamma}_\tau, j(\tilde{X}_0(\tau)), j(X_{+}(\tau))$ and $X_{-}^{\odot \star}(\tau)$. This yields the equations
\begin{align*}
    -\frac{\partial H_q^{00}(\tau,\xi^{(q)})}{\partial \tau} D_\xi H_q^{00}(\tau,\xi^{(q)}) \tilde{M}_0 \xi \hspace{-2pt} &= \hspace{-2pt} p_q(\tau,\xi^{(q)}) - \Pi_1 \tilde{H}_q^0(\tau,\xi^{(q)}) - R_q^{00}(\tau,\xi^{(q)}) \\
    -\frac{\partial \tilde{H}_q^0(\tau,\xi^{(q)})}{\partial \tau} + \tilde{M}_0 \tilde{H}_q^0(\tau,\xi^{(q)}) - D_\xi \tilde{H}_q^0(\tau,\xi^{(q)})\tilde{M}_0 \xi &= P_q(\tau,\xi^{(q)}) - \tilde{R}_q^0(\tau,\xi^{(q)})\\
        \bigg( -\frac{\partial}{\partial \tau} + A^{\odot \star}(\tau) \bigg) j(H_q^{\pm}(\tau,\xi^{(q)})) &= j(D_\xi H_q^{\pm}(\tau,\xi^{(q)})\tilde{M}_0 \xi) - R_q^{\pm}(\tau,\xi^{(q)}).
\end{align*}
We see that the equations for the $X_{\pm}^{\odot \star}(\tau)$-component are the same as in the proof of \Cref{thm:normalformI}. Hence, we obtain $H_q^{\pm}$ as in \eqref{eq:psi_q^+} and \eqref{eq:psi_q^-} respectively. To solve the remaining part of this hierarchy of equations, notice these equations are solvable in exactly the same way as the proof of \Cref{thm:normalformI} and the proposed normal forms follow. One should make the observation that $\tilde{H}_q^{0}$ has to be computed before $H_q^{00}$.
\end{proof}

Under these assumptions on the Floquet multipliers, we have proven that $\mathcal{W}_{\loc}^{c}(\Gamma)$ can also be parametrized as in \eqref{eq:Wloc2}. Moreover, $\mathcal{W}_{\loc}^{c}(\Gamma)$ has a $T$-periodic parametrization.

The last normal form theorem is more involved because we have to deal with the Floquet multiplier $-1$ that induces $T$-antiperiodic (generalized) eigenfunctions due to \Cref{prop:eigenfunctions2}. Introduce the decomposition
\begin{equation*}
    \tilde{X}_0(\tau) = \tilde{X}'_0(\tau) \oplus \tilde{X}''_0(\tau),
\end{equation*}
where $\tilde{X}'_0(\tau)$ is spanned by $T$-periodic maps $\varphi_1(\tau),\dots,\varphi_{n_0'}(\tau)$ and $\tilde{X}''_0(\tau)$ is spanned by $2T$-periodic maps $\varphi_{n_0' + 1}(\tau),\dots,\varphi_{n_0' + n_0''}(\tau)$ and $n_0' + n_0'' = n_0$, corresponding to all (generalized) eigenfunctions belonging to the Floquet multiplier $-1$. If we define the symmetry $\tilde{S}_0 : \mathbb{R}^{n_0} \to \mathbb{R}^{n_0}$ by
\begin{equation*}
    (\xi',\xi'') = \xi \mapsto \tilde{S}_0\xi = (\xi',-\xi''),
\end{equation*}
then we have the following normal form theorem. Such a normal form theorem would be helpful to study the period-doubling bifurcation. 

\begin{theorem}[Normal form III] \label{thm:normalformIII} Assume that $-1$ is a Floquet multiplier. Then the results of \Cref{thm:normalformI} or \Cref{thm:normalformII}, depending on the location and algebraic multiplicity of the other Floquet multipliers on the unit circle hold with the following modification: the maps $H,p$ and $P$, and $\mathcal{O}$-terms are $2T$-periodic in $\tau$ and satisfy additionally
\begin{equation*}
    H(\tau+T,\xi) = H(\tau,\tilde{S}_0 \xi), \quad p(\tau+T,\xi) = p(\tau,\tilde{S}_0 \xi), \quad P(\tau+T,\xi) = \tilde{S}_0  P(\tau,\xi),
\end{equation*}
for all $\tau \in \mathbb{R}$ and $\xi \in \mathbb{R}^{n_0}$.
\end{theorem}
\begin{proof}
The proof is analogous to that of \Cref{thm:normalformI} or \Cref{thm:normalformII} but in a $2T$-periodic setting. Hence, we obtain the results from \Cref{thm:normalformI} or \Cref{thm:normalformII}, depending on the location and algebraic multiplicity of the Floquet multipliers on the unit circle where now the maps $H,p$ and $P$ being $2T$-periodic in $\tau$. It remains to show the additional symmetries on the maps $H,p$ and $P$. Because the structure of the parts in the normal form are similar to that of the ODE case, treated in \cite[Theorem III.13]{Iooss1999}, this part will be omitted since the proof is identical by making the substitution of $\tau \mapsto \gamma(\tau)$ towards $\tau \mapsto \gamma_\tau$.
\end{proof}
Under the assumptions of \cref{thm:normalformIII}, $\mathcal{W}_{\loc}^{c}(\Gamma)$ can also be parametrized as performed in \eqref{eq:Wloc2}. However, compared to \cref{thm:normalformII} and \cref{thm:normalformIII}, $\mathcal{W}_{\loc}^{c}(\Gamma)$ now has a $2T$-periodic parametrization as for any $\xi \in \mathbb{R}^{n_0}$ the map $\tau \mapsto \gamma_\tau + \tilde{Q}_0(\tau)\xi + H(\tau,\xi)$ is $2T$-periodic due to the presence of a negative Floquet multiplier on the unit circle. 

\section{Conclusion and outlook} \label{sec:conclusions}
We have proven that the periodic normal forms from Iooss (and Adelmeyer) \cite{Iooss1988,Iooss1999} for bifurcations of limit cycles in finite-dimensional ODEs suit naturally in the framework of infinite-dimensional classical DDEs. This task has been accomplished by proving two principal results: the existence of a dual periodic smooth basis for the center eigenspace (\Cref{thm:eigenfunctions} and \Cref{thm:adjoint eigenfunctions}) and the existence of a special coordinate system on the periodic center manifold (\Cref{thm:normalformI}, \Cref{thm:normalformII} and \Cref{thm:normalformIII}) in which the local dynamics near the nonhyperbolic cycle can be described. A paper providing explicit computational formulas for the critical normal form coefficients for all codimension one bifurcations of limit cycles in classical DDEs, along the lines of the periodic normalization method \cite{Kuznetsov2005,Witte2014,Witte2013}, is in preparation.

In the present paper, we restricted ourselves to the setting of classical DDEs. However, our proof on the existence of a periodic smooth finite-dimensional center manifold near a nonhyperbolic cycle in \cite[Theorem 4]{Lentjes2023} is given in the general context of dual perturbation theory (sun-star calculus). As a consequence, the periodic center manifold theorem also holds for a much broader class of delay equations, as for example renewal equations \cite{Breda2020, Diekmann2008,Diekmann1995}. The natural question arises if the results from this paper can also be generalized towards the general context of sun-star calculus. Regarding the linear part, one can not expect in general that the (adjoint) (generalized) eigenfunctions are $C^{k+1}$-smooth as this was only a consequence of applying the method of steps to periodic linear DDEs (\cref{lemma:Dinvariant} and \cref{lemma:Dsuninvariant}). This loss of smoothness of the (generalized) eigenfunctions also causes a problem in the proof of the normal forms theorems (\cref{subsec: normal form}) as the parametrization of the periodic center manifold \eqref{eq:x_t manifold} is then not $C^{k}$-smooth anymore. However, if one assumes that the strongly continuous evolutionary systems $U$ and $U^\odot$ are \emph{eventually domain invariant} (there exists a $h \geq 0$ such that $U(t,s)$ maps $X$ into $\mathcal{D}(A(s))$ for $t \geq s + h$, and $U^\odot(s,t)$ maps $X^\odot$ into $\mathcal{D}(A^\odot(s))$ for $s \leq t - h$) and one still allows differentiability in the strong sense, then we believe that the generators of the strongly continuous evolution semigroups $\mathcal{U}$ and $\mathcal{U}^\odot$ are still given by the closure of the densely defined linear operators $\mathcal{A}$ and $\mathcal{A}^\odot$, respectively. It would then be interesting to see how the theory of evolution semigroups fits into the sun-star calculus framework. This would enable us to characterize a lot of relations between $\mathcal{U}^\star$ and $\mathcal{A}^\star$ defined on (a subspace of) $C_T(\mathbb{R},X^\star)$, and $\mathcal{U}^{\odot \star}$ and $\mathcal{A}^{\odot \star}$ defined on (a subspace of) $C_T(\mathbb{R},X^{\odot \star})$. Afterwards, it would be interesting to study the duality and spectral relations, as performed in \cref{subsec: duality relations} for classical DDEs, but in the general context of dual perturbation theory. When one does not allow differentiability in the strong sense, we can not introduce $\mathcal{A}$ on a (dense) subspace of $C_T(\mathbb{R},X)$. Therefore, we must study this operator on a (dense) subspace of, for example, the Bochner space $L_T^{p}(\mathbb{R},X)$ for some $1 \leq p \leq \infty$. This has been done recently in \cite{Arendt2009} since evolution semigroups and their (restricted) generators are related with $L^p$-maximal regularity of periodic boundary value problems, see \cite{Lunardi1995,Arendt2002,Dore1993,Chicone1999} for more information. There are some advantages regarding duality when one would work in the $L_T^{2}(\mathbb{R},X)$ setting, since then the pairing $\langle \cdot,\cdot \rangle_T$ from \eqref{eq:pairingT} is precisely the (canonical) dual pairing. However, it is not clear to the authors if the normal forms from \cref{subsec: normal form} in this setting would hold, as one loses the smoothness of the Jordan chains and the center manifold parametrization \eqref{eq:x_t manifold}.

If one is interested in bifurcations of limit cycles for systems consisting of infinite delay \cite{Diekmann2012,Diekmann2007,Hino1991,Liu2020} or abstract DDEs \cite{Janssens2019,Janssens2020,Spek2020,Spek2022,Gils2012,Spek2022a} that describe for example neural fields, it is known that $\odot$-reflexivity is in general lost \cite[Theorem 12]{Spek2020}, and therefore the center manifold theorem for nonhyperbolic cycles from \cite{Lentjes2023} does not directly apply. However, we believe that this technical difficulty can be resolved by employing similar techniques as in \cite{Janssens2020}. We are convinced that these techniques can also be used to prove the existence of the periodic normal forms in the setting of abstract DDEs and systems consisting of infinite delay because the proof on the existence of the periodic normal forms are written in a rather general setting of sun-star calculus.

\section*{Acknowledgements}
The authors express their sincere gratitude to Prof. Odo Diekmann (Utrecht University) for his invaluable recommendation to delve into the theory of evolution semigroups. The authors are thankful to Prof. Peter De Maesschalck (Hasselt University), Prof. Stephan van Gils (University of Twente) and Stein Meereboer (Radboud University) for helpful discussions and suggestions.

\appendix

\section{Variation-of-constants formula for the adjoint problem} \label{variation-of-constants adjoint}
In this section of the appendix, we will prove that solutions of an inhomogeneous perturbed abstract ordinary differential equation are precisely given by an abstract integral equation. This result is important in the proof of \Cref{thm:adjoint eigenfunctions}.
Let $J \subseteq \mathbb{R}$ be an interval and suppose that $(s,t) \in \Omega_J^\star$. Consider an inhomogeneous perturbation $f : J \to X^\star $ on the generator $A^\star(s)$ to the adjoint problem \cite[Equation (5.8)]{Clement1988}. This yields the initial value problem
\begin{equation} \label{eq:adjointODEf}
    \begin{dcases}
    d^\star u(s) = -A^\star(s)u(s) + f(s), \quad &s\leq t, \\
   u(t) = \psi, \quad &\psi \in X^\odot,
     \end{dcases}
\end{equation}
which suggests the variation-of-constants formula 
\begin{equation} \label{eq:adjointAIEf}
    u(s) = U^\odot(s,t)\psi + \int_t^s  U^\star(s,\tau)f(\tau)d\tau, \quad \psi \in  X^\odot ,
\end{equation}
for $t \leq s$, where the integral must be interpreted as a weak$^\star$ integral. This suggestion, with the additional assumptions on $f$, will be verified in this section. Let us first turn our attention towards the weak$^\star$ integral appearing in \eqref{eq:adjointAIEf}.

\begin{lemma} \label{lemma:wkstarintegral}
Let $g : J \to X^\star $ be a continuous function and denote the set $\{(s,r,t) \in J^3 : s \leq r \leq t \}$ by $\Theta_J^\star$. Then the map $v(\cdot,\cdot,\cdot,g) : \Theta_J^\star \to X^\star$ defined as the weak$^\star$ integral
\begin{equation*}
     v(s,r,t,g) := \int_t^s U^\star(r,\tau)g(\tau) d\tau,  \quad \forall (s,r,t)\in  \Theta_J^\star,
\end{equation*}
is continuous and takes values in $X^\odot$.
\end{lemma}
\begin{proof}
The statement of the theorem is a dual version of the first part of \cite[Lemma 1]{Lentjes2023}. The proof is along the same lines as the proof in \cite[Lemma 1]{Lentjes2023}, where one just has to work with the backward evolutionary system $U^\star$. See also \cite[Lemma 3.1]{Clement1989a} for a semigroup analogue of this lemma.
\end{proof}
As we have proven that the weak$^\star$ integral, appearing in \eqref{eq:adjointAIEf} is well-defined, we can turn our attention towards the verification of the variation-of-constants formula. The proof is inspired from \cite[Proposition 13]{Lentjes2023} and \cite[Proposition 21]{Janssens2020}.

\begin{proposition} \label{prop:solutionadjoint}
Let $f : J \to X^\star$ be a continuous function. If $u$ is a solution of \eqref{eq:adjointODEf} on $J$ then $u$ is given by \eqref{eq:adjointAIEf}.
\end{proposition}
\begin{proof}
Let $(s,t) \in \Omega_J^\star$ with $t > s$ be arbitrary.  Define the function $w : [s,t] \to X^\star$ as $w(\tau) := U^\star(s,\tau)u(\tau)$ for all $\tau \in [s,t]$. We claim that $w$ is weak$^\star$ differentiable with weak$^\star$ derivative
\begin{equation} \label{eq:dstarw}
     d^\star w(\tau) = U^\star(s,\tau)d^\star u(\tau) + U^\star(s,\tau)A^\star(\tau)u(\tau).
\end{equation}
To show this claim, let $\tau \in [s,t]$ and $x \in X$ be given. For any $h \in \mathbb{R}$ such that $\tau + h \in [s,t]$ we have that
\begin{align*}
     \langle w(\tau + h) - w(\tau),x \rangle &= \langle U^\star(s,\tau + h)u(\tau + h) - U^\star(s,\tau)u(\tau), x \rangle \\
     &= \langle U^\star(s,\tau + h)[u(\tau + h) - u(\tau)],x \rangle + \langle [U^\star(s,\tau + h)-U^\star(s,\tau)]u(\tau), x \rangle \\
     &= \langle u(\tau + h) - u(\tau),U(\tau + h,s)x \rangle +  \langle [U^\star(s,\tau + h)-U^\star(s,\tau)]u(\tau), x \rangle
\end{align*}
Because $U$ is a strongly continuous forward evolutionary system, we have that $U(\tau+h,s)x \to U(\tau,s)x$ in norm as $h \downarrow 0$. The definition of the weak$^\star$ derivative implies
\begin{equation*}
     \frac{1}{h}(u(\tau + h) -u(\tau)) \to d^\star u(\tau), \quad \mbox{ weakly}^\star \ h\downarrow 0,
\end{equation*}
if we can show that the difference quotients remains bounded in the limit. Since $u$ is a solution to \eqref{eq:adjointODEf}, we know that $u$ is weak$^\star$ continuously differentiable and so locally Lipschitz continuous by \cite[Remark 16]{Janssens2020}. Because $[s,t]$ is compact, $u$ is Lipschitz continuous on $[s,t]$ and so the difference quotient remains bounded in the limit. Combining these two facts, we get
\begin{equation*}
     \frac{1}{h}\langle u(\tau+h) - u(\tau),U(\tau+h,s)x  \rangle \to \langle d^\star u(\tau), U(s,\tau)x \rangle, \quad  h \downarrow 0.
\end{equation*}
Furthermore, since $u(\tau) \in \mathcal{D}(A^\star(\tau)) = \mathcal{D}(A_0^\star)$ (recall \cref{lemma:adjoint D(Astar)}), it follows from \cite[Theorem 5.7]{Clement1988} that 
\begin{equation*}
     \frac{1}{h} \langle [U^\star(s,\tau+h) - U^\star(s,\tau)]u(\tau), x \rangle \to \langle U^\star(s,\tau)A^\star(\tau)u(\tau),x \rangle, \quad  h \downarrow 0.
\end{equation*}
Hence, we get
\begin{equation*}
     \frac{1}{h} \langle w(\tau + h) - w(\tau),x \rangle \to \langle U^\star(s,\tau)d^\star u(\tau) + U^\star(s,\tau)A^\star(\tau)u(\tau), x \rangle, \quad  h \downarrow 0,
\end{equation*}
which proves \eqref{eq:dstarw}. Substituting the differential equation from \eqref{eq:adjointODEf} into \eqref{eq:dstarw} yields
\begin{equation*}
     d^\star w(\tau) = U^\star(s,\tau)f(\tau)
\end{equation*}
and so $d^\star w$ is weak$^\star$ continuous since $f$ was assumed to be norm continuous. For every $x \in X$ we have
\begin{align*}
     \langle u(s) - U^{\star}(s,t)u(t),x \rangle &= \langle w(s), x \rangle - \langle w(s), x \rangle = \int_t^s \langle d^\star w(\tau), x \rangle d\tau = \langle \int_t^s U^{ \star}(s,\tau)f(\tau) d\tau, x \rangle .
\end{align*}
Since $x$ and $(s,t) \in \Omega_J^\star$ were arbitrary, we conclude that
\begin{equation*}
     u(s) - U^{\star}(s,t)u(t) = \int_t^s U^{  \star}(s,\tau)f(\tau) d\tau,
\end{equation*}
or equivalently
\begin{equation*}
     u(s) = U^{\star}(s,t)\psi + \int_t^s U^{  \star}(s,\tau)f(\tau) d\tau,
\end{equation*}
since $u(t) = \psi$ by assumption. The continuity and range of $f$ ensures from \Cref{lemma:wkstarintegral} that the weak$^\star$ integral takes values in $X^\odot$. Since $\psi \in X^\odot$, we have that
\begin{equation*}
     u(s) = U^{\odot}(s,t)\psi + \int_t^s U^{ \star}(s,\tau)f(\tau) d\tau,
\end{equation*}
which completes the proof.
\end{proof}

\bibliographystyle{siamplain}
\bibliography{references}

\end{sloppypar}
\end{document}